\documentclass[A4j,11pt]{article}
\usepackage{amsfonts,amsmath,amscd}
\usepackage{graphics}
\usepackage{amssymb}
\usepackage[all]{xy}
\begin{document}

\title{{\bf Collapse of the mean curvature flow\\
for equifocal submanifolds}}
\author{{\bf Naoyuki Koike}}
\date{}
\maketitle

\begin{abstract}
In this paper, we investigate the mean curvature flows having an equifocal 
submanifold in a symmetric space of compact type and its focal submanifolds 
as initial data.  
It is known that an equifocal submanifold of codimension greater than one 
in an irreducible symmetric space of compact type occurs as a principal orbit 
of a Hermann action.  However, we investigate the flows 
conceptionally without use of this fact.  The investigation 
is performed by investigating the mean curvature flows having the lifts of 
the submanifolds to an (infinite dimensional separable) Hilbert space through 
a Riemannian submersion as initial data .  
\end{abstract}

\vspace{0.5truecm}

%\vspace{0.2truecm}

%\noindent
%$\overline{\qquad\qquad\qquad\qquad\qquad\qquad}$

%\vspace{0.1truecm}

%{\rm \footnotesize{2000 Mathematics Subject Classification. Primary 53C40, 
%53C35.}}

%{\rm \footnotesize{Key words and phrases. equifocal submanifold, polar 
%action.}}

\section{Introduction}
Let $f_t$'s ($t\in[0,T)$) be a one-parameter $C^{\infty}$-family of immersions 
of a manifold $M$ into a Riemannian manifold $N$, where $T$ is a positive 
constant or $T=\infty$.  Define a map 
$\widetilde f:M\times[0,T)\to N$ by $\widetilde f(x,t)=f_t(x)$ 
($(x,t)\in M\times[0,T)$).  If, for each $t\in[0,T)$,  
$\widetilde f_{\ast}((\frac{\partial}{\partial t})_{(x,t)})$ is the mean 
curvature vector of $f_t:M\hookrightarrow N$, then $f_t$'s ($t\in[0,T)$) is 
called a mean curvature flow.  
In particular, if $f_t$'s are embeddings, then we call $M_t:=f_t(M)$'s 
$(0\in[0,T))$ rather than $f_t$'s $(0\in[0,T))$ a mean curvature flow.  
Liu-Terng [LT] investigated the mean curvature flow having isoparametric 
submanifolds (or their focal submanifolds) in a Euclidean space 
as initial data and obtained the following facts.  

\vspace{0.5truecm}

\noindent
{\bf Fact 1([LT]).} {\sl Let $M$ be a compact isoparametric submanifold in a 
Euclidean space.  
Then the following statements ${\rm (i)}$ and ${\rm (ii)}$ hold:

${\rm (i)}$ The mean curvature flow $M_t$ having $M$ as initial data collapses 
to a focal submanifold of $M$ in finite time.  If a focal map of $M$ onto $F$ 
is spherical, then the mean curvature flow $M_t$ has type I singularity, 
that is, $\lim\limits_{t\to T-0}
{\rm max}_{v\in S^{\perp}M_t}\vert\vert A^t_v\vert\vert^2_{\infty}(T-t)\,<\,
\infty$, where $A^t_v$ is the shape operator of $M_t$ for $v$, 
$\vert\vert A^t_v\vert\vert_{\infty}$ is the sup norm of $A^t_v$ and 
$S^{\perp}M_t$ is the unit normal bundle of $M_t$.

{\rm(ii)} For any focal submanifold $F$ of $M$, there exists a parallel 
submanifold $M'$ of $M$ such that the mean curvature flow having $M'$ 
as initial data collapses to $F$ in finite time.}

\vspace{0.5truecm}

\noindent
{\bf Fact 2([LT]).} {\sl Let $M$ be as in Fact 1, $C$ be the Weyl domain of 
$M$ at $x_0\,(\in M)$ and $\sigma$ be a 
stratum of dimension greater than zero of $\partial C$.  
Then the following statements ${\rm (i)}$ and ${\rm (ii)}$ hold:

{\rm(i)} For any focal submanifold $F$ (of $M$) through $\sigma$, the maen 
curvature flow $F_t$ having $F$ as initial data collapses to a focal 
submanifold $F'$ (of $M$) through $\partial{\sigma}$ in finite time.  
If the fibration of $F$ onto $F'$ is spherical, then the mean curvature flow 
$F_t$ has type I singularity.

(ii) For any focal submanifold $F$ (of $M$) through $\partial\sigma$, there 
exists a focal submanifold $F'$ (of $M$) through $\sigma$ such that the mean 
curvature flow $F'_t$ having $F'$ as initial data collapses to $F$ 
in finite time.}

\vspace{0.5truecm}

As a generalized notion of compact isoparametric 
hypersurfaces in a sphere and a hyperbolic space, and a compact isoparametric 
submanifolds in a Euclidean space, Terng-Thorbergsson [TT] defined 
the notion of an equifocal submanifold in a symmetric space as a compact 
submanifold $M$ satisfying the following three conditions:

\vspace{0.2truecm}

(i) the normal holonomy group of $M$ is trivial,

(ii) $M$ has a flat section, that is, for each $x\in M$, 
$\Sigma_x:=\exp^{\perp}(T^{\perp}_xM)$ is totally geodesic and the induced 
metric on $\Sigma_x$ is flat, where $T^{\perp}_xM$ is the normal space of $M$ 
at $x$ and $\exp^{\perp}$ is the normal exponential map of $M$.  

(iii) for each parallel normal vector field $v$ of $M$, the focal radii of $M$ 
along the normal geodesic $\gamma_{v_x}$ (with $\gamma'_{v_x}(0)=v_x$) are 
independent of the choice of $x\in M$, where $\gamma'_{v_x}(0)$ is the 
velocity vector of $\gamma_{v_x}$ at $0$.  

\vspace{0.2truecm}

\noindent
On the other hand, Heintze-Liu-Olmos [HLO] defined the notion of an 
isoparametric submanifold with flat section in a general Riemannian manifold 
as a submanifold $M$ satisfying the above condition (i) and the 
following conditions (ii$'$) and (iii$'$):

\vspace{0.2truecm}

(ii$'$) for each $x\in M$, there exists a neighborhood $U_x$ of 
the zero vector (of $T^{\perp}_xM$) in $T^{\perp}_xM$ such that 
$\Sigma_x:=\exp^{\perp}(U_x)$ is totally geodesic and the induced metric on 
$\Sigma_x$ is flat, 

\vspace{0.2truecm}

(iii$'$) sufficiently close parallel submanifolds of $M$ are CMC with 
respect to the radial direction.  

\vspace{0.2truecm}

\noindent
In the case where the ambient space is a symmetric space $G/K$ of compact 
type, they showed that the notion of an isoparametric submanifold with flat 
section coincides with that of an equifocal submanifold.  
The proof was performed by investigating its lift to $H^0([0,1],\mathfrak g)$ 
through a Riemannian submersion $\pi\circ\phi$, where $\pi$ is the natural 
projection of $G$ onto $G/K$ and $\phi$ is the parallel transport map for 
$G$ (which is a Riemannian submersion of $H^0([0,1],\mathfrak g)$ onto $G$ 
($\mathfrak g:$the Lie algebra of $G$)).  
Let $M$ be an equifocal submanifold in $G/K$ and $v$ be a parallel normal 
vector field of $M$.  The end-point map $\eta_v(:M\mapsto G/K)$ for $v$ 
is defined by $\eta_v(x)=\exp^{\perp}(v_x)$ ($x\in M$).  
Set $M_v:=\eta_v(M)$.  We call $M_v$ 
a parallel submanifold of $M$ when ${\rm dim}\,M_v={\rm dim}\,M$ and a focal 
submanifold of $M$ when ${\rm dim}\,M_v<{\rm dim}\,M$.  The parallel 
submanifolds of $M$ are equifocal.  Let $f:M\times[0,T)\to G/K$ be 
the mean curvature flow having $M$ as initial data.  
Then, it is shown that, for each 
$t\in[0,T)$, $f_t:M\hookrightarrow G/K$ is a parallel submanifold of $M$ 
and hence it is equifocal (see Lemma 3.1).  Fix $x_0\in M$.  Let 
$\widetilde C\,(\subset T^{\perp}_{x_0}M)$ be the fundamental domain 
containing the zero vector (of $T^{\perp}_{x_0}M$) of the Coxeter group 
(which acts on $T^{\perp}_{x_0}M$) of $M$ at $x_0$ 
and set $C:=\exp^{\perp}(\widetilde C)$, where we note that 
$\exp^{\perp}\vert_{\widetilde C}$ is a diffeomorphism onto $C$.  Without loss 
of generality, we may assume that $G$ is simply connected.  
Set $\widetilde M:=(\pi\circ\phi)^{-1}(M)$, which is an isoparametric 
submanifold in $H^0([0,1],\mathfrak g)$.  
Fix $u_0\in(\pi\circ\phi)^{-1}(x_0)$.  
The normal space $T^{\perp}_{x_0}M$ is identified with the normal space 
$T^{\perp}_{u_0}\widetilde M$ of $\widetilde M$ at $u_0$ through 
$(\pi\circ\phi)_{\ast u_0}$.  Each parallel submanifold of $M$ intersects 
with $C$ at the only point and each focal submanifold of $M$ intersects with 
$\partial C$ at the only point, where $\partial C$ is the boundary of $C$.  
Hence, for the mean curvature flow $f:M\times [0,T)\to G/K$ having $M$ as 
initial data, each $M_t(:=f_t(M))$ intersects with $C$ at the only point.  
Denote by $x(t)$ this intersection point and define $\xi:[0,T)\to\widetilde C\,
(\subset T^{\perp}_{x_0}M=T^{\perp}_{u_0}\widetilde M)$ by 
$\exp^{\perp}(\xi(t))=x(t)$ ($t\in[0,T)$).  Set 
$\widetilde M_t:=(\pi\circ\phi)^{-1}(M_t)$ ($t\in[0,T)$).  It is shown that 
$\widetilde M_t$ ($t\in[0,T)$) is the mean curvature flow having 
$\widetilde M$ as initial data because the mean curvature vector of 
$\widetilde M_t$ is the horizontal lift of that of $M_t$ through 
$\pi\circ\phi$.  By investigating $\xi:[0,T)\to T^{\perp}_{u_0}\widetilde M$, 
we obtain the following fact corresponding to Fact 1.  

\vspace{0.5truecm}

\noindent
{\bf Theorem A.} {\sl Let $M$ be an equifocal submanifold in a symmetric space 
$G/K$ of compact type.  Then the following statements ${\rm (i)}$ and 
${\rm (ii)}$ hold:

${\rm (i)}$ If $M$ is not minimal, then the mean curvature flow $M_t$ having 
$M$ as initial data collapses to a focal submanifold $F$ of $M$ in finite time.
Furtheremore, if $M$ is irreducible, the codimension of $M$ is 
greater than one and if the fibration of $M$ onto $F$ is spherical, then 
the flow $M_t$ has type I singularity.

${\rm (ii)}$ For any focal submanifold $F$ of $M$, there exists a parallel 
submanifold of $M$ collapsing to $F$ along the mean curvature flow and 
the set of all parallel submanifolds collapsing to $F$ along the mean 
curvature flow is a one-parameter $C^{\infty}$-family.}

\vspace{0.5truecm}

Also, we obtain the following fact corresponding to Fact 2 
for the mean curvature flow having a focal submanifold of an equifocal 
submanifold as initial data.  

\vspace{0.5truecm}

\noindent
{\it Remark 1.1.} U. Christ ([Ch]) showed that all irreducible equifocal 
submanifolds of codimension greater than one on symmetric spaces of compact 
type are homogeneous and hence they 
occur as principal orbits of hyperpolar actions.  On the other hand, 
A. Kollross ([Kol]) showed that all hyperpolar actions of cohomogeneity 
greater than one on irrdeucible symmetric spaces of compact type are orbit 
equivalent to Hermann actions on the space.  Hence all equifocal 
submanifolds of codimension greater than one on irreducible symmetric spaces 
of compact type occurs as principal orbits of Hermann actions.  Therefore 
they are classified completely.  
Hence we can show the statement of Theorem A by using this classification.  
However, it is very important to prove the statement conceptionally without 
the use of this classification.  In fact, Liu-Terng ([LT]) prove the results 
in [LT] conceptionally without the use of the classification of isoparametric 
submanifolds in Euclidean spaces by J. Dadok ([D]).  
So, in this paper, we prove the statement of Theorem A conceptionally.  

\vspace{0.3truecm}

\noindent
{\bf Theorem B.} {\sl Let $M$ be as in the statement of Theorem A and 
$\sigma$ be a stratum of dimension greater than zero of $\partial C$ 
(which is a stratified space).  
Then the following statements ${\rm (i)}$ and ${\rm (ii)}$ hold:

${\rm (i)}$ For any non-minimal focal submanifold $F$ (of $M$) through 
$\sigma$, the maen curvature flow $F_t$ having $F$ as initial data collapses 
to a focal submanifold $F'$ (of $M$) through $\partial{\sigma}$ 
in finite time.  
Furthermore, if $M$ is irreducible, the codimension of $M$ is 
greater than one and if the fibration of $F$ onto $F'$ is spherical, then 
the flow $F_t$ has type I singularity.

${\rm (ii)}$ For any focal submanifold $F$ of $M$ through $\partial\sigma$, 
there exists a focal submanifold of $M$ through $\sigma$ collapsing to $F$ 
along the mean curvature flow, and 
the set of all focal submanifolds of $M$ through $\sigma$ collapsing to $F$ 
along the mean curvature flow is a one-parameter $C^{\infty}$-family.}

\vspace{0.5truecm}

Since focal submanifolds of $M$ through $0$-dimensional stratums of 
$\partial C$ are minimal, it follows from these theorems 
that $M$ collapses to a minimal focal submanifold of $M$ after finitely many 
times of collapses along the mean curvature flows.  \newpage
$$\begin{array}{c}
\displaystyle{
\begin{array}{llll}
\displaystyle{M_t\mathop{\longrightarrow}_{(t\to T_1)}}&
\displaystyle{\mathop{F^1}_{\rm non-min.}}&&\\
&\displaystyle{F^1_t\mathop{\longrightarrow}_{(t\to T_2)}
\mathop{F^2}_{\rm non-min.}}&&\\
&&\displaystyle{\ddots}&\\
&&&\displaystyle{F^{k-1}_t\mathop{\longrightarrow}_{(t\to T_k)}
\mathop{F^k}_{\rm min.}}
\end{array}
}\\
\displaystyle{
\left(
\begin{array}{l}
\displaystyle{F^1\,:\,{\rm a}\,\,{\rm focal}\,\,{\rm submanifold}\,\,
{\rm of}\,\,M}\\
\displaystyle{F^i\,:\,{\rm a}\,\,{\rm focal}\,\,{\rm submanifold}\,\,
{\rm of}\,\,F^{i-1}\,\,(i=2,\cdots,k)}
\end{array}
\right)}
\end{array}$$

\noindent
According to the homogeneity theorem for an equifocal submanifold by 
Christ [Ch], all irreducible equifocal submanifolds of codimension greater 
than one in symmetric spaces of compact type are homogeneous.  
Hence, according to the result by Heintze-Palais-Terng-Thorbergsson [HPTT], 
they are principal orbits of hyperpolar actions.  
Furthermore, according to the classification by Kollross [Kol] of hyperpolar 
actions on irreducible symmetric spaces of compact type, all hyperpolar 
actions of cohomogeneity greater than one on the symmetric spaces are Hermann 
actions.  
Therefore, all equifocal submanifolds of codimension greater than one in 
irreducible symmetric spaces of compact type are principal orbits of 
Hermann actions.  In the last section, we describe explicitly the mean 
curvature flows having orbits of Hermann actions of cohomogeneity two 
on irreducible symmetric spaces of compact type and rank two as initial data.  

\section{Preliminaries}
In this section, we briefly review the quantities associated with an 
isoparametric submanifold in an (infinite dimensional separable) Hilbert 
space, which was introduced by Terng [T2].  
Let $M$ be an isoparametric submanifold in a Hilbert space $V$.  

\vspace{0.3truecm}

\noindent
{\bf 2.1. Principal curvatures, curvature normals and curvature distributions}
Let $E_0$ and $E_i$ ($i\in I$) be all the curvature distributions of $M$, 
where $E_0$ is defined by 
$(E_0)_x=\displaystyle{\mathop{\cap}_{v\in T^{\perp}_xM}
{\rm Ker}\,A_v}\,(x\in M)$.  For each $x\in M$, we have 
$T_xM=\overline{(E_0)_x\oplus
\displaystyle{\left(\mathop{\oplus}_{i\in I}(E_i)_x\right)}}$, 
which is the common eigenspace decomposition of $A_v$'s 
($v\in T^{\perp}_xM$).  
Also, let $\lambda_i$ ($i\in I$) be the principal curvatures of $M$, that is, 
$\lambda_i$ is the section of the dual bundle $(T^{\perp}M)^{\ast}$ of 
$T^{\perp}M$ such that $A_v\vert_{(E_i)_x}=(\lambda_i)_x(v){\rm id}$ 
holds for any $x\in M$ and any $v\in T^{\perp}_xM$, and ${\bf n}_i$ be 
the curvature normal corresponding to $\lambda_i$, that is, 
$\lambda_i(\cdot)=\langle{\bf n}_i,\cdot\rangle$.  

\vspace{0.3truecm}

\noindent
{\bf 2.2. The Coxeter group associated with an isoparametric submanifold}
Denote by ${\it l}^x_i$ the affine hyperplane $(\lambda_i)_x^{-1}(1)$ in 
$T^{\perp}_xM$.  The focal set of $M$ at $x$ is equal to the sum 
$\displaystyle{\mathop{\cup}_{i\in I}(x+{\it l}^x_i)}$ of the affine 
hyperplanes $x+{\it l}_i^x$'s ($i\in I$) in the affine subspace 
$x+T^{\perp}_xM$ of $V$.  Each affine hyperplane ${\it l}_i^x$ is called a 
focal hyperplane of $M$ at $x$.  Let $W$ be the group generated by the 
reflection $R_i^x$'s ($i\in I$) with respect to ${\it l}_i^x$.  
This group is independent of the choice $x$ of $M$ up to group isomorphism.  
This group is called the Coxeter group associated with $M$.  
The fundamental domain of the Coxeter group containing the zero vector of 
$T^{\perp}_xM$ is given by 
$\{v\in T^{\perp}_xM\,\vert\,\lambda_i(v)<1\,(i\in I)\}$.  

\vspace{0.3truecm}

\noindent
{\bf 2.3. Principal curvatures of parallel submanifolds}
Let $M_w$ be the parallel submanifold of $M$ for a (non-focal) parallel normal 
vector field $w$, that is, $M_w=\eta_w(M)$, where $\eta_w$ is the end-point 
map for $w$.  Denote by $A^w$ the shape tensor of $M_w$.  
This submanifold $M_w$ also is isoparametric and 
$A^w_v\vert_{\eta_{w\ast}(E_i)_x}
=\frac{(\lambda_i)_x(v)}{1-(\lambda_i)_x(w_x)}{\rm id}
\,\,(i\in I)$ for any $v\in T^{\perp}_{\eta_w(x)}M_w$, that is, 
$\frac{\lambda_i}{1-\lambda_i(w)}$'s ($i\in I$) are the principal curvatures 
of $M_w$ and hence $\frac{{\bf n}_i}{1-\lambda_i(w)}$'s ($i\in I$) are 
the curvature normals of $M_w$, where we identify $T^{\perp}_{\eta_w(x)}M_w$ 
with $T^{\perp}_xM$.  

\vspace{0.5truecm}

Let $M$ be a (general) submanifold in a Hilbert space $V$.  

\vspace{0.3truecm}

\noindent
{\bf 2.4. The mean curvature vector of a regularizable submanifold}
Assume that $M$ is regularizable in sense of [HLO], that is, 
for each normal vector $v$ of $M$,  the regularizable trace ${\rm Tr}_r\,A_v$ 
and ${\rm Tr}\,A_v^2$ exist, where ${\rm Tr}_r\,A_v$ is defined by 
${\rm Tr}_r\,A_v:=\sum\limits_{i=1}^{\infty}(\mu^+_i+\mu^-_i)$ 
($\mu^-_1\leq\mu^-_2\leq\cdots\leq 0\leq\cdots\leq\mu^+_2\leq\mu^+_1\,:\,$
the spectrum of $A_v$).  
Then the mean curvature vector $H$ of $M$ is defined by 
$\langle H,v\rangle={\rm Tr}_r\,A_v\,\,(\forall\,v\in T^{\perp}M)$.  

\vspace{0.3truecm}

\noindent
{\bf 2.5. The mean curvature flow for a regularizable submanifold}
Let $f_t:M\hookrightarrow V$ ($0\leq t<T$) be a $C^{\infty}$-family of 
regularizable submanifold immersions into $V$.  Denote by $H_t$ the mean 
curvature vector of $f_t$.  Define a map $F:M\times[0,T)\to V$ by 
$F(x,t):=f_t(x)$ ($(x,t)\in M\times[0,T)$).  If 
$\frac{\partial F}{\partial t}=H_t$ holds, then we call $f_t$ ($0\leq t<T$) 
the {\it mean curvature flow}.  In particular, if $f_t$'s are embeddings, then 
we set $M_t:=f_t(M)$ and $M_t$ (rather than $f_t$) the 
{\it mean curvature flow}.  Note that, for a given regularizable submanifold 
immersion $f:M\hookrightarrow V$, there does not necessarily exist the mean 
curvature flow having $f$ as initial data in short time.  
Furthermore, we note that, for the restriction $f\vert_U$ of $f$ to 
a relative compact domain $U$ of $M$, the existence of the mean curvature flow 
having $f\vert_U$ as initial data in short time is not assured 
for the sake of absence of the infinite dimensional vector bundle version 
of the Hamilton's theorem (Theorem 5.1 of [Ha]) for the existence and the 
uniqueness of solutions of a certain kind of evolution equation 
(which includes the Ricci flow equation and the mean curvature flow equation 
(of finite dimensional case)) in short time.  

\vspace{0.3truecm}

Let $M$ be an equifocal submanifold in a symmetric space $G/K$ of compact type 
and set $\widetilde M:=(\pi\circ\phi)^{-1}(M)$, where $\pi$ is the natural 
projection of $G$ onto $G/K$ and $\phi:H^0([0,1],\mathfrak g)\to G$ is the 
parallel transport map for $G$.  

\vspace{0.3truecm}

\noindent
{\bf 2.6. The mean curvature vector of the lifted submanifold} 
Denote by $\widetilde H$ (resp. $H$) the mean 
curvature vector of $\widetilde M$ (resp. $M$).  Then $\widetilde M$ is 
a regularizable isoparametric submanifold and $\widetilde H$ 
is equal to the horizontal lift of $H^L$ of $H$ 
%%($n:={\rm dim}\,M$) 
(see Lemma 5.2 of [HLO]).  

\section{Proofs of Theorems A and B}
In this section, we prove Theorems A and B.  
Let $M$ be an equifocal submanifold in 
a symmetric space $G/K$ of compact type, 
$\pi:G\to G/K$ be the natural 
projection and $\phi$ be the parallel transport map for $G$.  Set 
$\widetilde M:=(\pi\circ\phi)^{-1}(M)$.  
Take $u_0\in\widetilde M$ and set $x_0:=(\pi\circ\phi)(u_0)$.  We identify 
$T^{\perp}_{x_0}M$ with $T^{\perp}_{u_0}\widetilde M$.  
Let $\widetilde C(\subset T^{\perp}_{u_0}\widetilde M=T^{\perp}_{x_0}M)$ 
be the fundamental domain of the Coxeter group of $\widetilde M$ at 
$u_0$ containing the zero vector ${\bf 0}$ of 
$T^{\perp}_{u_0}\widetilde M(=T^{\perp}_{x_0}M)$ and set 
$C:=\exp^{\perp}(\widetilde C)$, where $\exp^{\perp}$ is 
the normal exponential map of $M$.  
Denote by $H$ (resp. 
$\widetilde H$) the mean curvature vector of $M$ (resp. $\widetilde M$).  
The mean curvature vector $H$ and $\widetilde H$ are a parallel normal 
vector field of $M$ and $\widetilde M$, respectively.  
Let $w$ be a parallel normal vector field of $M$ and $w^L$ be 
the horizontal lift of $w$ to $H^0([0,1],\mathfrak g)$, which is a parallel 
normal vector field of $\widetilde M$.  Denote by $M_w$ (resp. 
$\widetilde M_{w^L}$) the parallel (or focal) submanifold $\eta_w(M)$ (resp. 
$\eta_{w^L}(\widetilde M)$) of $M$ (resp. $\widetilde M$), where $\eta_w$ 
(resp. $\eta_{w^L}$) is the end-point map for $w$ (resp. $w^L$).  
Then we have $\widetilde M_{w^L}=(\pi\circ\phi)^{-1}(M_w)$.  
Denote by $H^w$ (resp. $\widetilde H^{w^L}$) the mean curvature vector of 
$M_w$ (resp. $\widetilde M_{w^L}$).  
Define a vector field $X$ on 
$\widetilde C\,(\subset T^{\perp}_{u_0}\widetilde M=T^{\perp}_{x_0}M)$ by 
$X_w:=(\widetilde H^{\widetilde w})_{u_0+w}$ ($w\in \widetilde C$), 
where $\widetilde w$ is the parallel normal 
vector field of $\widetilde M$ with $\widetilde w_{u_0}=w$.  Let 
$\xi\,:\,(-S,T)\to\widetilde C$ be the maximal integral curve of $X$ with 
$\xi(0)={\bf 0}$.  Note that $S$ and $T$ are possible be equal to $\infty$.  
Let $\widetilde{\xi(t)}$ be the parallel normal vector field of $M$ with 
$\widetilde{\xi(t)}_{x_0}=\xi(t)$.  
%%Let $M_t$ (resp. $\widetilde M_t$) be the 
%%mean curvature flow having $M$ (resp. $\widetilde M$) as initial data.  

\vspace{0.3truecm}

\noindent
{\bf Lemma 3.1.} {\sl The family $\widetilde M_{\widetilde{\xi(t)}^L}$ 
($0\leq t<T$) is the mean curvature flow having $\widetilde M$ as initial data 
and $M_{\widetilde{\xi(t)}}$ ($0\leq t<T$) is the mean curvature flow having 
$M$ as initial data.}

\vspace{0.3truecm}

\noindent
{\it Proof.} Fix $t_0\in[0,T)$.  Define a flow $F:\widetilde M\times[0,T)\to 
H^0([0,1],\mathfrak g)$ by $F(u,t):=\eta_{\widetilde{\xi(t)}^L}(u)$ 
($(u,t)\in\widetilde M\times[0,T)$) and $f_t:\widetilde M\to 
H^0([0,1],\mathfrak g)$ ($0\leq t<T$) by $f_t(u):=F(u,t)$ 
($u\in\widetilde M$).  
Here we note that $f_t(\widetilde M)=\widetilde M_{\widetilde{\xi(t)}^L}$.  
For simplicity, denote by ${\widetilde H}^{t_0}$ the mean 
curvature vector of ${\widetilde M}_{\widetilde{\xi(t_0)}^L}$.  It is easy to 
show that $F_{\ast}((\frac{\partial}{\partial t})_{(\cdot,t_0)})$ is a 
parallel normal vector field of ${\widetilde M}_{\widetilde{\xi(t_0)}^L}$ and 
that $F_{\ast}((\frac{\partial}{\partial t})_{(u_0,t_0)})
=({\widetilde H}^{t_0})_{f_{t_0}(u_0)}$.  On the other hand, since 
${\widetilde M}_{\widetilde{\xi(t_0)}^L}$ is isoparametric, 
${\widetilde H}^{t_0}$ 
is also a parallel normal vector field of 
${\widetilde M}_{\widetilde{\xi(t_0)}^L}$.  Hence we have 
$F_{\ast}((\frac{\partial}{\partial t})_{(\cdot,t_0)})=\widetilde H^{t_0}$.  
Therefore, it follows from the arbitrariness of $t_0$ that 
$\widetilde M_{\widetilde{\xi(t)}^L}$ ($0\leq t<T$) is the mean curvature flow 
having $\widetilde M$ as initial data.  
Define a flow $\overline F:M\times [0,T)\to G/K$ by $\overline F(x,t):=
\eta_{\widetilde{\xi(t)}}(x)$ ($(x,t)\in M\times[0,T)$) and 
$\bar f_t:M\to G/K$ ($0\leq t<T$) by 
$\bar f_t(x):=\overline F(x,t)$ ($x\in M$).  
Here we note that $\bar f_t(M)=M_{\widetilde{\xi(t)}}$.  
For simplicity, denote by $H^t$ the mean curvature vector of 
$M_{\widetilde{\xi(t)}}$.  Fix $t_0\in[0,T)$.  Since 
$\widetilde M_{\widetilde{\xi(t_0)}^L}=(\pi\circ\phi)^{-1}
(M_{\widetilde{\xi(t_0)}})$, we have 
$(H^{t_0})^L=\widetilde H^{t_0}$.  On the other hand, we have 
$\overline F_{\ast}((\frac{\partial}{\partial t})_{(\cdot,t_0)})^L
=F_{\ast}((\frac{\partial}{\partial t})_{(\cdot,t_0)})(=\widetilde H^{t_0})$.  
Hence we have $\overline F_{\ast}((\frac{\partial}{\partial t})_{(\cdot,t_0)})
=H^{t_0}$.  Therefore, it follows from the arbitrariness of $t_0$ that 
$M_{\widetilde{\xi(t)}}$ ($0\leq t<T$) is the mean curvature flow having $M$ 
as initial data.  \hspace{1.5truecm}q.e.d.

\vspace{0.3truecm}

\noindent
{\it Proof of Theorem A.} Clearly we suffice to show the statement 
of Theorem A in the case where $M$ is full.  Hence, in the sequel, 
we assume that $M$ is full.  
Denote by $\Lambda$ the set of all principal curvatures of $\widetilde M$.  
Set $r:={\rm codim}\,M$.  It is shown that the set of all focal hyperplanes of 
$\widetilde M$ is given as the sum of finite pieces of infinite parallel 
families consisting of hyperplanes in $T^{\perp}_{u_0}\widetilde M$ which 
arrange at equal intervals.  Let 
$\{{\it l}_{aj}\,\vert\,j\in{\Bbb Z}\}$ ($1\leq a\leq \bar r$) be 
the finite pieces of infinite parallel families consisting of hyperplanes 
in $T^{\perp}_{u_0}\widetilde M$.  
Since ${\it l}_{aj}$'s ($j\in{\Bbb Z}$) arrange at equal intervals, 
we can express as 
$\displaystyle{
\Lambda=\mathop{\cup}_{a=1}^{\bar r}\{\frac{\lambda_a}{1+b_aj}\,\vert\,
j\in{\bf Z}\}}$, where $\lambda_a$'s and $b_a$'s are parallel sections of 
$(T^{\perp}\widetilde M)^{\ast}$ and positive constants greater than one, 
respectively, 
which are defined by $((\lambda_a)_{u_0})^{-1}(1+b_aj)={\it l}_{aj}$.  
For simplicity, we set $\lambda_{aj}:=\frac{\lambda_a}{1+b_aj}$.  
Denote by ${\bf n}_{aj}$ and $E_{aj}$ the curvature normal and the curvature 
distribution corresponding to $\lambda_{aj}$, respectively.  
The fundamental domain $\widetilde C$ of the Coxeter group of $\widetilde M$ 
at $u_0$ is given by 
$$\widetilde C=\{w\in T^{\perp}_{u_0}\widetilde M\,\vert\,\lambda_a(w)<1\,\,
(1\leq a\leq\bar r)\}.$$
It is shown that, for each $a$, $\lambda_{a,2j}$'s ($j\in{\Bbb Z}$) have the 
same multiplicity and so are also $\lambda_{a,2j+1}$'s ($j\in{\Bbb Z}$).  
Denote by $m_a^e$ and $m_a^o$ the multiplicities of $\lambda_{a,2j}$ and 
$\lambda_{a,2j+1}$, respectively.  Take a parallel normal vector field 
$w$ of $\widetilde M$ with $w_{u_0}\in\widetilde C$.  Denote by 
${\widetilde A}^w$ (resp. ${\widetilde H}^w$) the shape tensor (resp. 
the mean curvature vector) of the parallel submanifold ${\widetilde M}_w$.  
Since ${\widetilde A}^w_v\vert_{\eta_{w\ast}(E_{aj})_u}
=\frac{(\lambda_{aj})_u(v)}{1-(\lambda_{aj})_u(w_u)}\,{\rm id}$ 
$(v\in T^{\perp}_u\widetilde M$), we have 
$$\begin{array}{l}
\displaystyle{{\rm Tr}_r\widetilde A^w_v=\sum_{a=1}^{\bar r}\left(
\sum_{j\in{\bf Z}}\frac{m_a^e(\lambda_{a,2j})_u(v)}{1-(\lambda_{a,2j})_u(w_u)}+
\sum_{j\in{\bf Z}}\frac{m_a^o(\lambda_{a,2j+1})_u(v)}
{1-(\lambda_{a,2j+1})_u(w_u)}
\right)}\\
\hspace{1.1truecm}\displaystyle{=\sum_{a=1}^{\bar r}\left(
m_a^e\cot\frac{\pi}{2b_a}(1-(\lambda_a)_u(w_u))-m_a^o\tan\frac{\pi}{2b_a}
(1-(\lambda_a)_u(w_u))\right)\frac{\pi}{2b_a}(\lambda_a)_u(v),}
\end{array}$$
where we use the relation 
$\cot\,\frac{\theta}{2}=\sum\limits_{j\in{\bf Z}}
\frac{2}{\theta+2j\pi}$.  Therefore we have 
$$\begin{array}{l}
\displaystyle{\widetilde H^w=\sum_{a=1}^{\bar r}\left(
m_a^e\cot\frac{\pi}{2b_a}(1-\lambda_a(w))
-m_a^o\tan\frac{\pi}{2b_a}(1-\lambda_a(w))\right)\frac{\pi}{2b_a}{\bf n}_a,}
\end{array}\leqno{(3.1)}$$
where ${\bf n}_a$ is the curvature normal corresponding to $\lambda_a$.  
Denote by $\widetilde{\sigma}_a$ ($1\leq a\leq\bar r$) the maximal dimensional 
stratum $(\lambda_a)_{u_0}^{-1}(1)$ of $\partial\widetilde C$.  
Fix $a_0\in\{1,\cdots,\bar r\}$.  Take $w_0\in\widetilde{\sigma}_{a_0}$ and 
and $w'_0\in\widetilde C$ near $w_0$ such that $w_0-w_0'$ is 
normal to $\widetilde{\sigma}_{a_0}$.  Set $w_0^{\varepsilon}
:=\varepsilon w'_0+(1-\varepsilon)w_0$ for $\varepsilon\in(0,1)$.  
Then we have $\lim\limits_{\varepsilon\to+0}(\lambda_{a_0})_{u_0}
(w^{\varepsilon}_0)=1$ and 
$\displaystyle{\mathop{\sup}_{0<\varepsilon<1}
(\lambda_a)_{u_0}(w^{\varepsilon}_0)<1}$ for each 
$a\in\{1,\cdots,\bar r\}\setminus\{a_0\}$.  
Hence we have 
$
%\begin{array}{l}\displaystyle{
\lim_{\varepsilon\to+0}\cot\frac{\pi}{2b_{a_0}}
(1-(\lambda_{a_0})_{u_0}(w^{\varepsilon}_0))=\infty$
%}\\\displaystyle{
and
$\mathop{\sup}_{0<\varepsilon<1}\coth
(1-(\lambda_a)_{u_0}(w^{\varepsilon}_0))\,<\,\infty\,\,
(a\in\{1,\cdots,\bar r\}\setminus\{a_0\})$.  
%}\end{array}
Therefore, we have $\lim\limits_{\varepsilon\to+0}
\frac{X_{w_0^{\varepsilon}}}{\vert\vert X_{w_0^{\varepsilon}}\vert\vert}$ is 
the outward unit normal vector of $\widetilde{\sigma}_{a_0}$.  Also we have 
$\lim\limits_{\varepsilon\to +0}\vert\vert X_{w_0^{\varepsilon}}\vert\vert
=\infty$.  From these facts, $X$ is as in Fig. 1 on a sufficiently small 
collar neighborhood of $\widetilde{\sigma}_{a_0}$.  

\vspace{0.3truecm}

\centerline{
%WinTpicVersion3.08
\unitlength 0.1in
\begin{picture}( 19.2700, 18.7400)( 21.6900,-28.1400)
% POLYGON 2 0 3 0
% 5 2999 1404 2169 2814 3799 2814 3799 2814 2999 1404
% 
\special{pn 8}%
\special{pa 3000 1404}%
\special{pa 2170 2814}%
\special{pa 3800 2814}%
\special{pa 3800 2814}%
\special{pa 3000 1404}%
\special{fp}%
% VECTOR 1 0 3 0
% 2 3609 2574 4096 2292
% 
\special{pn 13}%
\special{pa 3610 2574}%
\special{pa 4096 2292}%
\special{fp}%
\special{sh 1}%
\special{pa 4096 2292}%
\special{pa 4028 2308}%
\special{pa 4050 2320}%
\special{pa 4048 2344}%
\special{pa 4096 2292}%
\special{fp}%
% STR 2 0 3 0
% 3 3500 1240 3500 1340 0 0
% 
\put(35.0000,-13.4000){\makebox(0,0)[lb]{}}%
% STR 2 0 3 0
% 3 4040 1730 4040 1830 2 0
% $X$
\put(40.4000,-18.3000){\makebox(0,0)[lb]{$X$}}%
% VECTOR 1 0 3 0
% 2 3440 2280 3927 1998
% 
\special{pn 13}%
\special{pa 3440 2280}%
\special{pa 3928 1998}%
\special{fp}%
\special{sh 1}%
\special{pa 3928 1998}%
\special{pa 3860 2014}%
\special{pa 3882 2026}%
\special{pa 3880 2050}%
\special{pa 3928 1998}%
\special{fp}%
% VECTOR 1 0 3 0
% 2 3290 2020 3777 1738
% 
\special{pn 13}%
\special{pa 3290 2020}%
\special{pa 3778 1738}%
\special{fp}%
\special{sh 1}%
\special{pa 3778 1738}%
\special{pa 3710 1754}%
\special{pa 3732 1766}%
\special{pa 3730 1790}%
\special{pa 3778 1738}%
\special{fp}%
% VECTOR 1 0 3 0
% 2 3130 1730 3617 1448
% 
\special{pn 13}%
\special{pa 3130 1730}%
\special{pa 3618 1448}%
\special{fp}%
\special{sh 1}%
\special{pa 3618 1448}%
\special{pa 3550 1464}%
\special{pa 3572 1476}%
\special{pa 3570 1500}%
\special{pa 3618 1448}%
\special{fp}%
% VECTOR 2 1 3 0
% 2 3330 1150 3080 1520
% 
\special{pn 8}%
\special{pa 3330 1150}%
\special{pa 3080 1520}%
\special{da 0.070}%
\special{sh 1}%
\special{pa 3080 1520}%
\special{pa 3134 1476}%
\special{pa 3110 1476}%
\special{pa 3102 1454}%
\special{pa 3080 1520}%
\special{fp}%
% STR 2 0 3 0
% 3 3340 1010 3340 1110 2 0
% $\widetilde{\sigma}_{a_0}$
\put(33.4000,-11.1000){\makebox(0,0)[lb]{$\widetilde{\sigma}_{a_0}$}}%
% STR 2 0 3 0
% 3 2950 2140 2950 2240 1 0
% $\widetilde C$
\put(29.5000,-22.4000){\makebox(0,0)[lt]{$\widetilde C$}}%
\end{picture}%
}

\vspace{0.6truecm}

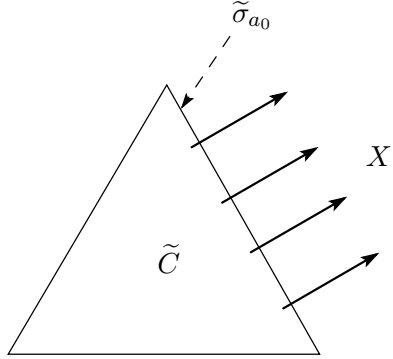
\centerline{{\bf Fig. 1.}}

\vspace{0.5truecm}

\noindent
Define a function $\rho$ over $\widetilde C$ by 
$$\begin{array}{l}
\displaystyle{\rho(w):=-\sum_{a=1}^{\bar r}\left(m_a^e
\log\sin\frac{\pi}{2b_a}(1-(\lambda_a)_{u_0}(w))\right.}\\
\hspace{2.4truecm}\displaystyle{
\left.+m_a^o\log\cos\frac{\pi}{2b_a}(1-(\lambda_a)_{u_0}(w))\right)
\qquad(w\in\widetilde C).}
\end{array}$$
Let $(x_1,\cdots,x_r)$ be the Euclidean coordinate of 
$T^{\perp}_{u_0}\widetilde M$.  For simplicity, set 
$\partial_i:=\frac{\partial}{\partial x_i}$ ($i=1,\cdots,r$).  
Then it follows from the definition of $X$ and $(3.1)$ that 
$(\partial_i\rho)(w)=\langle X_w,\partial_i\rangle$ 
($w\in\widetilde C,\,\,i=,\cdots,r$), that is, ${\rm grad}\,\rho=X$.  
Also we have 
$$\begin{array}{l}
\displaystyle{(\partial_i\partial_j\rho)(w)=
\sum_{a=1}^{\bar r}\left(\frac{m_a^e}
{\sin^2\frac{\pi}{2b_a}(1-(\lambda_a)_{u_0}(w))}
+\frac{m_a^o}{\cos^2\frac{\pi}{2b_a}(1-(\lambda_a)_{u_0}(w))}\right)}\\
\hspace{3.4truecm}\displaystyle{\times
\frac{\pi^2}{4b_a^2}(\lambda_a)_{u_0}(\partial_i)(\lambda_a)_{u_0}
(\partial_j).}
\end{array}$$
It follows from this relation that $\rho$ is downward convex.  
Also it is shown that $\rho(w)\to\infty$ as $w\to\partial\widetilde C$.  
Hence we see that $\rho$ has the only minimal point.  
Denote by $w_0$ this minimal point.  It is clear that $X_{w_0}=0$.  
From these facts and the fact that $X$ is as in Fig. 1 on a sufficently small 
collar neighborhood of each maximal dimensional stratum $\widetilde{\sigma}$ 
of $\partial\widetilde C$, $\rho$ and $X$ are as in Fig. 2.  

\vspace{0.3truecm}

\centerline{
%WinTpicVersion3.08
\unitlength 0.1in
\begin{picture}( 62.2000, 31.7000)(  1.9000,-39.5000)
% POLYGON 2 0 3 0
% 6 2390 3000 1180 3600 2800 3600 2400 3000 2400 3000 2390 3000
% 
\special{pn 8}%
\special{pa 2390 3000}%
\special{pa 1180 3600}%
\special{pa 2800 3600}%
\special{pa 2400 3000}%
\special{pa 2400 3000}%
\special{pa 2390 3000}%
\special{fp}%
% POLYGON 2 0 3 0
% 6 2400 1440 1190 2040 2810 2040 2410 1440 2410 1440 2400 1440
% 
\special{pn 8}%
\special{pa 2400 1440}%
\special{pa 1190 2040}%
\special{pa 2810 2040}%
\special{pa 2410 1440}%
\special{pa 2410 1440}%
\special{pa 2400 1440}%
\special{fp}%
% LINE 2 0 3 0
% 2 1180 3600 1180 1380
% 
\special{pn 8}%
\special{pa 1180 3600}%
\special{pa 1180 1380}%
\special{fp}%
% LINE 2 0 3 0
% 2 2800 3600 2800 1380
% 
\special{pn 8}%
\special{pa 2800 3600}%
\special{pa 2800 1380}%
\special{fp}%
% DOT 1 0 3 0
% 2 2200 3390 2200 3390
% 
\special{pn 13}%
\special{sh 1}%
\special{ar 2200 3390 10 10 0  6.28318530717959E+0000}%
\special{sh 1}%
\special{ar 2200 3390 10 10 0  6.28318530717959E+0000}%
% LINE 2 0 3 0
% 2 2400 3000 2400 2080
% 
\special{pn 8}%
\special{pa 2400 3000}%
\special{pa 2400 2080}%
\special{fp}%
% LINE 2 0 3 0
% 2 2400 1970 2400 1020
% 
\special{pn 8}%
\special{pa 2400 1970}%
\special{pa 2400 1020}%
\special{fp}%
% ELLIPSE 1 0 3 0
% 4 2970 3850 5420 1500 2100 1000 1210 1510
% 
\special{pn 13}%
\special{ar 2970 3850 2450 2350  4.0875027 4.4275224}%
% ELLIPSE 1 0 3 0
% 4 2120 1880 3460 1990 560 2170 3560 2190
% 
\special{pn 13}%
\special{ar 2120 1880 1340 110  1.2064650 1.9866020}%
% ELLIPSE 1 0 3 0
% 4 2220 2740 2770 1500 2840 1540 2480 1170
% 
\special{pn 13}%
\special{ar 2220 2740 550 1240  5.0698994 5.5740280}%
% ELLIPSE 1 0 3 0
% 4 2350 1630 2450 1670 2540 1480 2020 1440
% 
\special{pn 13}%
\special{ar 2350 1630 100 40  4.1052041 5.1813557}%
% ELLIPSE 1 0 3 0
% 4 1660 1960 1790 1930 1360 1940 1020 2090
% 
\special{pn 13}%
\special{ar 1660 1960 130 30  2.4201139 3.4238501}%
% ELLIPSE 1 0 3 0
% 4 2580 1940 2640 1980 2790 2240 3100 1780
% 
\special{pn 13}%
\special{ar 2580 1940 60 40  5.8507775 6.2831853}%
\special{ar 2580 1940 60 40  0.0000000 1.1341692}%
% LINE 2 2 3 0
% 2 2200 2810 2200 3380
% 
\special{pn 8}%
\special{pa 2200 2810}%
\special{pa 2200 3380}%
\special{dt 0.045}%
% ELLIPSE 1 0 3 0
% 4 2200 890 2730 2780 2200 3150 3750 2100
% 
\special{pn 13}%
\special{ar 2200 890 530 1890  0.2153192 1.5707963}%
% ELLIPSE 1 0 3 0
% 4 2200 990 3000 2780 10 1820 2200 3920
% 
\special{pn 13}%
\special{ar 2200 990 800 1790  1.5707963 2.9737795}%
% ELLIPSE 1 2 3 0
% 4 2200 200 2440 2780 2200 2900 2920 2500
% 
\special{pn 13}%
\special{ar 2200 200 240 2580  0.2889064 0.2974171}%
\special{ar 2200 200 240 2580  0.3229490 0.3314596}%
\special{ar 2200 200 240 2580  0.3569915 0.3655022}%
\special{ar 2200 200 240 2580  0.3910341 0.3995447}%
\special{ar 2200 200 240 2580  0.4250766 0.4335873}%
\special{ar 2200 200 240 2580  0.4591192 0.4676298}%
\special{ar 2200 200 240 2580  0.4931617 0.5016724}%
\special{ar 2200 200 240 2580  0.5272043 0.5357149}%
\special{ar 2200 200 240 2580  0.5612469 0.5697575}%
\special{ar 2200 200 240 2580  0.5952894 0.6038000}%
\special{ar 2200 200 240 2580  0.6293320 0.6378426}%
\special{ar 2200 200 240 2580  0.6633745 0.6718852}%
\special{ar 2200 200 240 2580  0.6974171 0.7059277}%
\special{ar 2200 200 240 2580  0.7314596 0.7399703}%
\special{ar 2200 200 240 2580  0.7655022 0.7740128}%
\special{ar 2200 200 240 2580  0.7995447 0.8080554}%
\special{ar 2200 200 240 2580  0.8335873 0.8420979}%
\special{ar 2200 200 240 2580  0.8676298 0.8761405}%
\special{ar 2200 200 240 2580  0.9016724 0.9101830}%
\special{ar 2200 200 240 2580  0.9357149 0.9442256}%
\special{ar 2200 200 240 2580  0.9697575 0.9782681}%
\special{ar 2200 200 240 2580  1.0038000 1.0123107}%
\special{ar 2200 200 240 2580  1.0378426 1.0463532}%
\special{ar 2200 200 240 2580  1.0718852 1.0803958}%
\special{ar 2200 200 240 2580  1.1059277 1.1144383}%
\special{ar 2200 200 240 2580  1.1399703 1.1484809}%
\special{ar 2200 200 240 2580  1.1740128 1.1825235}%
\special{ar 2200 200 240 2580  1.2080554 1.2165660}%
\special{ar 2200 200 240 2580  1.2420979 1.2506086}%
\special{ar 2200 200 240 2580  1.2761405 1.2846511}%
\special{ar 2200 200 240 2580  1.3101830 1.3186937}%
\special{ar 2200 200 240 2580  1.3442256 1.3527362}%
\special{ar 2200 200 240 2580  1.3782681 1.3867788}%
\special{ar 2200 200 240 2580  1.4123107 1.4208213}%
\special{ar 2200 200 240 2580  1.4463532 1.4548639}%
\special{ar 2200 200 240 2580  1.4803958 1.4889064}%
\special{ar 2200 200 240 2580  1.5144383 1.5229490}%
\special{ar 2200 200 240 2580  1.5484809 1.5569915}%
% ELLIPSE 1 0 3 0
% 4 2200 600 2520 2780 1470 2370 2200 3580
% 
\special{pn 13}%
\special{ar 2200 600 320 2180  1.5707963 2.7994368}%
% DOT 1 0 3 0
% 2 2200 2780 2200 2780
% 
\special{pn 13}%
\special{sh 1}%
\special{ar 2200 2780 10 10 0  6.28318530717959E+0000}%
\special{sh 1}%
\special{ar 2200 2780 10 10 0  6.28318530717959E+0000}%
% DOT 1 0 3 0
% 2 1950 1990 1950 1990
% 
\special{pn 13}%
\special{sh 1}%
\special{ar 1950 1990 10 10 0  6.28318530717959E+0000}%
\special{sh 1}%
\special{ar 1950 1990 10 10 0  6.28318530717959E+0000}%
% VECTOR 2 0 3 0
% 2 3110 2160 2480 2300
% 
\special{pn 8}%
\special{pa 3110 2160}%
\special{pa 2480 2300}%
\special{fp}%
\special{sh 1}%
\special{pa 2480 2300}%
\special{pa 2550 2306}%
\special{pa 2532 2288}%
\special{pa 2542 2266}%
\special{pa 2480 2300}%
\special{fp}%
% VECTOR 2 0 3 0
% 2 2980 3070 2520 3480
% 
\special{pn 8}%
\special{pa 2980 3070}%
\special{pa 2520 3480}%
\special{fp}%
\special{sh 1}%
\special{pa 2520 3480}%
\special{pa 2584 3452}%
\special{pa 2560 3446}%
\special{pa 2556 3422}%
\special{pa 2520 3480}%
\special{fp}%
% DOT 1 0 3 0
% 2 2400 1590 2400 1590
% 
\special{pn 13}%
\special{sh 1}%
\special{ar 2400 1590 10 10 0  6.28318530717959E+0000}%
\special{sh 1}%
\special{ar 2400 1590 10 10 0  6.28318530717959E+0000}%
% VECTOR 2 0 3 0
% 2 1180 3600 1180 780
% 
\special{pn 8}%
\special{pa 1180 3600}%
\special{pa 1180 780}%
\special{fp}%
\special{sh 1}%
\special{pa 1180 780}%
\special{pa 1160 848}%
\special{pa 1180 834}%
\special{pa 1200 848}%
\special{pa 1180 780}%
\special{fp}%
% POLYGON 2 0 3 0
% 5 4999 1794 4169 3204 5799 3204 5799 3204 4999 1794
% 
\special{pn 8}%
\special{pa 5000 1794}%
\special{pa 4170 3204}%
\special{pa 5800 3204}%
\special{pa 5800 3204}%
\special{pa 5000 1794}%
\special{fp}%
% DOT 1 0 3 0
% 2 4999 2804 4999 2804
% 
\special{pn 13}%
\special{sh 1}%
\special{ar 5000 2804 10 10 0  6.28318530717959E+0000}%
\special{sh 1}%
\special{ar 5000 2804 10 10 0  6.28318530717959E+0000}%
% VECTOR 1 0 3 0
% 2 5149 2714 5299 2614
% 
\special{pn 13}%
\special{pa 5150 2714}%
\special{pa 5300 2614}%
\special{fp}%
\special{sh 1}%
\special{pa 5300 2614}%
\special{pa 5232 2634}%
\special{pa 5256 2644}%
\special{pa 5256 2668}%
\special{pa 5300 2614}%
\special{fp}%
% VECTOR 1 0 3 0
% 2 5129 2114 5599 1804
% 
\special{pn 13}%
\special{pa 5130 2114}%
\special{pa 5600 1804}%
\special{fp}%
\special{sh 1}%
\special{pa 5600 1804}%
\special{pa 5532 1824}%
\special{pa 5554 1834}%
\special{pa 5554 1858}%
\special{pa 5600 1804}%
\special{fp}%
% VECTOR 1 0 3 0
% 2 5289 2404 5759 2094
% 
\special{pn 13}%
\special{pa 5290 2404}%
\special{pa 5760 2094}%
\special{fp}%
\special{sh 1}%
\special{pa 5760 2094}%
\special{pa 5692 2114}%
\special{pa 5714 2124}%
\special{pa 5714 2148}%
\special{pa 5760 2094}%
\special{fp}%
% VECTOR 1 0 3 0
% 2 5449 2684 5919 2374
% 
\special{pn 13}%
\special{pa 5450 2684}%
\special{pa 5920 2374}%
\special{fp}%
\special{sh 1}%
\special{pa 5920 2374}%
\special{pa 5852 2394}%
\special{pa 5874 2404}%
\special{pa 5874 2428}%
\special{pa 5920 2374}%
\special{fp}%
% VECTOR 1 0 3 0
% 2 5609 2964 6079 2654
% 
\special{pn 13}%
\special{pa 5610 2964}%
\special{pa 6080 2654}%
\special{fp}%
\special{sh 1}%
\special{pa 6080 2654}%
\special{pa 6012 2674}%
\special{pa 6034 2684}%
\special{pa 6034 2708}%
\special{pa 6080 2654}%
\special{fp}%
% VECTOR 2 1 3 0
% 2 5689 3114 6241 3001
% 
\special{pn 8}%
\special{pa 5690 3114}%
\special{pa 6242 3002}%
\special{da 0.070}%
\special{sh 1}%
\special{pa 6242 3002}%
\special{pa 6172 2996}%
\special{pa 6190 3012}%
\special{pa 6180 3034}%
\special{pa 6242 3002}%
\special{fp}%
% VECTOR 2 1 3 0
% 2 5729 3164 6221 3439
% 
\special{pn 8}%
\special{pa 5730 3164}%
\special{pa 6222 3440}%
\special{da 0.070}%
\special{sh 1}%
\special{pa 6222 3440}%
\special{pa 6174 3390}%
\special{pa 6174 3414}%
\special{pa 6154 3424}%
\special{pa 6222 3440}%
\special{fp}%
% VECTOR 2 1 3 0
% 2 5669 3174 5936 3670
% 
\special{pn 8}%
\special{pa 5670 3174}%
\special{pa 5936 3670}%
\special{da 0.070}%
\special{sh 1}%
\special{pa 5936 3670}%
\special{pa 5922 3602}%
\special{pa 5912 3624}%
\special{pa 5888 3622}%
\special{pa 5936 3670}%
\special{fp}%
% VECTOR 1 0 3 0
% 2 5519 3164 5514 3728
% 
\special{pn 13}%
\special{pa 5520 3164}%
\special{pa 5514 3728}%
\special{fp}%
\special{sh 1}%
\special{pa 5514 3728}%
\special{pa 5536 3662}%
\special{pa 5514 3676}%
\special{pa 5496 3662}%
\special{pa 5514 3728}%
\special{fp}%
% VECTOR 1 0 3 0
% 2 5219 3164 5214 3728
% 
\special{pn 13}%
\special{pa 5220 3164}%
\special{pa 5214 3728}%
\special{fp}%
\special{sh 1}%
\special{pa 5214 3728}%
\special{pa 5236 3662}%
\special{pa 5214 3676}%
\special{pa 5196 3662}%
\special{pa 5214 3728}%
\special{fp}%
% VECTOR 1 0 3 0
% 2 4849 3154 4844 3718
% 
\special{pn 13}%
\special{pa 4850 3154}%
\special{pa 4844 3718}%
\special{fp}%
\special{sh 1}%
\special{pa 4844 3718}%
\special{pa 4866 3652}%
\special{pa 4844 3666}%
\special{pa 4826 3652}%
\special{pa 4844 3718}%
\special{fp}%
% VECTOR 1 0 3 0
% 2 4529 3154 4524 3718
% 
\special{pn 13}%
\special{pa 4530 3154}%
\special{pa 4524 3718}%
\special{fp}%
\special{sh 1}%
\special{pa 4524 3718}%
\special{pa 4546 3652}%
\special{pa 4524 3666}%
\special{pa 4506 3652}%
\special{pa 4524 3718}%
\special{fp}%
% VECTOR 2 1 3 0
% 2 4349 3154 4075 3647
% 
\special{pn 8}%
\special{pa 4350 3154}%
\special{pa 4076 3648}%
\special{da 0.070}%
\special{sh 1}%
\special{pa 4076 3648}%
\special{pa 4126 3598}%
\special{pa 4102 3600}%
\special{pa 4090 3580}%
\special{pa 4076 3648}%
\special{fp}%
% VECTOR 2 1 3 0
% 2 4239 3154 3809 3518
% 
\special{pn 8}%
\special{pa 4240 3154}%
\special{pa 3810 3518}%
\special{da 0.070}%
\special{sh 1}%
\special{pa 3810 3518}%
\special{pa 3874 3490}%
\special{pa 3850 3484}%
\special{pa 3848 3460}%
\special{pa 3810 3518}%
\special{fp}%
% VECTOR 2 1 3 0
% 2 4299 3064 3760 3227
% 
\special{pn 8}%
\special{pa 4300 3064}%
\special{pa 3760 3228}%
\special{da 0.070}%
\special{sh 1}%
\special{pa 3760 3228}%
\special{pa 3830 3228}%
\special{pa 3812 3212}%
\special{pa 3818 3190}%
\special{pa 3760 3228}%
\special{fp}%
% VECTOR 1 0 3 0
% 2 4359 2964 3891 2650
% 
\special{pn 13}%
\special{pa 4360 2964}%
\special{pa 3892 2650}%
\special{fp}%
\special{sh 1}%
\special{pa 3892 2650}%
\special{pa 3936 2704}%
\special{pa 3936 2680}%
\special{pa 3958 2672}%
\special{pa 3892 2650}%
\special{fp}%
% VECTOR 1 0 3 0
% 2 4549 2684 4081 2370
% 
\special{pn 13}%
\special{pa 4550 2684}%
\special{pa 4082 2370}%
\special{fp}%
\special{sh 1}%
\special{pa 4082 2370}%
\special{pa 4126 2424}%
\special{pa 4126 2400}%
\special{pa 4148 2392}%
\special{pa 4082 2370}%
\special{fp}%
% VECTOR 1 0 3 0
% 2 4709 2394 4241 2080
% 
\special{pn 13}%
\special{pa 4710 2394}%
\special{pa 4242 2080}%
\special{fp}%
\special{sh 1}%
\special{pa 4242 2080}%
\special{pa 4286 2134}%
\special{pa 4286 2110}%
\special{pa 4308 2102}%
\special{pa 4242 2080}%
\special{fp}%
% VECTOR 1 0 3 0
% 2 4849 2114 4381 1800
% 
\special{pn 13}%
\special{pa 4850 2114}%
\special{pa 4382 1800}%
\special{fp}%
\special{sh 1}%
\special{pa 4382 1800}%
\special{pa 4426 1854}%
\special{pa 4426 1830}%
\special{pa 4448 1822}%
\special{pa 4382 1800}%
\special{fp}%
% VECTOR 2 1 3 0
% 2 4969 1924 4628 1475
% 
\special{pn 8}%
\special{pa 4970 1924}%
\special{pa 4628 1476}%
\special{da 0.070}%
\special{sh 1}%
\special{pa 4628 1476}%
\special{pa 4652 1540}%
\special{pa 4660 1518}%
\special{pa 4684 1516}%
\special{pa 4628 1476}%
\special{fp}%
% VECTOR 2 1 3 0
% 2 4999 1874 5003 1310
% 
\special{pn 8}%
\special{pa 5000 1874}%
\special{pa 5004 1310}%
\special{da 0.070}%
\special{sh 1}%
\special{pa 5004 1310}%
\special{pa 4984 1378}%
\special{pa 5004 1364}%
\special{pa 5024 1378}%
\special{pa 5004 1310}%
\special{fp}%
% VECTOR 2 1 3 0
% 2 5039 1934 5355 1467
% 
\special{pn 8}%
\special{pa 5040 1934}%
\special{pa 5356 1468}%
\special{da 0.070}%
\special{sh 1}%
\special{pa 5356 1468}%
\special{pa 5302 1512}%
\special{pa 5326 1512}%
\special{pa 5334 1534}%
\special{pa 5356 1468}%
\special{fp}%
% VECTOR 1 0 3 0
% 2 4999 2474 5003 2294
% 
\special{pn 13}%
\special{pa 5000 2474}%
\special{pa 5004 2294}%
\special{fp}%
\special{sh 1}%
\special{pa 5004 2294}%
\special{pa 4982 2360}%
\special{pa 5002 2348}%
\special{pa 5022 2362}%
\special{pa 5004 2294}%
\special{fp}%
% VECTOR 1 0 3 0
% 2 4849 2704 4696 2608
% 
\special{pn 13}%
\special{pa 4850 2704}%
\special{pa 4696 2608}%
\special{fp}%
\special{sh 1}%
\special{pa 4696 2608}%
\special{pa 4742 2660}%
\special{pa 4742 2636}%
\special{pa 4764 2626}%
\special{pa 4696 2608}%
\special{fp}%
% VECTOR 1 0 3 0
% 2 4739 2904 4594 3011
% 
\special{pn 13}%
\special{pa 4740 2904}%
\special{pa 4594 3012}%
\special{fp}%
\special{sh 1}%
\special{pa 4594 3012}%
\special{pa 4660 2988}%
\special{pa 4638 2980}%
\special{pa 4636 2956}%
\special{pa 4594 3012}%
\special{fp}%
% VECTOR 1 0 3 0
% 2 4999 2974 4999 3154
% 
\special{pn 13}%
\special{pa 5000 2974}%
\special{pa 5000 3154}%
\special{fp}%
\special{sh 1}%
\special{pa 5000 3154}%
\special{pa 5020 3088}%
\special{pa 5000 3102}%
\special{pa 4980 3088}%
\special{pa 5000 3154}%
\special{fp}%
% VECTOR 1 0 3 0
% 2 5289 2914 5437 3014
% 
\special{pn 13}%
\special{pa 5290 2914}%
\special{pa 5438 3014}%
\special{fp}%
\special{sh 1}%
\special{pa 5438 3014}%
\special{pa 5394 2960}%
\special{pa 5394 2984}%
\special{pa 5372 2994}%
\special{pa 5438 3014}%
\special{fp}%
% STR 2 0 3 0
% 3 3500 1240 3500 1340 0 0
% 
\put(35.0000,-13.4000){\makebox(0,0)[lb]{}}%
% STR 2 0 3 0
% 3 3140 2880 3140 2980 4 0
% $\widetilde C$
\put(31.4000,-29.8000){\makebox(0,0)[rt]{$\widetilde C$}}%
% STR 2 0 3 0
% 3 2260 3340 2260 3440 4 0
% $w_0$
\put(22.6000,-34.4000){\makebox(0,0)[rt]{$w_0$}}%
% STR 2 0 3 0
% 3 3310 1980 3310 2080 4 0
% $\rho$
\put(33.1000,-20.8000){\makebox(0,0)[rt]{$\rho$}}%
% STR 2 0 3 0
% 3 1090 820 1090 920 4 0
% ${\Bbb R}$
\put(10.9000,-9.2000){\makebox(0,0)[rt]{${\Bbb R}$}}%
% STR 2 0 3 0
% 3 5090 3850 5090 3950 4 0
% $X$
\put(50.9000,-39.5000){\makebox(0,0)[rt]{$X$}}%
% STR 2 0 3 0
% 3 5210 2680 5210 2780 4 0
% $w_0$
\put(52.1000,-27.8000){\makebox(0,0)[rt]{$w_0$}}%
% STR 2 0 3 0
% 3 4960 1100 4960 1200 2 0
% ?
\put(49.6000,-12.0000){\makebox(0,0)[lb]{?}}%
% STR 2 0 3 0
% 3 3600 3560 3600 3660 2 0
% ?
\put(36.0000,-36.6000){\makebox(0,0)[lb]{?}}%
% STR 2 0 3 0
% 3 6410 3580 6410 3680 2 0
% ?
\put(64.1000,-36.8000){\makebox(0,0)[lb]{?}}%
\end{picture}%
\hspace{2.5truecm}}

\vspace{0.5truecm}

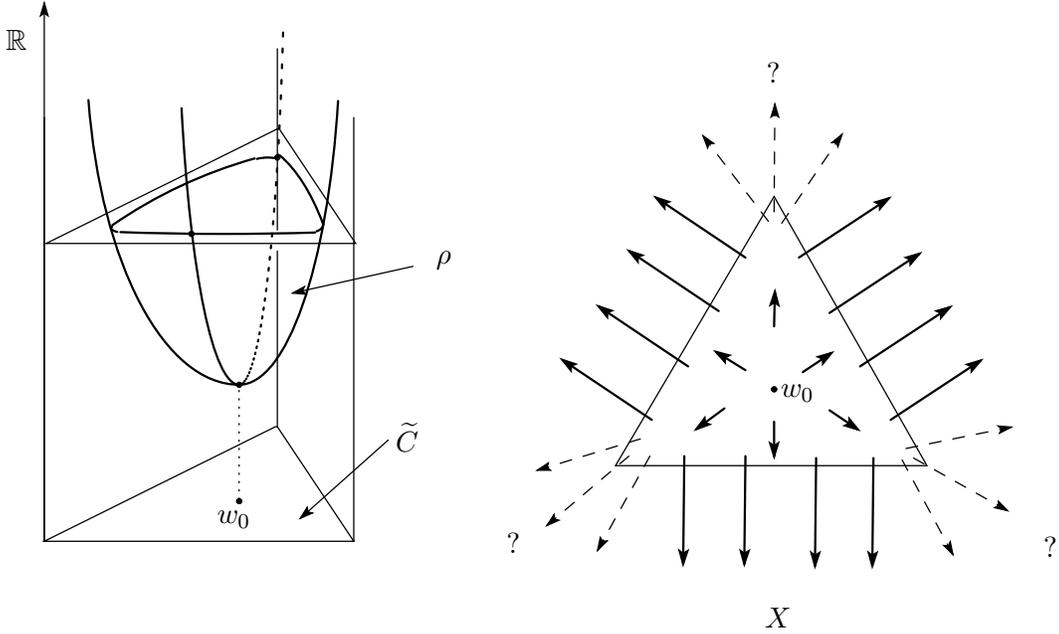
\centerline{{\bf Fig. 2.}}

\vspace{0.5truecm}

\noindent
Let $W$ be the Coxeter group of $\widetilde M$ at $u_0$ and $W_0$ the isotropy 
group of $W$ at ${\bf 0}$.  Also, let $\Gamma$ be the lattice of translations 
of $W$ and set $\Gamma^{\ast}:=\{\omega\in(T^{\perp}_{u_0}\widetilde M)^{\ast}
\,\vert\,\omega(\Gamma)\subset{\bf Z}\}$.  
Also, let $C^{\triangle}(T^{\perp}_{u_0}\widetilde M)^W$ be the space of 
all finite sums of $W$-invariant eigenfunctions of $\triangle$, where 
$\triangle$ is the Laplace operator of $T^{\perp}_{u_0}\widetilde M$.  
Then we can show 
$$\begin{array}{l}
\displaystyle{C^{\triangle}(T^{\perp}_{u_0}\widetilde M)^W\otimes{\Bbb C}=
\{f=\sum_{\omega\in\Gamma^{\ast}}a_{\omega}e^{2\pi\sqrt{-1}\omega}\,\vert\,
a_{\omega}\in{\Bbb C}\,\,(a_{\omega}=0\,\,{\rm except}\,\,{\rm for}\,\,
{\rm finite}\,\,\omega)}\\
\hspace{8.3truecm}\displaystyle{\&\,\,f\,:\,W-{\rm invariant}\}.}
\end{array}$$
Hence, according to Chapter VI,$\S$ 3, Theorem 1 of [B], we have 
$C^{\triangle}(T^{\perp}_{u_0}\widetilde M)^W\otimes{\Bbb C}
={\Bbb C}[\phi_1,\cdots,\phi_r]$ (polynomial ring) for some 
$\phi_1,\cdots,\phi_r(\in C^{\triangle}(T^{\perp}_{u_0}\widetilde M)^W)$.  
Here we note that $r={\rm dim}\,T^{\perp}_{u_0}\widetilde M$ because $G/K$ is 
irreducible and hence $M$ is full and irreducible.  
By reordering $\phi_1,\cdots,\phi_r$ suitably, it is shown that 
$\{{\rm Re}\,\phi_1,\cdots,{\rm Re}\,\phi_{r_1+r_2},\,{\rm Im}\,\phi_{r_1+1},$
\newline
$\cdots,{\rm Im}\,\phi_{r_1+r_2}\}$ is a base of 
$C^{\triangle}(T^{\perp}_{u_0}(\widetilde M))^W$ (see proof of Theorem 7.6 of 
[HLO]), where ${\rm Re}(\cdot)$ (resp. ${\rm Im}(\cdot)$) is the real part 
(resp. the imaginary part) of $(\cdot)$, and $r_1$ and $r_2$ are positive 
integers with $r_1+2r_2=r$.  Denote by $\{\psi_1,\cdots,\psi_r\}$ this base 
for simplicity.  Set $\Psi:=(\psi_1,\cdots,\psi_r)$, which is a 
$C^{\infty}$-map from $T^{\perp}_{u_0}\widetilde M$ onto ${\Bbb R}^r$.  
It is shown that $\Psi$ is injective and that 
$\Psi\vert_{\overline{\widetilde C}}$ is a homeomorphism of 
$\overline{\widetilde C}$ onto $\Psi(\overline{\widetilde C})$, where 
$\overline{\widetilde C}$ is the closure of $\widetilde C$.  
Set $\xi_w(t):=\psi_t(w)$ and $\bar{\xi}_w(t):=\Psi(\psi_t(w))$, where 
$w\in\widetilde C$.  
Also, we set $\xi^i_w(t):=x_i(\xi_w(t))$ and 
$\bar{\xi}^i_w(t):=y_i(\bar{\xi}_w(t))$, where $(y_1,\cdots,y_r)$ is 
the natural coordinate of ${\Bbb R}^r$.  Then we have 
$$\begin{array}{l}
\hspace{0.6truecm}
\displaystyle{(\bar{\xi}^i_w)'(t)=\langle{\rm grad}(y_i\circ\Psi)_{\xi_w(t)},
X_{\xi_w(t)}\rangle}\\
\displaystyle{=\sum_{a=1}^{\bar r}\left(\left(
m_a^e\cot(\frac{\pi}{2b_a}(1-(\lambda_a)_{u_0}(\xi_w(t))))
-m_a^o\tan(\frac{\pi}{2b_a}(1-(\lambda_a)_{u_0}(\xi_w(t))))\right)\right.}\\
\hspace{1.2truecm}
\displaystyle{\left.\times\frac{\pi}{2b_a}(\lambda_a)_{u_0}
({\rm grad}(y_i\circ\Psi)_{\xi_w(t)})\right).}
\end{array}$$
Denote by $W$ the Coxeter group of $\widetilde M$ at $u_0$.  
Let $f_i$ be the $W$-invariant $C^{\infty}$-function over 
$T^{\perp}_{u_0}\widetilde M$ such that 
$$\begin{array}{l}
\displaystyle{f_i(v):=
\sum_{a=1}^{\bar r}\left(\left(
m_a^e\cot(\frac{\pi}{2b_a}(1-(\lambda_a)_{u_0}(v)))
-m_a^o\tan(\frac{\pi}{2b_a}(1-(\lambda_a)_{u_0}(v)))\right)\right.}\\
\hspace{2truecm}
\displaystyle{\left.\times\frac{\pi}{2b_a}(\lambda_a)_{u_0}
({\rm grad}(y_i\circ\Psi)_v)\right)}
\end{array}$$
for all $v\in W\cdot\widetilde C$.  
It is easy to show that such a $W$-invariant $C^{\infty}$-function exists 
uniquely.  Hence, we can describe 
$f_i$ as $f_i=Y_i\circ\Psi$ in terms of some $C^{\infty}$-function 
$Y_i$ over ${\Bbb R}^r$.  Set $Y:=(Y_1,\cdots,Y_r)$, which is regarded as 
a $C^{\infty}$-vector field on ${\Bbb R}^r$.  
Then we have $Y_{\Psi(w)}=\Psi_{\ast}(X_w)$ ($w\in\widetilde C$), that is, 
$Y\vert_{\Psi(\widetilde C)}=\Psi_{\ast}(X)$.  Also we can show that 
$Y\vert_{\partial\Psi(\widetilde C)}$ has no zero point.  
From these facts and the fact that $X$ is as in Fig. 2, we see that 
the flow of $X$ starting any point of $\widetilde C$ 
other than $w_0$ ($w_0:$ the only zero point of $X$) converges to a point of 
$\partial\widetilde C$ in finite time, and that, furthermore, for each point 
of $\partial\widetilde C$, there exists a unique flow of $X$ converging to 
the point.  
Since $M$ is not minimal by the assumption, we have $X_{\bf 0}\not =0$, 
that is, ${\bf 0}\not=w_0$.  Hence we have $T<\infty$ and 
$\lim\limits_{t\to T}\xi(t)\in\partial\widetilde C$, where 
$\xi(t)=\psi_t({\bf 0})$ and $T$ is the supremum of the domain of 
$\xi$.  Set 
$w_1:=\lim\limits_{t\to T-0}\xi(t)$.  Therefore, since 
$M_t=M_{\widetilde{\xi(t)}}$, the mean curvature flow $M_t$ collapses to 
the focal submanifold $F:=M_{\widetilde w_1}$ in time $T$, where 
$\widetilde w_1$ is the parallel normal vector field of $M$ with 
$(\widetilde w_1)_{x_0}=w_1$.  
Also, $M_t$'s ($-S<t<T$) are all of parallel submanifolds of $M$ collapsing to 
$F$ along the mean curvature flow, where $-S$ is the infimum of the domain of 
$\xi$.  Thus the first-half part of the statement (i) and the statement (ii) 
are shown.  
Next we shall show the second-half part of the statement (i).  
Assume that $M$ is irreducible and the codimension of $M$ is greater than 
one and that the fibration of $M$ onto $F$ is spherical.  Set 
$\widetilde M:=(\pi\circ\phi)^{-1}(M)$ and $\widetilde F
:=(\pi\circ\phi)^{-1}(F)$.  Since the fibration of $M$ onto $F$ is 
spherical, $\widetilde F$ passes through a highest dimensional stratum 
$\widetilde{\sigma}$ of $\partial\widetilde C$.  Let $a_0$ be the element of 
$\{1,\cdots,\bar r\}$ 
with $\widetilde{\sigma}\subset(\lambda_{a_0})_{u_0}^{-1}
(1)$.  Set $\widetilde M_t:=(\pi\circ\phi)^{-1}(M_t)$ ($t\in[0,T)$), which 
is the mean curvature flow having $\widetilde M$ as initial data.  
Denote by $A^t$ (resp. $\widetilde A^t$) the shape tensor of $M_t$ (resp. 
$\widetilde M_t$).  Then, since $\widetilde M_t$ is the parallel submanifold 
of $\widetilde M$ for ${\widetilde{\xi(t)}}^L$, we have 
$${\rm Spec}\,\widetilde A^t_v\setminus\{0\}=\{\frac{(\lambda_{aj})_{u_0}(v)}
{1-(\lambda_{aj})_{u_0}(\xi(t))}\,\vert\,a=1,\cdots,\bar r,
\,\,j\in{\Bbb Z}\}$$
for each $v\in T^{\perp}_{u_0+\xi(t)}\widetilde M_t=T^{\perp}_{u_0}
\widetilde M$.  
Since $\lim\limits_{t\to T-0}\xi(t)\in(\lambda_{a_0})_{u_0}^{-1}(1)$ and 
$\lim\limits_{t\to T-0}\xi(t)\notin(\lambda_a)_{u_0}^{-1}(1)$ 
($a\in\{1,\cdots,\bar r\}\setminus\{a_0\}$), we have 
$\lim\limits_{t\to T-0}(\lambda_{a_0})_{u_0}(\xi(t))=1$ and 
$\lim\limits_{t\to T-0}(\lambda_a)_{u_0}(\xi(t))<1$ ($a\not=a_0$).  
From these facts, $\xi'(t)=(\widetilde H^{\widetilde{\xi(t)}^L})_{u_0+\xi(t)}$ 
and $(3.1)$, we have 
$$\begin{array}{l}
\hspace{0.6truecm}
\displaystyle{\lim_{t\to T-0}\vert\vert\widetilde A^t_v\vert\vert_{\infty}^2
(T-t)}\\
\displaystyle{=\lim_{t\to T-0}\frac{(\lambda_{a_0})_{u_0}(v)^2}
{(1-(\lambda_{a_0})_{u_0}(\xi(t)))^2}(T-t)}\\
\displaystyle{=\frac12(\lambda_{a_0})_{u_0}(v)^2
\lim_{t\to T-0}\frac{1}{(1-(\lambda_{a_0})_{u_0}(\xi(t)))(\lambda_{a_0})_{u_0}
(\xi'(t))}}\\
\displaystyle{=\frac{(\lambda_{a_0})_{u_0}(v)^2}
{2m_{a_0}^e\vert\vert({\bf n}_{a_0})_{u_0}\vert\vert^2}.}
\end{array}\leqno{(3.2)}$$
Hence we have 
$$\lim_{t\to T-0}\mathop{\max}_{v\in S^{\perp}_{u_0+\xi(t)}\widetilde M_t}
\vert\vert\widetilde A^t_v\vert\vert_{\infty}^2(T-t)
=\frac{1}{2m_{a_0}^e}.$$
Thus the mean curvature flow $\widetilde M_t$ has type I singularity.  Set 
$\bar v_t:=(\pi\circ\phi)_{\ast u_0+\xi(t)}(v)$ and let 
$\{\lambda^t_1,\cdots,\lambda^t_n\}\,\,(\lambda^t_1\leq\cdots\leq\lambda^t_n)$ 
(resp. $\{\mu^t_1,\cdots,\mu^t_n\}\,\,(0\leq\mu^t_1\leq\cdots\leq\mu^t_n)$) 
be all the eigenvalues of $A^t_{\bar v_t}$ 
(resp. $R(\cdot,\bar v_t)\bar v_t$), where $n:={\rm dim}\,M$.  
Since $M$ is an irreducible equifocal 
submanifold of codimension greater than one by the assumption, it is 
homogeneous by the homogeneity theorem of Christ (see [Ch]) and hence it is 
a principal orbit of a Hermann action by the result of 
Heintze-Palais-Terng-Thorbergsson (see [HPTT]) and the classification of 
hyperpolar actions by Kollross (see [Kol]).  Furthermore, $M$ and its parallel 
submanifolds are curvature-adapted by the result of Goertsches-Thorbergsson 
(see [GT]).  
Therefore, $A^t_{\bar v_t}$ and $R(\cdot,\bar v_t)\bar v_t$ commute and hence 
we have 
$$\sum_{i=1}^n\sum_{j=1}^n\left(
{\rm Ker}(A^t_{\bar v_t}-\lambda^t_i\,{\rm id})\cap{\rm Ker}(R(\cdot,\bar v_t)
\bar v_t-\mu_j^t\,{\rm id})\right)=T_{(\pi\circ\phi)(u_0+\xi(t))}M_t.$$
Set $\bar E_{ij}^t:={\rm Ker}(A^t_{\bar v_t}-\lambda_i^t\,{\rm id})\cap 
{\rm Ker}(R(\cdot,\bar v_t)\bar v_t-\mu_j^t\,{\rm id})$ 
($i,j\in\{1,\cdots,n\}$) and 
$I_t:=\{(i,j)\in\{1,\cdots,n\}^2\,\vert\,\bar E_{ij}^t\not=\{0\}\}$.  
For each $(i,j)\in I_t$, we have 
$${\rm Spec}(\widetilde A^t_v\vert_{(\pi\circ\phi)_{\ast}^{-1}(\bar E_{ij}^t)})
=
\left\{
\begin{array}{ll}
\displaystyle{\left\{\frac{\sqrt{\mu_j^t}}{\arctan\frac{\sqrt{\mu_j^t}}
{\lambda_i^t}+k\pi}\,\vert\,k\in{\Bbb Z}\right\}} & 
\displaystyle{(\mu_j^t\not=0)}\\
\displaystyle{\{\lambda_i^t\}} & \displaystyle{(\mu_j^t=0)}
\end{array}\right.$$
in terms of Proposition 3.2 of [Koi1] and hence 
$$\vert\vert\widetilde A^t_v\vert\vert_{\infty}=\max\left(\left\{
\frac{\sqrt{\mu_j^t}}{\arctan\frac{\sqrt{\mu_j^t}}{\vert\lambda_i^t\vert}}\,
\vert\,(i,j)\in I_t\,\,{\rm s.t.}\,\,\mu_j^t\not=0\right\}\cup
\{\vert\lambda_i^t\vert\,\vert\,(i,j)\in I_t\,\,{\rm s.t.}\,\,\mu_j^t=0\}
\right).$$
It is clear that 
$\displaystyle{\mathop{{\rm sup}}_{0\leq t<T}\mu_n^t\,<\,\infty}$.  
If $\lim\limits_{t\to T-0}\vert\lambda_i^t\vert=\infty$, then we have 
$\displaystyle{\lim_{t\to T-0}\left(\frac{\sqrt{\mu_j^t}}
{\arctan\frac{\sqrt{\mu_j^t}}{\lambda_i^t}}\right)/\lambda_i^t}$\newline
$=1$.  Hence we have 
$$\begin{array}{l}
\hspace{0.6truecm}\displaystyle{\lim\limits_{t\to T-0}
\vert\vert\widetilde A^t_v\vert\vert^2_{\infty}(T-t)}\\
\displaystyle{=\max\left\{\lim\limits_{t\to T-0}
(\lambda_i^t)^2(T-t)\,\vert\,i=1,\cdots,n\right\}}\\
\displaystyle{=\lim\limits_{t\to T-0}\max\{(\lambda_i^t)^2(T-t)\,\vert\,
i=1,\cdots,n\}}\\
\displaystyle{=\lim\limits_{t\to T-0}
\vert\vert A^t_{\bar v_t}\vert\vert^2_{\infty}(T-t),}
\end{array}$$
which together with $(3.2)$ deduces
$$\lim\limits_{t\to T-0}\vert\vert A^t_{\bar v_t}\vert\vert^2_{\infty}(T-t)=
\frac{(\lambda_{a_0})_{u_0}(v)^2}{2m_{a_0}^e\vert\vert
({\bf n}_{a_0})_{u_0}\vert\vert^2}.$$
Therefore we obtain 
$$\lim\limits_{t\to T-0}\mathop{\max}_{v\in S^{\perp}_{\exp^{\perp}(\xi(t))}
M_t}\vert\vert A^t_v\vert\vert^2_{\infty}(T-t)=\frac{1}{2m_{a_0}^e}<\infty.$$
Thus the mean curvature flow $M_t$ has type I singularity.  
\hspace{3.4truecm}q.e.d.

\vspace{0.5truecm}

Next we prove Theorem B.  

\vspace{0.5truecm}

\noindent
{\it Proof of Theorem B.} 
For simplicity, set $I:=\{1,\cdots,\bar r\}$.  
Let $\widetilde{\sigma}$ be a stratum of dimension greater than zero of 
$\partial\widetilde C$ and 
$I_{\widetilde{\sigma}}:=\{a\in I\,\vert\,\widetilde{\sigma}\subset 
(\lambda_a)_{u_0}^{-1}(1)\}$.  Let $w_1\in\widetilde{\sigma}$.  
Denote by $F$ (resp. $\widetilde F$) the focal submanifold 
of $M$ (resp. $\widetilde M$) for $\widetilde w_1$ (resp. 
${\widetilde w_1}^L$).  
Assume that $F$ is not minimal.  Then, since 
$\displaystyle{{\rm Ker}\,\eta_{\widetilde w_{\ast}}
=\mathop{\oplus}_{a\in I_{\widetilde{\sigma}}}(E_{a0})_{u_0}}$, we have 
$$T_{u_0+w_1}\widetilde F
=\left(\mathop{\oplus}_{a\in I\setminus I_{\widetilde{\sigma}}}
\mathop{\oplus}_{j\in{\bf Z}}\eta_{\widetilde w_1\ast}((E_{aj})_{u_0})\right)
\oplus\left(\mathop{\oplus}_{a\in I_{\widetilde{\sigma}}}
\mathop{\oplus}_{j\in{\bf Z}\setminus\{0\}}\eta_{\widetilde w_1\ast}
((E_{aj})_{u_0})
\right).$$
Also we have 
$$T^{\perp}_{u_0+w_1}\widetilde F
=\left(\mathop{\oplus}_{a\in I_{\widetilde{\sigma}}}
(E_{a0})_{u_0}\right)\oplus T^{\perp}_{u_0}\widetilde M,$$
where we identify $T_{u_0+w_1}H^0([0,1],\mathfrak g)$ with 
$T_{u_0}H^0([0,1],\mathfrak g)$.  
For $v\in T^{\perp}_{u_0}\widetilde M\,(\subset 
T^{\perp}_{u_0+w_1}\widetilde F)$, 
we have 
$$\widetilde A^{\widetilde w_1^L}_v\vert_{\eta_{\widetilde w\ast}
((E_{aj})_{u_0})}=\frac{(\lambda_{aj})_{u_0}(v)}{1-(\lambda_{aj})_{u_0}(w_1)}
{\rm id}\,\,\,\,
((a,j)\in((I\setminus I_{\widetilde{\sigma}})\times{\Bbb Z})\cup
(I_{\widetilde{\sigma}}\times({\Bbb Z}\setminus\{0\}))).$$
Hence we have 
$$\begin{array}{l}
\hspace{0.6truecm}
\displaystyle{{\rm Tr}_r{\widetilde A}^{{\widetilde w}_1^L}_v}\\
\displaystyle{=\sum_{a\in I\setminus I_{\widetilde{\sigma}}}\left(
\sum_{j\in{\bf Z}}\frac{m_a^e(\lambda_{a,2j})_{u_0}(v)}
{1-(\lambda_{a,2j})_{u_0}(w_1)}+\sum_{j\in{\bf Z}}
\frac{m_a^o(\lambda_{a,2j+1})_{u_0}(v)}
{1-(\lambda_{a,2j+1})_{u_0}(w_1)}\right)}\\
\hspace{0.6truecm}\displaystyle{+\sum_{a\in I_{\widetilde{\sigma}}}
\left(\sum_{j\in{\bf Z}\setminus\{0\}}\frac{m_a^e(\lambda_{a,2j})_{u_0}(v)}
{1-(\lambda_{a,2j})_{u_0}(w_1)}+\sum_{j\in{\bf Z}}
\frac{m_a^o(\lambda_{a,2j+1})_{u_0}(v)}
{1-(\lambda_{a,2j+1})_{u_0}(w_1)}\right)}\\
%%\displaystyle{=\sum_{a\in I\setminus I_{\widetilde{\sigma}}}\left(
%%(m_a^e+m_a^o)\cot(\frac{\pi(1-(\lambda_a)_{u_0}(w))}{b_a})\right.}\\
%%\hspace{2truecm}\displaystyle{
%%\left.+(m_a^e-m_a^o){\rm cosec}(\frac{\pi(1-(\lambda_a)_{u_0}(w))}{b_a})\right%%)
%%\frac{\pi}{2b_a}(\lambda_a)_{u_0}(v)}\\
%%\hspace{0.6truecm}\displaystyle{+\sum_{a\in I_{\widetilde{\sigma}}}
%%\left(\frac{m_a^e(\lambda_a)_{u_0}(v)}{b_a}
%%\sum_{j\in{\bf Z}\setminus\{0\}}\frac{1}{2j}+
%%\frac{m_a^o(\lambda_a)_{u_0}(v)}{b_a}
%%\sum_{j\in{\bf Z}}\frac{1}{2j+1}\right)}\\
\displaystyle{=\sum_{a\in I\setminus I_{\widetilde{\sigma}}}\left(
m_a^e\cot\frac{\pi}{2b_a}(1-(\lambda_a)_{u_0}(w_1))
-m_a^o\tan\frac{\pi}{2b_a}(1-(\lambda_a)_{u_0}(w_1))\right)\frac{\pi}{2b_a}
(\lambda_a)_{u_0}(v),}
\end{array}$$
that is, the $T_{u_0}^{\perp}\widetilde M$-component 
$((\widetilde H^{{\widetilde w}_1^L})_{u_0+w_1})_{T^{\perp}_{u_0}
\widetilde M}$ 
of $(\widetilde H^{{\widetilde w}_1^L})_{u_0+w_1}$ is equal to 
$$\sum_{a\in I\setminus I_{\widetilde{\sigma}}}\left(
m_a^e\cot\frac{\pi}{2b_a}(1-(\lambda_a)_{u_0}(w_1))
-m_a^o\tan\frac{\pi}{2b_a}(1-(\lambda_a)_{u_0}(w_1))\right)\frac{\pi}{2b_a}
({\bf n}_a)_{u_0}.$$
Denote by $\Phi_{u_0+w_1}$ the normal 
holonomy group of $\widetilde F$ at $u_0+w_1$ and $L_{u_0}$ be the focal leaf 
through $u_0$ for $\widetilde w_1$.  Since $L_{u_0}=\Phi_{u_0+w_1}\cdot u_0$, 
there exists $\mu\in\Phi_{u_0+w_1}$ such that 
$\mu(T^{\perp}_{u_0}\widetilde M)=T^{\perp}_{u_1}\widetilde M$ for any point 
$u_1$ of $L_{u_0}$.  
On the other hand, since $\widetilde F$ has constant principal curvatures 
in the sense of [HOT], $({\widetilde H}^{{\widetilde w}_1^L})_{u_0+w_1}$ is 
$\Phi_{u_0+w_1}$-invariant.  
Hence we have 
$\displaystyle{(\widetilde H^{{\widetilde w}_1^L})_{u_0+w_1}
\in\mathop{\cap}_{u\in L_{u_0}}T^{\perp}_u\widetilde M}$, 
where we note that 
$\displaystyle{\mathop{\cap}_{u\in L_{u_0}}T^{\perp}_u\widetilde M}$ contains 
$\widetilde{\sigma}$ as an open subset.  Therefore, we obtain 
$$\begin{array}{l}
\displaystyle{
(\widetilde H^{{\widetilde w}_1^L})_{u_0+w_1}=
\sum_{a\in I\setminus I_{\widetilde{\sigma}}}\left(
m_a^e\cot\frac{\pi}{2b_a}(1-(\lambda_a)_{u_0}(w_1))\right.}\\
\hspace{2.4truecm}\displaystyle{\left.
-m_a^o\tan\frac{\pi}{2b_a}(1-(\lambda_a)_{u_0}(w_1))\right)\frac{\pi}{2b_a}
({\bf n}_a)_{u_0}\,\,\,\,(\in T\widetilde{\sigma}).}
\end{array}\leqno{(3.3)}$$
Define a tangent vector field $X^{\widetilde{\sigma}}$ on 
$\widetilde{\sigma}$ by $X^{\widetilde{\sigma}}_w
:=(\widetilde H^{{\widetilde w}^L})_{u_0+w}$ ($w\in\widetilde{\sigma}$).  
Let $\xi:(-S,T)\to\widetilde{\sigma}$ be the maximal integral curve of 
$X^{\widetilde{\sigma}}$ with $\xi(0)=w_1$.  
Define a function $\rho_{\widetilde{\sigma}}$ over $\widetilde{\sigma}$ by 
$$\begin{array}{l}
\displaystyle{\rho_{\widetilde{\sigma}}(w)
:=-\sum_{a\in I\setminus I_{\widetilde{\sigma}}}
\left(m_a^e\log\sin\frac{\pi}{2b_a}(1-(\lambda_a)_{u_0}(w))\right.}\\
\hspace{2.4truecm}\displaystyle{
\left.+m_a^o\log\cos\frac{\pi}{2b_a}(1-(\lambda_a)_{u_0}(w))\right)
\qquad(w\in\widetilde{\sigma}).}
\end{array}$$
It follows from the definition of $X^{\widetilde{\sigma}}$ and $(3.3)$ that 
${\rm grad}\,\rho_{\widetilde{\sigma}}
=X^{\widetilde{\sigma}}$.  Also we can show that 
$\rho_{\widetilde{\sigma}}$ is downward convex and that 
$\rho_{\widetilde{\sigma}}(w)\to\infty$ as $w\to\partial\widetilde{\sigma}$.  
Hence we see that $\rho_{\widetilde{\sigma}}$ has the only minimal point.  
Denote by $w_0$ this minimal point.  It is clear that 
$X^{\widetilde{\sigma}}_{w_0}=0$.  
Also, by imitating the proof of Theorem A, we can show that the flow of 
$X^{\widetilde{\sigma}}$ starting any point of $\widetilde{\sigma}$ other than 
$w_0$ converges to a point of $\partial\widetilde{\sigma}$ in finite time, and 
that, furthermore, for each point of $\partial\widetilde{\sigma}$, there 
exists a unique flow of $X^{\widetilde{\sigma}}$ converging to the point.  
Since $F$ is not minimal by the assumption, we have 
$X^{\widetilde{\sigma}}_{w_1}\not=0$, that is, $w_1\not=w_0$.  Hence we have 
$T<\infty$ and $\lim\limits_{t\to T-0}\xi(t)\in\partial\widetilde{\sigma}$.  
Set $w_2:=\lim\limits_{t\to T-0}\xi(t)$.  Therefore, since 
$F_t=M_{\widetilde{\xi(t)}}$, the mean curvature flow $F_t$ collapses to the 
lower dimensional focal submanifold $F':=M_{\widetilde w_2}$ in time $T$, 
where $\widetilde w_2$ is the parallel normal vector field of $M$ with 
$(\widetilde w_2)_{x_0}=w_2$.  
Also, $F_t$'s ($-S<t<T$) are all of focal submanifolds of $M$ through 
$\widetilde{\sigma}$ collapsing to $F'$ along the mean curvature flow.  
Thus the first-half part of the statement (i) and the statement (ii) are 
shown.  Also, by imitating the proof of Theorem A, we can show 
the second-half part of the statement (i).  \hspace{10.5truecm}q.e.d.

\section{Hermann actions of cohomogeneity two}
According to the homogeneity theorem for an equifocal submanifold 
in a symmetric space of compact type by Christ ([Ch]), 
equifocal submanifolds of codimension greater than one 
in an irreducible compact type symmetric space are homogeneous.  
Hence, according to the result by Heintze-Palais-Terng-Thorbergsson 
([HPTT]), they occur as principal orbits of hyperpolar actions on 
the symmetric space.  Furthermore, by using the classification of hyperpolar 
actions on irreducible compact type symmetric spaces by Kollross ([Kol]), 
we see that they occur as principal orbits of Hermann actions on 
the symmetric spaces.  
We have only to analyze the vector field $X$ defined in the previous secton 
to analyze the mean curvature flows having parallel submanifolds of an 
equifocal submanifold $M$ as initial data.  
Also, we have only to analyze the vector fields $X^{\widetilde{\sigma}}$'s 
($\widetilde{\sigma}\,:\,$a simplex of $\partial\widetilde C$) defined in 
the proof of Theorem B to analyze the mean curvature flows having focal 
submanifolds of $M$ as initial data.  
In this section, we shall explicitly describe the vector field $X$ 
defined for principal orbits of all Hermann actions of cohomogeneity two on 
all irreducible symmetric spaces of compact type and rank two (see Table 3).  
Let $G/K$ be a symmetric space of compact type and $H$ be 
a symmetric subgroup of $G$.  Also, let $\theta$ be an involution of 
$G$ with $({\rm Fix}\,\theta)_0\subset K\subset{\rm Fix}\,\theta$ and $\tau$ 
be an invloution of $G$ with $({\rm Fix}\,\tau)_0\subset H\subset{\rm Fix}\,
\tau$, where ${\rm Fix}\,\theta$ (resp. ${\rm Fix}\,\tau$) is the fixed point 
group of $\theta$ (resp. $\tau$) and $({\rm Fix}\,\theta)_0$ (resp. 
$({\rm Fix}\,\tau)_0$) is the identity component of ${\rm Fix}\,\theta$ 
(resp. ${\rm Fix}\,\tau$).  In the sequel, we assume that 
$\tau\circ\theta=\theta\circ\tau$.  Set $L:={\rm Fix}(\theta\circ\tau)$.  
Denote by the same symbol $\theta$ (resp. $\tau$) the involution of the Lie 
algebra $\mathfrak g$ of $G$ induced from $\theta$ (resp. $\tau$).  Set 
$\mathfrak k:={\rm Ker}(\theta-{\rm id}),\,\mathfrak p:={\rm Ker}(\theta+
{\rm id}),\,\mathfrak h:={\rm Ker}(\tau-{\rm id})$ and $\mathfrak q:=
{\rm Ker}(\tau+{\rm id})$.  The space $\mathfrak p$ is identified with 
$T_{eK}(G/K)$.  From $\theta\circ\tau=\tau\circ\theta$, we have 
$\mathfrak p=\mathfrak p\cap\mathfrak h+\mathfrak p\cap\mathfrak q$.  Take 
a maximal abelian subspace $\mathfrak b$ of $\mathfrak p\cap\mathfrak q$ and 
let $\mathfrak p=\mathfrak z_{\mathfrak p}(\mathfrak b)
+\sum\limits_{\beta\in\triangle'_+}\mathfrak p_{\beta}$ be the root space 
decomposition with respect to $\mathfrak b$, where 
$\mathfrak z_{\mathfrak p}(\mathfrak b)$ is the centralizer of $\mathfrak b$ 
in $\mathfrak p$, $\triangle'_+$ is the positive root system of 
$\triangle':=\{\beta\in\mathfrak b^{\ast}\,\vert\,\exists\,X(\not=0)\in
\mathfrak p\,\,{\rm s.t.}\,\,{\rm ad}(b)^2(X)=-\beta(b)^2X\,\,(\forall\,b\in
\mathfrak b)\}$ under some lexicographic ordering of $\mathfrak b^{\ast}$ and 
$\mathfrak p_{\beta}:=\{X\in\mathfrak p\,\vert\,{\rm ad}(b)^2(X)=-\beta(b)^2X
\,\,(\forall\,b\in\mathfrak b)\}$ ($\beta\in\triangle'_+$). Also, let 
${\triangle'}^V_+:=\{\beta\in\triangle'_+\,\vert\,\mathfrak p_{\beta}\cap
\mathfrak q\not=\{0\}\}$ and ${\triangle'}^H_+:=\{\beta\in\triangle'_+\,\vert\,
\mathfrak p_{\beta}\cap\mathfrak h\not=\{0\}\}$.  
Then we have $\mathfrak q=\mathfrak b+\sum\limits_{\beta\in{\triangle'}^V_+}
(\mathfrak p_{\beta}\cap\mathfrak q)$ and $\mathfrak h=
\mathfrak z_{\mathfrak h}(\mathfrak b)+\sum\limits_{\beta\in{\triangle'}^H_+}
(\mathfrak p_{\beta}\cap\mathfrak h)$, where $\mathfrak z_{\mathfrak h}
(\mathfrak b)$ is the centralizer of $\mathfrak b$ in $\mathfrak h$.  The 
orbit $H(eK)$ is a reflective submanifold and it is isometric to the symmetric 
space $H/H\cap K$ (equipped with a metric scaled suitably).  Also, 
$\exp^{\perp}(T^{\perp}_{eK}(H(eK)))$ is also a reflective submanifold and it 
is isometric to the symmetric space $L/H\cap K$ (equipped with a metric scaled 
suitably), where $\exp^{\perp}$ is the normal exponential map of $H(eK)$.  
The system ${\triangle'}^V:={\triangle'}^V_+\cup(-{\triangle'}^V_+)$ is the 
root system of $L/H\cap K$.  Define a subset $\widetilde C$ of $\mathfrak b$ 
by 
$$\begin{array}{l}
\displaystyle{\widetilde C:=\{b\in\mathfrak b\,\vert\,0<\beta<\pi\,(\forall\,
\beta\in{\triangle'}^V_+),\,\,-\frac{\pi}{2}<\beta<
\frac{\pi}{2}\,(\forall\,\beta\in{\triangle'}^H_+)\}.}
\end{array}$$
%Let $\Pi$ be the simple root system of $\triangle'_+$, and set $\Pi_V:=\Pi\cap
%{\triangle'}^V_+$ and $\Pi_H:=\Pi\cap{\triangle'}^H_+$.  Also, let $\delta$ 
%be the highest root of ${\triangle'}^V_+\cup2{\triangle'}^H_+$.  
%Then we have 
%$$\widetilde C=\{b\in\mathfrak b\,\vert\,\beta(b)>0\,(\forall\,\beta\in\Pi_V),
%\,\,\beta(b)>-\frac{\pi}{2}\,(\forall\,\beta\in\Pi_H\setminus\Pi_V),\,\,
%\delta(b)<\pi\}.$$
Set $C:={\rm Exp}(\widetilde C)$, where ${\rm Exp}$ is the exponential map 
of $G/K$ at $eK$.  
Let $P(G,H\times K):=\{g\in H^1([0,1],G)\,\vert\,(g(0),g(1))\in H\times K\}$, 
where $H^1([0,1],G)$ is the Hilbert Lie group of all $H^1$-paths in $G$.  
This group acts on $H^0([0,1],\mathfrak g)$ as gauge action.  The orbits of 
the $P(G,H\times K)$-action are the inverse images of orbits of the $H$-action 
by $\pi\circ\phi$.  The set $\Sigma:={\rm Exp}(\mathfrak b)$ is a section of 
the $H$-action and $\mathfrak b$ is a section of the $P(G,H\times K)$-action 
on $H^0([0,1],\mathfrak g)$, where $\mathfrak b$ is identified with the 
horizontal lift of $\mathfrak b$ to the zero element $\hat 0$ of 
$H^0([0,1],\mathfrak g)$ ($\hat0\,:\,$the constant path at the zero element 
$0$ of $\mathfrak g$).  
The set $\widetilde C$ is the fundamental domain of the Coxeter group of 
a principal $P(G,H\times K)$-orbit and 
each prinicipal $H$-orbit meets $C$ at one point and each singular $H$-orbit 
meets $\partial C$ at one point.  
The focal set of the principal orbit $P(G,H\times K)\cdot Z_0$ 
($Z_0\in\widetilde C$) consists of the hyperplanes 
$\beta^{-1}(j\pi)$'s 
($\beta\in{\triangle'}^V_+\setminus{\triangle'}^H_+,\,\,j\in{\Bbb Z}$), 
$\beta^{-1}((j+\frac12)\pi)$'s 
($\beta\in{\triangle'}^H_+\setminus{\triangle'}^V_+,\,\,j\in{\Bbb Z}$), 
$\beta^{-1}(\frac{j\pi}{2})$'s 
($\beta\in{\triangle'}^V_+\cap{\triangle'}^H_+,\,\,j\in{\Bbb Z}$) in 
$\mathfrak b(=T^{\perp}_{Z_0}(P(G,H\times K)\cdot Z_0))$.  
%%Since 
%%the reflections with respect to these focal hyperplanes preserves the sum of 
%%these focal hyperplanes (see ?? of [T-PF]), we see that, if 
%%${\triangle'}V_+\cap{\triangle'}^H_+\not=\emptyset$, then 
%%${\triangle'}V_+={\triangle'}^H_(=\triangle'_+)$ holds.  
Denote by $\exp^G$ the exponential map of $G$.  Note that 
$\pi\circ\exp^G\vert_{\mathfrak p}={\rm Exp}$.  
Let $Y_0\in\widetilde C$ and $M(Y_0):=H({\rm Exp}(Y_0))$.  
Then we have $T^{\perp}_{{\rm Exp}(Y_0)}M(Y_0)
=(\exp^G(Y_0))_{\ast}(\mathfrak b)$.  Denote by $A^{Y_0}$ the shape tensor of 
$M(Y_0)$.  Take $v\in T^{\perp}_{{\rm Exp}(Y_0)}M(Y_0)$ and set 
$\bar v:=(\exp^G(Y_0))_{\ast}^{-1}(v)$.  By scaling the metric of $G/K$ by 
a suitable positive constant, we have 
$$A^{Y_0}_v\vert_{\exp^G(Y_0)_{\ast}(\mathfrak p_{\beta}\cap\mathfrak q)}=
-\frac{\beta(\bar v)}{\tan\beta(Y_0)}{\rm id}\,\,\,\,
(\beta\in{\triangle'}^V_+)\leqno{(4.1)}$$
and 
$$A^{Y_0}_v\vert_{\exp^G(Y_0)_{\ast}(\mathfrak p_{\beta}\cap\mathfrak h)}=
\beta(\bar v)\tan\beta(Y_0){\rm id}\,\,\,\,(\beta\in{\triangle'}^H_+).
\leqno{(4.2)}$$
Set $m_{\beta}^V:={\rm dim}(\mathfrak p_{\beta}\cap\mathfrak q)$ 
($\beta\in{\triangle'}^V_+$) and 
$m_{\beta}^H:={\rm dim}(\mathfrak p_{\beta}\cap\mathfrak h)$ 
($\beta\in{\triangle'}^H_+$).  
Set 
$\widetilde M(Y_0):=(\pi\circ\phi)^{-1}(M(Y_0))(=P(G,H\times K)\cdot Y_0)$.  
We can show $(\pi\circ\phi)(Y_0)={\rm Exp}(Y_0)$.  
Denote by ${\widetilde A}^{Y_0}$ the shape tensor of $\widetilde M(Y_0)$.  
According to Proposition 3.2 of [Koi1], we have 
$$\begin{array}{l}
\displaystyle{{\rm Spec}({\widetilde A}^{Y_0}_{\bar v}
\vert_{(\pi\circ\phi)^{-1}_{\ast Y_0}(\exp^G(Y_0)_{\ast}
(\mathfrak p_{\beta}\cap\mathfrak q))})\setminus\{0\}=
\{\frac{-\beta(\bar v)}{\beta(Y_0)+j\pi}\,\vert\,j\in{\Bbb Z}\}\,\,\,\,
(\beta\in{\triangle'}^V_+),}\\
\displaystyle{{\rm Spec}({\widetilde A}^{Y_0}_{\bar v}
\vert_{(\pi\circ\phi)^{-1}_{\ast Y_0}(\exp^G(Y_0)_{\ast}
(\mathfrak p_{\beta}\cap\mathfrak h))})\setminus\{0\}=
\{\frac{-\beta(\bar v)}{\beta(Y_0)+(j+\frac12)\pi}\,\vert\,j\in{\Bbb Z}\}\,\,\,
\,(\beta\in{\triangle'}^H_+),}
\end{array}$$
and 
$${\rm Spec}({\widetilde A}^{Y_0}_{\bar v}
\vert_{(\pi\circ\phi)^{-1}_{\ast Y_0}(\exp^G(Y_0)_{\ast}
(\mathfrak z_{\mathfrak h}(\mathfrak b)))})=\{0\}.$$
Hence the set ${\cal PC}_{\widetilde M(Y_0)}$ of all principal curvatures of 
$\widetilde M(Y_0)$ is given by 
$${\cal PC}_{\widetilde M(Y_0)}=\{\frac{-\widetilde{\beta}}{\beta(Y_0)+j\pi}\,
\vert\,\beta\in{\triangle'}^V_+,\,\,j\in\Bbb Z\}\cup
\{\frac{-\widetilde{\beta}}{\beta(Y_0)+(j+\frac12)\pi}\,\vert\,
\beta\in{\triangle'}^H_+,\,\,j\in\Bbb Z\},$$
where $\widetilde{\beta}$ is the parallel section of $(T^{\perp}\widetilde M
(Y_0))^{\ast}$ with $\widetilde{\beta}_{u_0}
=\beta\circ\exp^G(Y_0)_{\ast}^{-1}$.  Also, we can show that the multiplicity 
of $\frac{-\widetilde{\beta}}
{\beta(Y_0)+j\pi}$ ($\beta\in{\triangle'}^V_+$) is equal to $m_{\beta}^V$ 
and that of 
$\frac{-\widetilde{\beta}}{\beta(Y_0)+(j+\frac12)\pi}$ 
($\beta\in{\triangle'}^H_+$) is equal to $m_{\beta}^H$.  
Define $\lambda_{\beta}^{Y_0}$ and $b_{\beta}^{Y_0}$ 
($\beta\in{\triangle'}_+$) by 
$$(\lambda^{Y_0}_{\beta},b^{Y_0}_{\beta}):=\left\{
\begin{array}{ll}
\displaystyle{(\frac{-\widetilde{\beta}}{\beta(Y_0)},\,\frac{\pi}
{\beta(Y_0)})} & \displaystyle{(\beta\in{\triangle'}^V_+\setminus
{\triangle'}^H_+)}\\
\displaystyle{(\frac{-\widetilde{\beta}}{\beta(Y_0)+\frac{\pi}{2}},\,
\frac{\pi}{\beta(Y_0)+\frac{\pi}{2}})} & \displaystyle{
(\beta\in{\triangle'}^H_+\setminus{\triangle'}^V_+)}\\
\displaystyle{(\frac{-\widetilde{\beta}}{\beta(Y_0)},\,\frac{\pi}
{2\beta(Y_0)})} & 
\displaystyle{(\beta\in{\triangle'}^V_+\cap{\triangle'}^H_+).}\end{array}
\right.$$
Then we have $\frac{-\widetilde{\beta}}{\beta(Y_0)+j\pi}
=\frac{\lambda_{\beta}^{Y_0}}{1+jb_{\beta}^{Y_0}}$ when 
$\beta\in{\triangle'}^V_+\setminus{\triangle'}^H_+$, 
$\frac{-\widetilde{\beta}}{\beta(Y_0)+(j+\frac12)\pi}
=\frac{\lambda_{\beta}^{Y_0}}{1+jb_{\beta}^{Y_0}}$ when 
$\beta\in{\triangle'}^H_+\setminus{\triangle'}^V_+$ and 
$(\frac{-\widetilde{\beta}}{\beta(Y_0)+j\pi},\frac{-\widetilde{\beta}}
{\beta(Y_0)+(j+\frac12)\pi})
=(\frac{\lambda_{\beta}^{Y_0}}{1+2jb_{\beta}^{Y_0}},\,
\frac{\lambda_{\beta}^{Y_0}}{1+(2j+1)b_{\beta}^{Y_0}})$ when 
$\beta\in{\triangle'}^V_+\cap{\triangle'}^H_+$.  That is, we have 
$${\cal PC}_{\widetilde M(Y_0)}=\{\frac{\lambda_{\beta}^{Y_0}}
{1+jb_{\beta}^{Y_0}}\,\vert\,\beta\in\triangle'_+,\,j\in{\Bbb Z}\}.$$
Denote by $m_{\beta j}$ the multiplicity of 
$\frac{\lambda_{\beta}}{1+jb_{\beta}}$.  Then we have 
$$
%%\begin{array}{l}
m_{\beta,2j}=\left\{
\begin{array}{ll}
\displaystyle{m_{\beta}^V} & \displaystyle{(\beta\in{\triangle'}^V_+\setminus
{\triangle'}^H_+)}\\
\displaystyle{m_{\beta}^H} & \displaystyle{(\beta\in{\triangle'}^H_+
\setminus{\triangle'}^V_+)}\\
\displaystyle{m_{\beta}^V} & \displaystyle{(\beta\in{\triangle'}^V_+\cap
{\triangle'}^H_+),}
\end{array}\right.\qquad
m_{\beta,2j+1}=\left\{
\begin{array}{ll}
\displaystyle{m_{\beta}^V} & \displaystyle{(\beta\in{\triangle'}^V_+
\setminus{\triangle'}^H_+)}\\
\displaystyle{m_{\beta}^H} & \displaystyle{(\beta\in{\triangle'}^H_+\setminus
{\triangle'}^V_+)}\\
\displaystyle{m_{\beta}^H} & \displaystyle{(\beta\in{\triangle'}^V_+\cap
{\triangle'}^H_+),}
\end{array}\right.$$
where $j\in{\Bbb Z}$.  
Denote by $\widetilde H^{Y_0}$ the mean curvature vector of 
$\widetilde M(Y_0)$ and ${\bf n}^{Y_0}_{\beta}$ the curvature normal 
corresponding to $\lambda^{Y_0}_{\beta}$.  
Define $\beta^{\sharp}\,(\in\mathfrak b)$ by $\beta(\cdot)=\langle
\beta^{\sharp},\cdot\rangle$ and let $\widetilde{\beta^{\sharp}}^{Y_0}$ be 
the parallel normal vector field of $\widetilde M(Y_0)$ with 
$(\widetilde{\beta^{\sharp}}^{Y_0})_{Y_0}=\beta^{\sharp}$, where we identify 
$\mathfrak b$ with $T^{\perp}_{Y_0}\widetilde M(Y_0)$.  
From $(3.1)$ (the case of $w=0$), 
we have 
$$\begin{array}{l}
\displaystyle{\widetilde H^{Y_0}=\sum_{\beta\in{\triangle'}^V_+}
m_{\beta}^V\cot\frac{\pi}{2b^{Y_0}_{\beta}}\cdot\frac{\pi}{2b^{Y_0}_{\beta}}
{\bf n}^{Y_0}_{\beta}}\\
\hspace{1.2truecm}\displaystyle{
-\sum_{\beta\in{\triangle'}^H_+}m_{\beta}^H
\tan\frac{\pi}{2b^{Y_0}_{\beta}}\cdot\frac{\pi}{2b_{\beta}^{Y_0}}
{\bf n}_{\beta}^{Y_0}}\\
\hspace{0.6truecm}\displaystyle{=
-\sum_{\beta\in{\triangle'}^V_+}m_{\beta}^V\cot\beta(Y_0)
{\widetilde{\beta^{\sharp}}}^{Y_0}+\sum_{\beta\in{\triangle'}^H_+}
m_{\beta}^H\tan\beta(Y_0){\widetilde{\beta^{\sharp}}}^{Y_0}.}
\end{array}
\leqno{(4.3)}$$
Define a tangent vector field $X$ on $\widetilde C$ by assigning 
$({\widetilde H}^{Y_0})_{Y_0}\,(\in T^{\perp}_{Y_0}\widetilde M(Y_0)
=\mathfrak b(\subset V))$ to each $Y_0\in\widetilde C$.  
From $(4.3)$, we have 
$$X_{Y_0}=-\sum_{\beta\in{\triangle'}^V_+}m_{\beta}^V\cot\beta(Y_0)
\beta^{\sharp}
+\sum_{\beta\in{\triangle'}^H_+}m_{\beta}^H\tan\beta(Y_0)\beta^{\sharp}.
\leqno{(4.4)}$$
By using this description, we can explicitly describe this vector field $X$ 
for all Hermann actions of cohomogeneity two on all irreducible symmetric 
spaces of compact type and rank two.  
All Hermann actions of cohomogeneity two on all irreducible symmetric spaces 
of compact type and rank two are given in Table 1.  The systems 
${\triangle'}^V_+$ and ${\triangle'}^H_+$ for the Hermann actions are given in 
Table 2 and the explicit descriptions of $X$ for the Hermann actions are given 
in Table 3.  In Table 1, $H^{\ast}\curvearrowright G^{\ast}/K$ implies 
the dual action of $H\curvearrowright G/K$ and $L^{\ast}/H\cap K$ is the dual 
of $L/H\cap K$.  
In Table 2, $\{\alpha,\beta,\alpha+\beta\}$ implies a positive root system 
of the root system of (${\mathfrak a}_2$)-type 
($\alpha=(2,0),\beta=(-1,\sqrt3)$), 
$\{\alpha,\beta,\alpha+\beta, 2\alpha+\beta\}$ implies a positive root system 
of the root system of (${\mathfrak b}_2$)(=(${\mathfrak c}_2$))-type 
($\alpha=(1,0),\beta=(-1,1)$) and 
$\{\alpha,\beta,\alpha+\beta,\alpha+2\beta,\alpha+3\beta,2\alpha+3\beta\}$ 
implies a positive root system 
of the root system of (${\mathfrak g}_2$)-type 
($\alpha=(2\sqrt3,0),\beta=(-\sqrt3,1)$).  
In Table $1\sim3$, $\rho_i$ ($i=1,\cdots,16$) imply 
automorphisms of $G$ and $(\cdot)^2$ implies the product Lie group 
$(\cdot)\times(\cdot)$ of a Lie group $(\cdot)$.  In Table 2, 
$\displaystyle{\mathop{\alpha}_{(m)}}$ and so on imply that 
the multiplicity of $\alpha$ is equal to $m$.  

\vspace{0.5truecm}

$$\begin{tabular}{|c|c|c|}
\hline
{\scriptsize$H\curvearrowright G/K$} & {\scriptsize$H^{\ast}\curvearrowright 
G^{\ast}/K$} & {\scriptsize$L^{\ast}/H\cap K$}\\
\hline
{\scriptsize$\rho_1(SO(3))\curvearrowright SU(3)/SO(3)$} & 
{\scriptsize$SO_0(1,2)\curvearrowright SL(3,{\Bbb R})/SO(3)$} & 
{\scriptsize$(SL(2,{\Bbb R})/SO(2))\times{\Bbb R}$}\\
\hline
{\scriptsize$SO(6)\curvearrowright SU(6)/Sp(3)$} & 
{\scriptsize$SO^{\ast}(6)\curvearrowright SU^{\ast}(6)/Sp(3)$} & 
{\scriptsize$SL(3,{\Bbb C})/SU(3)$}\\
\hline
{\scriptsize$\rho_2(Sp(3))\curvearrowright SU(6)/Sp(3)$} & 
{\scriptsize$Sp(1,2)\curvearrowright SU^{\ast}(6)/Sp(3)$} & 
{\scriptsize$(SU^{\ast}(4)/Sp(2))\times U(1)$}\\
\hline
{\scriptsize$SO(q+2)\curvearrowright$} 
& {\scriptsize$SO_0(2,q)\curvearrowright$} & 
{\scriptsize$SO_0(2,q)/SO(2)\times SO(q)$}\\
{\scriptsize$SU(q+2)/S(U(2)\times U(q))$} 
& {\scriptsize$SU(2,q)/S(U(2)\times U(q))$} & 
{\scriptsize}\\
\hline
{\scriptsize$S(U(j+1)\times U(q-j+1))\curvearrowright$} 
& {\scriptsize$S(U(1,j)\times U(1,q-j))\curvearrowright$} & 
{\scriptsize$(SU(1,j)/S(U(1)\times U(j)))\times$}\\
{\scriptsize$SU(q+2)/S(U(2)\times U(q))$} 
& {\scriptsize$SU(2,q)/S(U(2)\times U(q))$} & 
{\scriptsize$(SU(1,q-j)/S(U(1)\times U(q-j)))$}\\
\hline
{\scriptsize$SO(j+1)\times SO(q-j+1)\curvearrowright$} 
& {\scriptsize$SO(1,j)\times SO(1,q-j)\curvearrowright$} & 
{\scriptsize$(SO_0(1,j)/SO(j))\times$}\\
{\scriptsize$SO(q+2)/SO(2)\times SO(q)$} 
& {\scriptsize$SO(2,q)/SO(2)\times SO(q)$} & 
{\scriptsize$(SO_0(1,q-j)/SO(q-j))$}\\
\hline
{\scriptsize$SO(4)\times SO(4)\curvearrowright$} 
& {\scriptsize$SO^{\ast}(4)\times SO^{\ast}(4)\curvearrowright$} & 
{\scriptsize$SU(2,2)/S(U(2)\times U(2))$}\\
{\scriptsize$SO(8)/U(4)$} 
& {\scriptsize$SO^{\ast}(8)/U(4)$} & 
{\scriptsize}\\
\hline
{\scriptsize$\rho_3(SO(4)\times SO(4))\curvearrowright$} 
& {\scriptsize$SO(4,{\Bbb C})\curvearrowright SO^{\ast}(8)/U(4)$} & 
{\scriptsize$SO(4,{\Bbb C})/SO(4)$}\\
{\scriptsize$SO(8)/U(4)$} 
& {\scriptsize} & {\scriptsize}\\
\hline
{\scriptsize$\rho_4(U(4))\curvearrowright SO(8)/U(4)$} 
& {\scriptsize$U(2,2)\curvearrowright SO^{\ast}(8)/U(4)$} & 
{\scriptsize$(SO^{\ast}(4)/U(2))\times(SO^{\ast}(4)/U(2))$}\\
\hline
{\scriptsize$SO(4)\times SO(6)\curvearrowright$} 
& {\scriptsize$SO^{\ast}(4)\times SO^{\ast}(6)\curvearrowright$} & 
{\scriptsize$SU(2,3)/S(U(2)\times U(3))$}\\
{\scriptsize$SO(10)/U(5)$} 
& {\scriptsize$SO^{\ast}(10)/U(5)$} & 
{\scriptsize}\\
\hline
{\scriptsize$SO(5)\times SO(5)\curvearrowright$} 
& {\scriptsize$SO^(5,{\Bbb C})\curvearrowright 
SO^{\ast}(10)/U(5)$} & 
{\scriptsize$SO(5,{\Bbb C})/SO(5)$}\\
{\scriptsize$SO(10)/U(5)$} 
& {\scriptsize} & {\scriptsize}\\
\hline
{\scriptsize$\rho_5(U(5))\curvearrowright SO(10)/U(5)$} 
& {\scriptsize$U(2,3)\curvearrowright 
SO^{\ast}(10)/U(5)$} & 
{\scriptsize$(SO^{\ast}(4)/U(2))\times(SO^{\ast}(6)/U(3))$}\\
\hline
{\scriptsize$SO(2)^2\times SO(3)^2\curvearrowright$} 
& {\scriptsize$SO(2,{\Bbb C})\times SO(3,{\Bbb C})\curvearrowright$} & 
{\scriptsize$SO_0(2,3)/SO(2)\times SO(3)$}\\
{\scriptsize$(SO(5)\times SO(5))/SO(5)$} 
& {\scriptsize$SO(5,{\Bbb C})/SO(5)$} & 
{\scriptsize}\\
\hline
{\scriptsize$\rho_6(SO(5))\curvearrowright$} 
& {\scriptsize$SO_0(2,3)\curvearrowright SO(5,{\Bbb C})/SO(5)$} & 
{\scriptsize$(SO(2,{\Bbb C})/SO(2))$}\\
{\scriptsize$(SO(5)\times SO(5))/SO(5)$} 
& {\scriptsize} & 
{\scriptsize$\times(SO(3,{\Bbb C})/SO(3))$}\\
\hline
{\scriptsize$\rho_7(U(2))\curvearrowright Sp(2)/U(2)$} 
& {\scriptsize$U(1,1)\curvearrowright Sp(2,{\Bbb R})/U(2)$} & 
{\scriptsize$(Sp(1,{\Bbb R})/U(1))$}\\
{\scriptsize} & {\scriptsize} & {\scriptsize$\times(Sp(1,{\Bbb R})/U(1))$}\\
\hline
{\scriptsize$SU(q+2)\curvearrowright$} 
& {\scriptsize$SU(2,q)\curvearrowright$} & 
{\scriptsize$SU(2,q)/S(U(2)\times U(q))$}\\
{\scriptsize$Sp(q+2)/Sp(2)\times Sp(q)$} 
& {\scriptsize$Sp(2,q)/Sp(2)\times Sp(q)$} & 
{\scriptsize}\\
\hline
\end{tabular}$$

\vspace{0.3truecm}

\centerline{{\bf Table 1.}}

\newpage

$$\begin{tabular}{|c|c|c|}
\hline
{\scriptsize$H\curvearrowright G/K$} & {\scriptsize$H^{\ast}\curvearrowright 
G^{\ast}/K$} & {\scriptsize$L^{\ast}/H\cap K$}\\
\hline
{\scriptsize$U(4)\curvearrowright$} 
& {\scriptsize$U^{\ast}(4)\curvearrowright$} & 
{\scriptsize$Sp(2,{\Bbb C})/Sp(2)$}\\
{\scriptsize$Sp(4)/Sp(2)\times Sp(2)$} 
& {\scriptsize$Sp(2,2)/Sp(2)\times Sp(2)$} & 
{\scriptsize}\\
\hline
{\scriptsize$Sp(j+1)\times Sp(q-j+1)\curvearrowright$} 
& {\scriptsize$Sp(1,j)\times Sp(1,q-j)\curvearrowright$} & 
{\scriptsize$(Sp(1,j)/Sp(1)\times Sp(j))\times$}\\
{\scriptsize$Sp(q+2)/Sp(2)\times Sp(q)$} 
& {\scriptsize$Sp(2,q)/Sp(2)\times Sp(q)$} & 
{\scriptsize$(Sp(1,q-j)/Sp(1)\times Sp(q-j))$}\\
\hline
{\scriptsize$SU(2)^2\cdot SO(2)^2\curvearrowright$} 
& {\scriptsize$SL(2,{\Bbb C})\cdot SO(2,{\Bbb C})\curvearrowright$} & 
{\scriptsize$Sp(2,{\Bbb R})/U(2)$}\\
{\scriptsize$(Sp(2)\times Sp(2))/Sp(2)$} 
& {\scriptsize$Sp(2,{\Bbb C})/Sp(2)$} & {\scriptsize}\\
\hline
{\scriptsize$\rho_8(Sp(2))\curvearrowright$} 
& {\scriptsize$Sp(2,{\Bbb R})\curvearrowright Sp(2,{\Bbb C})/Sp(2)$} & 
{\scriptsize$(SL(2,{\Bbb C})/SU(2))$}\\
{\scriptsize$(Sp(2)\times Sp(2))/Sp(2)$} 
& {\scriptsize} & {\scriptsize$\times(SO(2,{\Bbb C})/SO(2))$}\\
\hline
{\scriptsize$\rho_9(Sp(2))\curvearrowright$} 
& {\scriptsize$Sp(1,1)\curvearrowright Sp(2,{\Bbb C})/Sp(2)$} & 
{\scriptsize$(Sp(1,{\Bbb C})/Sp(1))$}\\
{\scriptsize$(Sp(2)\times Sp(2))/Sp(2)$} 
& {\scriptsize} & {\scriptsize$\times(Sp(1,{\Bbb C})/Sp(1))$}\\
\hline
{\scriptsize$Sp(4)\curvearrowright E_6/Spin(10)\cdot U(1)$} 
& {\scriptsize$Sp(2,2)\curvearrowright E_6^{-14}/Spin(10)\cdot U(1)$} & 
{\scriptsize$Sp(2,2)/Sp(2)\times Sp(2)$}\\
\hline
{\scriptsize$SU(6)\cdot SU(2)\curvearrowright$} 
& {\scriptsize$SU(2,4)\cdot SU(2)\curvearrowright$} & 
{\scriptsize$SU(2,4)/S(U(2)\times U(4))$}\\
{\scriptsize$ E_6/Spin(10)\cdot U(1)$} 
& {\scriptsize$E_6^{-14}/Spin(10)\cdot U(1)$} & 
{\scriptsize}\\
\hline
{\scriptsize$\rho_{10}(SU(6)\cdot SU(2))\curvearrowright$} 
& {\scriptsize$SU(1,5)\cdot SL(2,{\Bbb R})\curvearrowright$} & 
{\scriptsize$SO^{\ast}(10)/U(5)$}\\
{\scriptsize$ E_6/Spin(10)\cdot U(1)$} 
& {\scriptsize$E_6^{-14}/Spin(10)\cdot U(1)$} & 
{\scriptsize}\\
\hline
{\scriptsize$\rho_{11}(Spin(10)\cdot U(1))\curvearrowright$} 
& {\scriptsize$SO^{\ast}(10)\cdot U(1)\curvearrowright$} & 
{\scriptsize$(SU(1,5)/S(U(1)\times U(5)))$}\\
{\scriptsize$ E_6/Spin(10)\cdot U(1)$} 
& {\scriptsize$E_6^{-14}/Spin(10)\cdot U(1)$} & 
{\scriptsize$\times(SL(2,{\Bbb R})/SO(2))$}\\
\hline
{\scriptsize$\rho_{12}(Spin(10)\cdot U(1))\curvearrowright$} 
& {\scriptsize$SO_0(2,8)\cdot U(1)\curvearrowright$} & 
{\scriptsize$SO_0(2,8)/SO(2)\times SO(8)$}\\
{\scriptsize$ E_6/Spin(10)\cdot U(1)$} 
& {\scriptsize$E_6^{-14}/Spin(10)\cdot U(1)$} & 
{\scriptsize}\\
\hline
{\scriptsize$Sp(4)\curvearrowright E_6/F_4$} 
& {\scriptsize$Sp(1,3)\curvearrowright E_6^{-26}/F_4$} & 
{\scriptsize$SU^{\ast}(6)/Sp(3)$}\\
\hline
{\scriptsize$\rho_{13}(F_4)\curvearrowright E_6/F_4$} 
& {\scriptsize$F_4^{-20}\curvearrowright E_6^{-26}/F_4$} & 
{\scriptsize$(SO_0(1,9)/SO(9))\times U(1)$}\\
\hline
{\scriptsize$\rho_{14}(SO(4))\curvearrowright G_2/SO(4)$} 
& {\scriptsize$SL(2,{\Bbb R})\times SL(2,{\Bbb R})
\curvearrowright G_2^2/SO(4)$} & 
{\scriptsize$SO(4)/SO(2)\times SO(2)$}\\
\hline
{\scriptsize$\rho_{15}(SO(4))\curvearrowright G_2/SO(4)$} 
& {\scriptsize$\rho^{\ast}_{15}(SO(4))\curvearrowright G_2^2/SO(4)$} & 
{\scriptsize$(SL(2,{\Bbb R})/SO(2))$}\\
{\scriptsize} & {\scriptsize} & 
{\scriptsize$\times(SL(2,{\Bbb R})/SO(2))$}\\
\hline
{\scriptsize$\rho_{16}(G_2)\curvearrowright(G_2\times G_2)/G_2$} 
& {\scriptsize$G_2^2\curvearrowright G_2^{\bf C}/G_2$} & 
{\scriptsize$(SL(2,{\Bbb C})/SU(2))$}\\
{\scriptsize} & {\scriptsize} & 
{\scriptsize$\times(SL(2,{\Bbb C})/SU(2))$}\\
\hline
{\scriptsize$SU(2)^4\curvearrowright(G_2\times G_2)/G_2$} 
& {\scriptsize$SL(2,{\Bbb C})\times SL(2,{\Bbb C})\curvearrowright 
G_2^{\bf C}/G_2$} & 
{\scriptsize$G_2^2/SO(4)$}\\
\hline
\end{tabular}$$

\vspace{0.3truecm}

\centerline{{\bf Table 1(continued).}}

\vspace{0.5truecm}

$$\begin{tabular}{|c|c|c|c|}
\hline
{\scriptsize$H\curvearrowright G/K$} & {\scriptsize$\triangle_+=\triangle'_+$} 
& {\scriptsize${\triangle'}^V_+$} & {\scriptsize${\triangle'}^H_+$}\\
\hline
{\scriptsize$\rho_1(SO(3))\curvearrowright SU(3)/SO(3)$} & 
{\scriptsize$\displaystyle{\{\mathop{\alpha}_{(1)},
\mathop{\beta}_{(1)},\mathop{\alpha+\beta}_{(1)}\}}$} & 
{\scriptsize$\displaystyle{\{\mathop{\alpha}_{(1)}\}}$} & 
{\scriptsize$\displaystyle{\{\mathop{\beta}_{(1)},
\mathop{\alpha+\beta}_{(1)}\}}$}\\
\hline
{\scriptsize$SO(6)\curvearrowright SU(6)/Sp(3)$} & 
{\scriptsize$\displaystyle{\{\mathop{\alpha}_{(4)},
\mathop{\beta}_{(4)},\mathop{\alpha+\beta}_{(4)}\}}$} & 
{\scriptsize$\displaystyle{\{\mathop{\alpha}_{(2)},\mathop{\beta}_{(2)},
\mathop{\alpha+\beta}_{(2)}\}}$} & 
{\scriptsize$\displaystyle{\{\mathop{\alpha}_{(2)},\mathop{\beta}_{(2)},
\mathop{\alpha+\beta}_{(2)}\}}$}\\
\hline
{\scriptsize$\rho_2(Sp(3))\curvearrowright SU(6)/Sp(3)$} & 
{\scriptsize$\displaystyle{\{\mathop{\alpha}_{(4)},
\mathop{\beta}_{(4)},\mathop{\alpha+\beta}_{(4)}\}}$} & 
{\scriptsize$\displaystyle{\{\mathop{\alpha}_{(4)}\}}$} & 
{\scriptsize$\displaystyle{\{\mathop{\beta}_{(4)},
\mathop{\alpha+\beta}_{(4)}\}}$}\\
\hline
{\scriptsize$SO(q+2)\curvearrowright$} & 
{\scriptsize$\displaystyle{\{\mathop{\alpha}_{(2q-4)},
\mathop{\beta}_{(2)},\mathop{\alpha+\beta}_{(2q-4)},}$} & 
{\scriptsize$\displaystyle{\{\mathop{\alpha}_{(q-2)},\mathop{\beta}_{(1)},
\mathop{\alpha+\beta}_{(q-2)},}$} & 
{\scriptsize$\displaystyle{\{\mathop{\alpha}_{(q-2)},\mathop{\beta}_{(1)},
\mathop{\alpha+\beta}_{(q-2)},}$}\\
{\scriptsize$SU(q+2)/S(U(2)\times U(q))$} & 
{\scriptsize$\displaystyle{\mathop{2\alpha+\beta}_{(2)},\mathop{2\alpha}_{(1)},
\mathop{2\alpha+2\beta}_{(1)}\}}$} & 
{\scriptsize$\displaystyle{\mathop{2\alpha+\beta}_{(1)}\}}$} & 
{\scriptsize$\displaystyle{\mathop{2\alpha+\beta}_{(1)},\mathop{2\alpha}_{(1)},
\mathop{2\alpha+2\beta}_{(1)}\}}$}\\
{\scriptsize$(q\,>\,2)$}&&&\\
\hline
{\scriptsize$S(U(j+1)\times U(q-j+1))\curvearrowright$} & 
{\scriptsize$\displaystyle{\{\mathop{\alpha}_{(2q-4)},
\mathop{\beta}_{(2)},\mathop{\alpha+\beta}_{(2q-4)},}$} & 
{\scriptsize$\displaystyle{\{\mathop{\alpha}_{(2j-2)},
\mathop{\alpha+\beta}_{(2q-2j-2)},}$} & 
{\scriptsize$\displaystyle{\{\mathop{\alpha}_{(2q-2j-6)},
\mathop{\beta}_{(2)},}$}\\
{\scriptsize$SU(q+2)/S(U(2)\times U(q))$} & 
{\scriptsize$\displaystyle{\mathop{2\alpha+\beta}_{(2)},\mathop{2\alpha}_{(1)},
\mathop{2\alpha+2\beta}_{(1)}\}}$} & 
{\scriptsize$\displaystyle{\mathop{2\alpha}_{(1)},
\mathop{2\alpha+2\beta}_{(1)}\}}$} & 
{\scriptsize$\displaystyle{\mathop{\alpha+\beta}_{(2j-2)},
\mathop{2\alpha+\beta}_{(2)}\}}$}\\
{\scriptsize$(q>2)$} & & & \\
\hline
{\scriptsize$S(U(2)\times U(2))\curvearrowright$} & 
{\scriptsize$\displaystyle{\{\mathop{\alpha}_{(2)},
\mathop{\beta}_{(1)},\mathop{\alpha+\beta}_{(2)},}$} & 
{\scriptsize$\displaystyle{\{\mathop{\alpha}_{(1)},
\mathop{\alpha+\beta}_{(1)}\}}$} & 
{\scriptsize$\displaystyle{\{\mathop{\alpha}_{(1)},
\mathop{\beta}_{(1)},}$}\\
{\scriptsize$SU(4)/S(U(2)\times U(2))$} & 
{\scriptsize$\displaystyle{\mathop{2\alpha+\beta}_{(1)}\}}$} & 
{\scriptsize} & 
{\scriptsize$\displaystyle{\mathop{\alpha+\beta}_{(1)},
\mathop{2\alpha+\beta}_{(1)}\}}$}\\
{\scriptsize (non-isotropy gr. act.)} & & &\\
\hline
{\scriptsize$SO(j+1)\times SO(q-j+1)\curvearrowright$} & 
{\scriptsize$\displaystyle{\{\mathop{\alpha}_{(q-2)},
\mathop{\beta}_{(1)},\mathop{\alpha+\beta}_{(q-2)},}$} & 
{\scriptsize$\displaystyle{\{\mathop{\alpha}_{(j-1)},
\mathop{\alpha+\beta}_{(q-j-1)}\}}$} & 
{\scriptsize$\displaystyle{\{\mathop{\alpha}_{(q-j-1)},
\mathop{\beta}_{(1)},}$}\\
{\scriptsize$SO(q+2)/SO(2)\times SO(q)$} & 
{\scriptsize$\displaystyle{\mathop{2\alpha+\beta}_{(1)}\}}$} & 
{\scriptsize} & 
{\scriptsize$\displaystyle{\mathop{\alpha+\beta}_{(j-1)},
\mathop{2\alpha+\beta}_{(1)}\}}$}\\
\hline
{\scriptsize$SO(4)\times SO(4)\curvearrowright$} & 
{\scriptsize$\displaystyle{\{\mathop{\alpha}_{(4)},
\mathop{\beta}_{(1)},\mathop{\alpha+\beta}_{(4)},}$} & 
{\scriptsize$\displaystyle{\{\mathop{\alpha}_{(2)},\mathop{\beta}_{(1)},
\mathop{\alpha+\beta}_{(2)},}$} & 
{\scriptsize$\displaystyle{\{\mathop{\alpha}_{(2)},
\mathop{\alpha+\beta}_{(2)}\}}$}\\
{\scriptsize$SO(8)/U(4)$} & 
{\scriptsize$\displaystyle{\mathop{2\alpha+\beta}_{(1)}\}}$} & 
{\scriptsize$\displaystyle{\mathop{2\alpha+\beta}_{(1)}\}}$} & 
{\scriptsize}\\
\hline
{\scriptsize$\rho_3(SO(4)\times SO(4))\curvearrowright$} 
& {\scriptsize$\displaystyle{\{\mathop{\alpha}_{(4)},
\mathop{\beta}_{(1)},\mathop{\alpha+\beta}_{(4)},}$} & 
{\scriptsize$\displaystyle{\{\mathop{\alpha}_{(2)},\mathop{\alpha+\beta}_{(2)}
\}}$} & 
{\scriptsize$\displaystyle{\{\mathop{\alpha}_{(2)},\mathop{\beta}_{(1)},
\mathop{\alpha+\beta}_{(2)},}$}\\
{\scriptsize$SO(8)/U(4)$} & 
{\scriptsize$\displaystyle{\mathop{2\alpha+\beta}_{(1)}\}}$} & 
{\scriptsize} & 
{\scriptsize$\displaystyle{\mathop{2\alpha+\beta}_{(1)}}$}\\
\hline
{\scriptsize$\rho_4(U(4))\curvearrowright SO(8)/U(4)$} 
& {\scriptsize$\displaystyle{\{\mathop{\alpha}_{(4)},
\mathop{\beta}_{(1)},\mathop{\alpha+\beta}_{(4)},}$} & 
{\scriptsize$\displaystyle{\{\mathop{\alpha}_{(1)},\mathop{\alpha+\beta}_{(1)}
\}}$} & 
{\scriptsize$\displaystyle{\{\mathop{\alpha}_{(3)},\mathop{\beta}_{(1)},
\mathop{\alpha+\beta}_{(3)},}$}\\
{\scriptsize} & 
{\scriptsize$\displaystyle{\mathop{2\alpha+\beta}_{(1)}\}}$} & 
{\scriptsize} & 
{\scriptsize$\displaystyle{\mathop{2\alpha+\beta}_{(1)}}$}\\
\hline
{\scriptsize$SO(4)\times SO(6)\curvearrowright$} 
& {\scriptsize$\displaystyle{\{\mathop{\alpha}_{(4)},
\mathop{\beta}_{(4)},\mathop{\alpha+\beta}_{(4)},}$} & 
{\scriptsize$\displaystyle{\{\mathop{\alpha}_{(2)},\mathop{\beta}_{(2)},
\mathop{\alpha+\beta}_{(2)},}$} & 
{\scriptsize$\displaystyle{\{\mathop{\alpha}_{(2)},\mathop{\beta}_{(2)},
\mathop{\alpha+\beta}_{(2)},}$}\\
{\scriptsize$SO(10)/U(5)$} & 
{\scriptsize$\displaystyle{\mathop{2\alpha+\beta}_{(4)},
\mathop{2\alpha}_{(1)},\mathop{2\alpha+2\beta}_{(1)}\}}$} & 
{\scriptsize$\displaystyle{\mathop{2\alpha+\beta}_{(2)},
\mathop{2\alpha}_{(1)},\mathop{2\alpha+2\beta}_{(1)}\}}$} & 
{\scriptsize$\displaystyle{\mathop{2\alpha+\beta}_{(2)}}$}\\
\hline
{\scriptsize$SO(5)\times SO(5)\curvearrowright$} 
& {\scriptsize$\displaystyle{\{\mathop{\alpha}_{(4)},
\mathop{\beta}_{(4)},\mathop{\alpha+\beta}_{(4)},}$} & 
{\scriptsize$\displaystyle{\{\mathop{\alpha}_{(2)},\mathop{\beta}_{(2)},
\mathop{\alpha+\beta}_{(2)},}$} & 
{\scriptsize$\displaystyle{\{\mathop{\alpha}_{(2)},\mathop{\beta}_{(2)},
\mathop{\alpha+\beta}_{(2)},}$}\\
{\scriptsize$SO(10)/U(5)$} & 
{\scriptsize$\displaystyle{\mathop{2\alpha+\beta}_{(4)},
\mathop{2\alpha}_{(1)},\mathop{2\alpha+2\beta}_{(1)}\}}$} & 
{\scriptsize$\displaystyle{\mathop{2\alpha+\beta}_{(2)}\}}$} & 
{\scriptsize$\displaystyle{\mathop{2\alpha+\beta}_{(2)},
\mathop{2\alpha}_{(1)},\mathop{2\alpha+2\beta}_{(1)}}$}\\
\hline
\end{tabular}$$

\vspace{0.3truecm}

\centerline{{\bf Table 2.}}

\newpage

%\vspace{0.3truecm}

$$\begin{tabular}{|c|c|c|c|}
\hline
{\scriptsize$H\curvearrowright G/K$} & {\scriptsize$\triangle_+=\triangle'_+$} 
& {\scriptsize${\triangle'}^V_+$} & {\scriptsize${\triangle'}^H_+$}\\
\hline
{\scriptsize$\rho_5(U(5))\curvearrowright SO(10)/U(5)$} 
& {\scriptsize$\displaystyle{\{\mathop{\alpha}_{(4)},
\mathop{\beta}_{(4)},\mathop{\alpha+\beta}_{(4)},}$} & 
{\scriptsize$\displaystyle{\{\mathop{\alpha}_{(4)},\mathop{2\alpha}_{(1)},
\mathop{2\alpha+2\beta}_{(1)}\}}$} & 
{\scriptsize$\displaystyle{\{\mathop{\beta}_{(4)},\mathop{\alpha+\beta}_{(4)},
\mathop{2\alpha+\beta}_{(4)}\}}$}\\
{\scriptsize} & 
{\scriptsize$\displaystyle{\mathop{2\alpha+\beta}_{(4)},
\mathop{2\alpha}_{(1)},\mathop{2\alpha+2\beta}_{(1)}\}}$} & 
{\scriptsize} & {\scriptsize}\\
\hline
{\scriptsize$SO(2)^2\times SO(3)^2\curvearrowright$} 
& {\scriptsize$\displaystyle{\{\mathop{\alpha}_{(2)},
\mathop{\beta}_{(2)},\mathop{\alpha+\beta}_{(2)},}$} & 
{\scriptsize$\displaystyle{\{\mathop{\alpha}_{(1)},
\mathop{\beta}_{(1)},\mathop{\alpha+\beta}_{(1)},}$} & 
{\scriptsize$\displaystyle{\{\mathop{\alpha}_{(1)},\mathop{\beta}_{(1)},
\mathop{\alpha+\beta}_{(1)},}$}\\
{\scriptsize$(SO(5)\times SO(5))/SO(5)$} & 
{\scriptsize$\displaystyle{\{\mathop{2\alpha+\beta}_{(2)}\}}$} & 
{\scriptsize$\displaystyle{\mathop{2\alpha+\beta}_{(1)}\}}$} & 
{\scriptsize$\displaystyle{\mathop{2\alpha+\beta}_{(1)}\}}$}\\
\hline
{\scriptsize$\rho_6(SO(5))\curvearrowright$} 
& {\scriptsize$\displaystyle{\{\mathop{\alpha}_{(2)},
\mathop{\beta}_{(2)},\mathop{\alpha+\beta}_{(2)},}$} & 
{\scriptsize$\displaystyle{\{\mathop{\alpha}_{(2)}\}}$} & 
{\scriptsize$\displaystyle{\{\mathop{\beta}_{(2)},
\mathop{\alpha+\beta}_{(2)},\mathop{2\alpha+\beta}_{(2)}\}}$}\\
{\scriptsize$(SO(5)\times SO(5))/SO(5)$} & 
{\scriptsize$\displaystyle{\{\mathop{2\alpha+\beta}_{(2)}\}}$} & 
{\scriptsize} & {\scriptsize}\\
\hline
{\scriptsize$\rho_7(U(2))\curvearrowright Sp(2)/U(2)$} 
& {\scriptsize$\displaystyle{\{\mathop{\alpha}_{(1)},
\mathop{\beta}_{(1)},\mathop{\alpha+\beta}_{(1)},}$} & 
{\scriptsize$\displaystyle{\{\mathop{\alpha}_{(1)},
\mathop{\alpha+\beta}_{(1)}\}}$} & 
{\scriptsize$\displaystyle{\{\mathop{\beta}_{(1)},
\mathop{2\alpha+\beta}_{(1)}\}}$}\\
{\scriptsize} & 
{\scriptsize$\displaystyle{\{\mathop{2\alpha+\beta}_{(1)}\}}$} & 
{\scriptsize} & {\scriptsize}\\
\hline
{\scriptsize$SU(q+2)\curvearrowright$} 
& {\scriptsize$\displaystyle{\{\mathop{\alpha}_{(4q-8)},
\mathop{\beta}_{(4)},\mathop{\alpha+\beta}_{(4q-8)},}$} & 
{\scriptsize$\displaystyle{\{\mathop{\alpha}_{(2q-4)},\mathop{\beta}_{(2)},
\mathop{\alpha+\beta}_{(2q-4)}\}}$} & 
{\scriptsize$\displaystyle{\{\mathop{\alpha}_{(2q-4)},\mathop{\beta}_{(2)},
\mathop{\alpha+\beta}_{(2q-4)}\}}$}\\
{\scriptsize$Sp(q+2)/Sp(2)\times Sp(q)$} & 
{\scriptsize$\displaystyle{\mathop{2\alpha+\beta}_{(4)},
\mathop{2\alpha}_{(3)},\mathop{2\alpha+2\beta}_{(3)}\}}$} & 
{\scriptsize$\displaystyle{\mathop{2\alpha+\beta}_{(2)},
\mathop{2\alpha}_{(1)},\mathop{2\alpha+2\beta}_{(1)}\}}$} & 
{\scriptsize$\displaystyle{\mathop{2\alpha+\beta}_{(2)},
\mathop{2\alpha}_{(2)},\mathop{2\alpha+2\beta}_{(2)}\}}$}\\
{\scriptsize$(q>2)$}&&&\\
\hline
{\scriptsize$SU(4)\curvearrowright$} 
& {\scriptsize$\displaystyle{\{\mathop{\alpha}_{(4)},
\mathop{\beta}_{(3)},\mathop{\alpha+\beta}_{(3)},}$} & 
{\scriptsize$\displaystyle{\{\mathop{\alpha}_{(2)},\mathop{\beta}_{(1)},
\mathop{\alpha+\beta}_{(2)}\}}$} & 
{\scriptsize$\displaystyle{\{\mathop{\alpha}_{(2)},\mathop{\beta}_{(2)},
\mathop{\alpha+\beta}_{(1)}\}}$}\\
{\scriptsize$Sp(4)/Sp(2)\times Sp(2)$} & 
{\scriptsize$\displaystyle{\mathop{2\alpha+\beta}_{(4)}\}}$} & 
{\scriptsize$\displaystyle{\mathop{2\alpha+\beta}_{(1)}\}}$} & 
{\scriptsize$\displaystyle{\mathop{2\alpha+\beta}_{(3)}\}}$}\\
\hline
{\scriptsize$U(4)\curvearrowright$} 
& {\scriptsize$\displaystyle{\{\mathop{\alpha}_{(4)},
\mathop{\beta}_{(3)},\mathop{\alpha+\beta}_{(3)},}$} & 
{\scriptsize$\displaystyle{\{\mathop{\alpha}_{(2)},\mathop{\beta}_{(2)},
\mathop{\alpha+\beta}_{(2)}\}}$} & 
{\scriptsize$\displaystyle{\{\mathop{\alpha}_{(2)},\mathop{\beta}_{(1)},
\mathop{\alpha+\beta}_{(1)}\}}$}\\
{\scriptsize$Sp(4)/Sp(2)\times Sp(2)$} & 
{\scriptsize$\displaystyle{\mathop{2\alpha+\beta}_{(4)}\}}$} & 
{\scriptsize$\displaystyle{\mathop{2\alpha+\beta}_{(2)}\}}$} & 
{\scriptsize$\displaystyle{\mathop{2\alpha+\beta}_{(2)}\}}$}\\
\hline
{\scriptsize$Sp(j+1)\times Sp(q-j+1)\curvearrowright$} 
& {\scriptsize$\displaystyle{\{\mathop{\alpha}_{(4q-8)},
\mathop{\beta}_{(4)},\mathop{\alpha+\beta}_{(4q-8)},}$} & 
{\scriptsize$\displaystyle{\{\mathop{\alpha}_{(2j-4)},\mathop{2\alpha}_{(3)},
\mathop{\alpha+\beta}_{(4q-4j-4)}\}}$} & 
{\scriptsize$\displaystyle{\{\mathop{\alpha}_{(4q-4j-4)},\mathop{\beta}_{(4)},
\mathop{\alpha+\beta}_{(4j-4)}\}}$}\\
{\scriptsize$Sp(q+2)/Sp(2)\times Sp(q)$} & 
{\scriptsize$\displaystyle{\mathop{2\alpha+\beta}_{(4)},
\mathop{2\alpha}_{(3)},\mathop{2\alpha+2\beta}_{(3)}\}}$} & 
{\scriptsize$\displaystyle{\mathop{2\alpha+2\beta}_{(3)}\}}$} & 
{\scriptsize$\displaystyle{\mathop{2\alpha+\beta}_{(4)}\}}$}\\
{\scriptsize$(q>2)$}&&&\\
\hline
{\scriptsize$Sp(2)\times Sp(2)\curvearrowright$} 
& {\scriptsize$\displaystyle{\{\mathop{\alpha}_{(4)},
\mathop{\beta}_{(3)},\mathop{\alpha+\beta}_{(3)},}$} & 
{\scriptsize$\displaystyle{\{\mathop{\alpha}_{(3)},
\mathop{\alpha+\beta}_{(3)}\}}$} & 
{\scriptsize$\displaystyle{\{\mathop{\alpha}_{(1)},\mathop{\beta}_{(3)},
\mathop{2\alpha+\beta}_{(4)}\}}$}\\
{\scriptsize$Sp(4)/Sp(2)\times Sp(2)$} & 
{\scriptsize$\displaystyle{\mathop{2\alpha+\beta}_{(4)}\}}$} & 
{\scriptsize} & {\scriptsize}\\
\hline
{\scriptsize$SU(2)^2\cdot SO(2)^2\curvearrowright$} 
& {\scriptsize$\displaystyle{\{\mathop{\alpha}_{(2)},
\mathop{\beta}_{(2)},\mathop{\alpha+\beta}_{(2)},}$} & 
{\scriptsize$\displaystyle{\{\mathop{\alpha}_{(1)},\mathop{\beta}_{(1)},
\mathop{\alpha+\beta}_{(1)},}$} & 
{\scriptsize$\displaystyle{\{\mathop{\alpha}_{(1)},\mathop{\beta}_{(1)},
\mathop{\alpha+\beta}_{(1)},}$}\\
{\scriptsize$(Sp(2)\times Sp(2))/Sp(2)$} & 
{\scriptsize$\displaystyle{\mathop{2\alpha+\beta}_{(2)}\}}$} & 
{\scriptsize$\displaystyle{\mathop{2\alpha+\beta}_{(1)}\}}$} & 
{\scriptsize$\displaystyle{\mathop{2\alpha+\beta}_{(1)}\}}$}\\
\hline
{\scriptsize$\rho_8(Sp(2))\curvearrowright$} 
& {\scriptsize$\displaystyle{\{\mathop{\alpha}_{(2)},
\mathop{\beta}_{(2)},\mathop{\alpha+\beta}_{(2)},}$} & 
{\scriptsize$\displaystyle{\{\mathop{\alpha}_{(2)},
\mathop{\alpha+\beta}_{(2)}\}}$} & 
{\scriptsize$\displaystyle{\{\mathop{\beta}_{(2)},
\mathop{2\alpha+\beta}_{(2)}\}}$}\\
{\scriptsize$(Sp(2)\times Sp(2))/Sp(2)$} & 
{\scriptsize$\displaystyle{\mathop{2\alpha+\beta}_{(2)}\}}$} & 
{\scriptsize} & {\scriptsize}\\
\hline
\end{tabular}$$

\vspace{0.3truecm}

\centerline{{\bf Table 2(continued).}}

\newpage

$$\begin{tabular}{|c|c|c|c|}
\hline
{\scriptsize$H\curvearrowright G/K$} & {\scriptsize$\triangle_+=\triangle'_+$} 
& {\scriptsize${\triangle'}^V_+$} & {\scriptsize${\triangle'}^H_+$}\\
\hline
{\scriptsize$\rho_9(Sp(2))\curvearrowright$} 
& {\scriptsize$\displaystyle{\{\mathop{\alpha}_{(2)},
\mathop{\beta}_{(2)},\mathop{\alpha+\beta}_{(2)},}$} & 
{\scriptsize$\displaystyle{\{\mathop{\alpha}_{(2)},
\mathop{\alpha+\beta}_{(2)}\}}$} & 
{\scriptsize$\displaystyle{\{\mathop{\beta}_{(2)},
\mathop{2\alpha+\beta}_{(2)}\}}$}\\
{\scriptsize$(Sp(2)\times Sp(2))/Sp(2)$} & 
{\scriptsize$\displaystyle{\mathop{2\alpha+\beta}_{(2)}\}}$} & 
{\scriptsize} & {\scriptsize}\\
\hline
{\scriptsize$Sp(4)\curvearrowright$} 
& {\scriptsize$\displaystyle{\{\mathop{\alpha}_{(8)},
\mathop{\beta}_{(6)},\mathop{\alpha+\beta}_{(9)},}$} & 
{\scriptsize$\displaystyle{\{\mathop{\alpha}_{(4)},\mathop{\beta}_{(3)},
\mathop{\alpha+\beta}_{(3)},}$} & 
{\scriptsize$\displaystyle{\{\mathop{\alpha}_{(4)},\mathop{\beta}_{(3)},
\mathop{\alpha+\beta}_{(6)},}$}\\
{\scriptsize$\displaystyle{E_6/Spin(10)\cdot U(1)}$} & 
{\scriptsize$\displaystyle{\mathop{2\alpha+\beta}_{(5)},
\mathop{2\alpha}_{(1)},\mathop{2\alpha+2\beta}_{(1)}\}}$} & 
{\scriptsize$\displaystyle{\mathop{2\alpha+\beta}_{(4)}\}}$} & 
{\scriptsize$\displaystyle{\mathop{2\alpha+\beta}_{(1)}\}}$}\\
\hline
{\scriptsize$SU(6)\cdot SU(2)\curvearrowright$} 
& {\scriptsize$\displaystyle{\{\mathop{\alpha}_{(8)},
\mathop{\beta}_{(6)},\mathop{\alpha+\beta}_{(9)},}$} & 
{\scriptsize$\displaystyle{\{\mathop{\alpha}_{(4)},\mathop{\beta}_{(2)},
\mathop{\alpha+\beta}_{(4)},}$} & 
{\scriptsize$\displaystyle{\{\mathop{\alpha}_{(4)},\mathop{\beta}_{(4)},
\mathop{\alpha+\beta}_{(5)},}$}\\
{\scriptsize$E_6/Spin(10)\cdot U(1)$} & 
{\scriptsize$\displaystyle{\mathop{2\alpha+\beta}_{(5)},
\mathop{2\alpha}_{(1)},\mathop{2\alpha+2\beta}_{(1)}\}}$} & 
{\scriptsize$\displaystyle{\mathop{2\alpha+\beta}_{(2)},
\mathop{2\alpha}_{(1)},\mathop{2\alpha+2\beta}_{(1)}\}}$} & 
{\scriptsize$\displaystyle{\mathop{2\alpha+\beta}_{(3)}\}}$}\\
\hline
{\scriptsize$\rho_{10}(SU(6)\cdot SU(2))\curvearrowright$} 
& {\scriptsize$\displaystyle{\{\mathop{\alpha}_{(8)},
\mathop{\beta}_{(6)},\mathop{\alpha+\beta}_{(9)},}$} & 
{\scriptsize$\displaystyle{\{\mathop{\alpha}_{(4)},\mathop{\beta}_{(4)},
\mathop{\alpha+\beta}_{(4)},}$} & 
{\scriptsize$\displaystyle{\{\mathop{\alpha}_{(4)},\mathop{\beta}_{(2)},
\mathop{\alpha+\beta}_{(5)},}$}\\
{\scriptsize$E_6/Spin(10)\cdot U(1)$} & 
{\scriptsize$\displaystyle{\mathop{2\alpha+\beta}_{(5)},
\mathop{2\alpha}_{(1)},\mathop{2\alpha+2\beta}_{(1)}\}}$} & 
{\scriptsize$\displaystyle{\mathop{2\alpha+\beta}_{(4)},
\mathop{2\alpha}_{(1)},\mathop{2\alpha+2\beta}_{(1)}\}}$} & 
{\scriptsize$\displaystyle{\mathop{2\alpha+\beta}_{(1)}\}}$}\\
\hline
{\scriptsize$\rho_{11}(Spin(10)\cdot U(1))\curvearrowright$} 
& {\scriptsize$\displaystyle{\{\mathop{\alpha}_{(8)},
\mathop{\beta}_{(6)},\mathop{\alpha+\beta}_{(9)},}$} & 
{\scriptsize$\displaystyle{\{\mathop{\alpha}_{(8)},\mathop{2\alpha}_{(1)},
\mathop{2\alpha+2\beta}_{(1)}\}}$} & 
{\scriptsize$\displaystyle{\{\mathop{\beta}_{(6)},\mathop{\alpha+\beta}_{(9)},
\mathop{2\alpha+\beta}_{(5)}\}}$}\\
{\scriptsize$E_6/Spin(10)\cdot U(1)$} & 
{\scriptsize$\displaystyle{\mathop{2\alpha+\beta}_{(5)},
\mathop{2\alpha}_{(1)},\mathop{2\alpha+2\beta}_{(1)}\}}$} & 
{\scriptsize} & {\scriptsize}\\
\hline
{\scriptsize$\rho_{12}(Spin(10)\cdot U(1))\curvearrowright$} 
& {\scriptsize$\displaystyle{\{\mathop{\alpha}_{(8)},
\mathop{\beta}_{(6)},\mathop{\alpha+\beta}_{(9)},}$} & 
{\scriptsize$\displaystyle{\{\mathop{\alpha}_{(6)},\mathop{\beta}_{(1)},
\mathop{\alpha+\beta}_{(6)},}$} & 
{\scriptsize$\displaystyle{\{\mathop{\alpha}_{(2)},\mathop{\beta}_{(5)},
\mathop{\alpha+\beta}_{(3)},}$}\\
{\scriptsize$E_6/Spin(10)\cdot U(1)$} & 
{\scriptsize$\displaystyle{\mathop{2\alpha+\beta}_{(5)},
\mathop{2\alpha}_{(1)},\mathop{2\alpha+2\beta}_{(1)}\}}$} & 
{\scriptsize$\displaystyle{\mathop{2\alpha+\beta}_{(1)}\}}$} & 
{\scriptsize$\displaystyle{\mathop{2\alpha+\beta}_{(4)},
\mathop{2\alpha}_{(1)},\mathop{2\alpha+2\beta}_{(1)}\}}$}\\
\hline
{\scriptsize$Sp(4)\curvearrowright E_6/F_4$} 
& {\scriptsize$\displaystyle{\{\mathop{\alpha}_{(8)},
\mathop{\beta}_{(8)},\mathop{\alpha+\beta}_{(8)}\}}$} 
& {\scriptsize$\displaystyle{\{\mathop{\alpha}_{(4)},\mathop{\beta}_{(4)},
\mathop{\alpha+\beta}_{(4)}\}}$} 
& {\scriptsize$\displaystyle{\{\mathop{\alpha}_{(4)},\mathop{\beta}_{(4)},
\mathop{\alpha+\beta}_{(4)}\}}$}\\
\hline
{\scriptsize$\rho_{13}(F_4)\curvearrowright E_6/F_4$} 
& {\scriptsize$\displaystyle{\{\mathop{\alpha}_{(8)},
\mathop{\beta}_{(8)},\mathop{\alpha+\beta}_{(8)}\}}$} 
& {\scriptsize$\displaystyle{\{\mathop{\alpha}_{(8)}\}}$} 
& {\scriptsize$\displaystyle{\{\mathop{\beta}_{(8)},
\mathop{\alpha+\beta}_{(8)}\}}$}\\
\hline
{\scriptsize$\rho_{14}(SO(4))\curvearrowright$} 
& {\scriptsize$\displaystyle{\{\mathop{\alpha}_{(1)},
\mathop{\beta}_{(1)},\mathop{\alpha+\beta}_{(1)},}$} & 
{\scriptsize$\displaystyle{\{\mathop{\alpha}_{(1)},
\mathop{3\alpha+2\beta}_{(1)}\}}$} & 
{\scriptsize$\displaystyle{\{\mathop{\beta}_{(1)},\mathop{\alpha+\beta}_{(1)},
\mathop{2\alpha+\beta}_{(1)},}$}\\
{\scriptsize$\displaystyle{G_2/SO(4)}$} & 
{\scriptsize$\displaystyle{\mathop{2\alpha+\beta}_{(1)},
\mathop{3\alpha+\beta}_{(1)},\mathop{3\alpha+2\beta}_{(1)}\}}$} & 
{\scriptsize} & {\scriptsize$\displaystyle{\mathop{3\alpha+\beta}_{(1)}\}}$}\\
\hline
{\scriptsize$\rho_{15}(SO(4))\curvearrowright$} 
& {\scriptsize$\displaystyle{\{\mathop{\alpha}_{(1)},
\mathop{\beta}_{(1)},\mathop{\alpha+\beta}_{(1)},}$} & 
{\scriptsize$\displaystyle{\{\mathop{\alpha}_{(1)},
\mathop{3\alpha+2\beta}_{(1)}\}}$} & 
{\scriptsize$\displaystyle{\{\mathop{\beta}_{(1)},\mathop{\alpha+\beta}_{(1)},
\mathop{2\alpha+\beta}_{(1)},}$}\\
{\scriptsize$\displaystyle{G_2/SO(4)}$} & 
{\scriptsize$\displaystyle{\mathop{2\alpha+\beta}_{(1)},
\mathop{3\alpha+\beta}_{(1)},\mathop{3\alpha+2\beta}_{(1)}\}}$} & 
{\scriptsize} & {\scriptsize$\displaystyle{\mathop{3\alpha+\beta}_{(1)}\}}$}\\
\hline
{\scriptsize$\rho_{16}(G_2)\curvearrowright$} 
& {\scriptsize$\displaystyle{\{\mathop{\alpha}_{(2)},
\mathop{\beta}_{(2)},\mathop{\alpha+\beta}_{(2)},}$} & 
{\scriptsize$\displaystyle{\{\mathop{\alpha}_{(2)},
\mathop{3\alpha+2\beta}_{(2)}\}}$} & 
{\scriptsize$\displaystyle{\{\mathop{\beta}_{(2)},\mathop{\alpha+\beta}_{(2)},
\mathop{2\alpha+\beta}_{(2)},}$}\\
{\scriptsize$\displaystyle{(G_2\times G_2)/G_2}$} & 
{\scriptsize$\displaystyle{\mathop{2\alpha+\beta}_{(2)},
\mathop{3\alpha+\beta}_{(2)},\mathop{3\alpha+2\beta}_{(2)}\}}$} & 
{\scriptsize} & {\scriptsize$\displaystyle{\mathop{3\alpha+\beta}_{(2)}\}}$}\\
\hline
{\scriptsize$SU(2)^4\curvearrowright$} 
& {\scriptsize$\displaystyle{\{\mathop{\alpha}_{(2)},
\mathop{\beta}_{(2)},\mathop{\alpha+\beta}_{(2)},}$} & 
{\scriptsize$\displaystyle{\{\mathop{\alpha}_{(1)},\mathop{\beta}_{(1)},
\mathop{\alpha+\beta}_{(1)},}$} & 
{\scriptsize$\displaystyle{\{\mathop{\alpha}_{(1)},\mathop{\beta}_{(1)},
\mathop{\alpha+\beta}_{(1)},}$}\\
{\scriptsize$\displaystyle{(G_2\times G_2)/G_2}$} & 
{\scriptsize$\displaystyle{\mathop{2\alpha+\beta}_{(2)},
\mathop{3\alpha+\beta}_{(2)},\mathop{3\alpha+2\beta}_{(2)}\}}$} & 
{\scriptsize$\displaystyle{\mathop{2\alpha+\beta}_{(1)},
\mathop{3\alpha+\beta}_{(1)},\mathop{3\alpha+2\beta}_{(1)}\}}$} & 
{\scriptsize$\displaystyle{\mathop{2\alpha+\beta}_{(1)},
\mathop{3\alpha+\beta}_{(1)},\mathop{3\alpha+2\beta}_{(1)}\}}$}\\
\hline
\end{tabular}$$

\vspace{0.3truecm}

\centerline{{\bf Table 2(continued${}^2$).}}

\newpage

\vspace{0.5truecm}

$$\begin{tabular}{|c|c|}
\hline
{\scriptsize$H\curvearrowright G/K$} & {\scriptsize$X\qquad(\widetilde C)$}\\
\hline
{\scriptsize$\rho_1(SO(3))\curvearrowright SU(3)/SO(3)$} & 
{\scriptsize$X_{(x_1,x_2)}=(\tan(x_1+\sqrt3x_2)-2\cot2x_1
+\tan(x_1-\sqrt3x_2),$}\\
{\scriptsize} & 
{\scriptsize$\sqrt3\tan(x_1+\sqrt3x_2)-\sqrt3\tan(x_1-\sqrt3x_2))$}\\
{\scriptsize} & 
{\scriptsize$(\widetilde C\,:\,x_1>0,
x_2>\frac{1}{\sqrt3}x_1-\frac{\pi}{2\sqrt3},
x_2<-\frac{1}{\sqrt3}x_1+\frac{\pi}{2\sqrt3})$}\\
\hline
{\scriptsize$SO(6)\curvearrowright SU(6)/Sp(3)$} & 
{\scriptsize$X_{(x_1,x_2)}=(-4\cot2x_1-2\cot(x_1-\sqrt3x_2)
-2\cot(x_1+\sqrt3x_2)$}\\
{\scriptsize} & 
{\scriptsize$+4\tan2x_1+2\tan(x_1-\sqrt3x_2)
+2\tan(x_1+\sqrt3x_2),$}\\
{\scriptsize} & 
{\scriptsize$2\sqrt3\cot(x_1-\sqrt3x_2)-2\sqrt3\cot(x_1+\sqrt3x_2)$}\\
{\scriptsize} & 
{\scriptsize$-2\sqrt3\tan(x_1-\sqrt3x_2)+2\sqrt3\tan(x_1+\sqrt3x_2))$}\\
{\scriptsize} & 
{\scriptsize$(\widetilde C\,:\,x_1>0,
x_2>\frac{1}{\sqrt3}x_1,
x_2<-\frac{1}{\sqrt3}x_1+\frac{\pi}{2\sqrt3})$}\\
\hline
{\scriptsize$\rho_2(Sp(3))\curvearrowright SU(6)/Sp(3)$} & 
{\scriptsize$X_{(x_1,x_2)}=(-8\cot2x_1+4\tan(x_1-\sqrt3x_2)
+4\cot(x_1+\sqrt3x_2),$}\\
{\scriptsize} & 
{\scriptsize$4\sqrt3\tan(x_1+\sqrt3x_2)
-4\sqrt3\tan(x_1-\sqrt3x_2))$}\\
{\scriptsize} & 
{\scriptsize$(\widetilde C\,:\,x_1>0,
x_2>\frac{1}{\sqrt3}x_1-\frac{\pi}{2\sqrt3},
x_2<-\frac{1}{\sqrt3}x_1+\frac{\pi}{2\sqrt3})$}\\
\hline
{\scriptsize$SO(q+2)\curvearrowright$} & 
{\scriptsize$X_{(x_1,x_2)}=(-(q-2)\cot x_1+\cot(x_1-x_2)
-\cot(x_1+x_2)$}\\
{\scriptsize$SU(q+2)/S(U(2)\times U(q))$} & 
{\scriptsize$+(q-2)\tan x_1+\tan(x_1-x_2)+\tan(x_1+x_2)+2\tan 2x_1,$}\\
{\scriptsize$(q>2)$} & 
{\scriptsize$\cot(x_1-x_2)-(q-2)\cot x_2-\cot(x_1+x_2)$}\\
{\scriptsize} & 
{\scriptsize$-\tan(x_1-x_2)+(q-2)\tan x_2+\tan(x_1+x_2)+2\tan2x_2)$}\\
{\scriptsize} & 
{\scriptsize$(\widetilde C\,:\,x_1>0,
x_2>x_1,x_2<\frac{\pi}{4})$}\\
\hline
{\scriptsize$SO(4)\curvearrowright$} & 
{\scriptsize$X_{(x_1,x_2)}=(-\cot x_1+\tan x_1+\tan(x_1-x_2)+\tan(x_1+x_2),$}\\
{\scriptsize$SU(4)/S(U(2)\times U(2))$} & 
{\scriptsize$-\cot x_2-\tan(x_1-x_2)+\tan x_2+\tan(x_1+x_2))$}\\
{\scriptsize} & 
{\scriptsize$(\widetilde C\,:\,x_1>0,x_2>0,x_1+x_2<\frac{\pi}{2})$}\\
\hline
{\scriptsize$S(U(j+1)\times U(q-j+1))\curvearrowright$} & 
{\scriptsize$X_{(x_1,x_2)}=(-2(j-1)\cot x_1-2\cot2x_1
+2(q-j-1)\tan x_1$}\\
{\scriptsize$SU(q+2)/S(U(2)\times U(q))$} & 
{\scriptsize$+2\tan(x_1-x_2)+2\tan(x_1+x_2),$}\\
{\scriptsize$(q>2)$} & 
{\scriptsize$-2(q-j-1)\cot x_2-2\cot2x_2-2\tan(x_1-x_2)$}\\
{\scriptsize} & 
{\scriptsize$+2(j-1)\tan x_2+2\tan(x_1+x_2))$}\\
{\scriptsize} & 
{\scriptsize$(\widetilde C\,:\,x_1>0,x_2>0,x_1+x_2<\frac{\pi}{2})$}\\
\hline
{\scriptsize$S(U(2)\times U(2))\curvearrowright$} & 
{\scriptsize$X_{(x_1,x_2)}=(-\cot x_1+\tan x_1+\tan(x_1-x_2)+\tan(x_1+x_2),$}\\
{\scriptsize$SU(4)/S(U(2)\times U(2))$} & 
{\scriptsize$-\cot x_2-\tan(x_1-x_2)+\tan x_2+\tan(x_1+x_2))$}\\
{\scriptsize(non-isotropy gr. act.)} & 
{\scriptsize$(\widetilde C\,:\,x_1>0,x_2>0,x_1+x_2<\frac{\pi}{2})$}\\
\hline
{\scriptsize$SO(j+1)\times SO(q-j+1)\curvearrowright$} & 
{\scriptsize$X_{(x_1,x_2)}=(-(j-1)\cot x_1+(q-j-1)\tan x_1$}\\
{\scriptsize$SO(q+2)/SO(2)\times SO(q)$} & 
{\scriptsize$+\tan(x_1-x_2)+\tan(x_1+x_2),$}\\
{\scriptsize$(q>2)$} & 
{\scriptsize$-(q-j-1)\cot x_2-\tan(x_1-x_2)$}\\
{\scriptsize} & 
{\scriptsize$+(j-1)\tan x_2+\tan(x_1+x_2))$}\\
{\scriptsize} & 
{\scriptsize$(\widetilde C\,:\,x_1>0,x_2>0,x_1+x_2<\frac{\pi}{2})$}\\
\hline
\end{tabular}$$

\vspace{0.3truecm}

\centerline{{\bf Table 3.}}

\newpage

$$\begin{tabular}{|c|c|}
\hline
{\scriptsize$H\curvearrowright G/K$} & {\scriptsize$X\qquad(\widetilde C)$}\\
\hline
{\scriptsize$SO(4)\times SO(4)\curvearrowright$} & 
{\scriptsize$X_{(x_1,x_2)}=(-2\cot x_1-\cot(x_1-x_2)$}\\
{\scriptsize$SO(8)/U(4)$} & 
{\scriptsize$-2\cot(x_1+x_2)+2\tan x_1,$}\\
{\scriptsize} & {\scriptsize$\cot(x_1-x_2)-2\cot x_2$}\\
{\scriptsize} & {\scriptsize$-2\cot(x_1+x_2)+2\tan x_2)$}\\
{\scriptsize} & 
{\scriptsize$(\widetilde C\,:\,x_1>0,x_2>0,x_2<x_1,x_1+x_2<\pi)$}\\
\hline
{\scriptsize$\rho_3(SO(4)\times SO(4))\curvearrowright$} & 
{\scriptsize$X_{(x_1,x_2)}=(-2\cot x_1+2\tan x_1$}\\
{\scriptsize$SO(8)/U(4)$} & 
{\scriptsize$+\tan(x_1-x_2)+\tan(x_1+x_2),$}\\
{\scriptsize} & {\scriptsize$-2\cot x_2-\tan(x_1-x_2)$}\\
{\scriptsize} & {\scriptsize$+2\tan x_2+\tan(x_1+x_2))$}\\
{\scriptsize} & 
{\scriptsize$(\widetilde C\,:\,x_1>0,x_2<0,x_1+x_2<\frac{\pi}{2})$}\\
\hline
{\scriptsize$\rho_4(U(4))\curvearrowright SO(8)/U(4)$} & 
{\scriptsize$X_{(x_1,x_2)}=(-\cot x_1+3\tan x_1$}\\
{\scriptsize} & 
{\scriptsize$+\tan(x_1-x_2)+\tan(x_1+x_2),$}\\
{\scriptsize} & {\scriptsize$-\cot x_2-\tan(x_1-x_2)$}\\
{\scriptsize} & {\scriptsize$+3\tan x_2+\tan(x_1+x_2))$}\\
{\scriptsize} & 
{\scriptsize$(\widetilde C\,:\,x_1>0,x_2>0,x_1+x_2<\frac{\pi}{2})$}\\
\hline
{\scriptsize$SO(4)\times SO(6)\curvearrowright$} & 
{\scriptsize$X_{(x_1,x_2)}=2(-\cot x_1-\cot(x_1-x_2)-\cot(x_1+x_2)$}\\
{\scriptsize$SO(10)/U(5)$} & 
{\scriptsize$-\cot2x_1+\tan x_1$}\\
{\scriptsize} & {\scriptsize$+\tan(x_1-x_2)+\tan(x_1+x_2),$}\\
{\scriptsize} & {\scriptsize$\cot(x_1-x_2)-\cot x_2-\cot(x_1+x_2)$}\\
{\scriptsize} & {\scriptsize$-\cot2x_2-\tan(x_1-x_2)$}\\
{\scriptsize} & {\scriptsize$+\tan x_2+\tan(x_1+x_2))$}\\
{\scriptsize} & 
{\scriptsize$(\widetilde C\,:\,x_1>0,x_2>x_1,x_1+x_2<\frac{\pi}{2})$}\\
\hline
{\scriptsize$SO(5)\times SO(5)\curvearrowright$} & 
{\scriptsize$X_{(x_1,x_2)}=2(-\cot x_1-\cot(x_1-x_2)-\cot(x_1+x_2)$}\\
{\scriptsize$SO(10)/U(5)$} & 
{\scriptsize$+\tan x_1+\tan(x_1-x_2)$}\\
{\scriptsize} & {\scriptsize$+\tan(x_1+x_2)+\tan2x_1,$}\\
{\scriptsize} & {\scriptsize$\cot(x_1-x_2)-\cot x_2-\cot(x_1+x_2)$}\\
{\scriptsize} & {\scriptsize$-\tan(x_1-x_2)+\tan x_2$}\\
{\scriptsize} & {\scriptsize$+\tan(x_1+x_2)+\tan2x_2)$}\\
{\scriptsize} & 
{\scriptsize$(\widetilde C\,:\,x_1>0,x_2>0,x_2>x_1,x_2<\frac{\pi}{4})$}\\
\hline
{\scriptsize$\rho_5(U(5))\curvearrowright SO(10)/U(5)$} & 
{\scriptsize$X_{(x_1,x_2)}=2(-2\cot x_1-\cot2x_1$}\\
{\scriptsize} & 
{\scriptsize$+2\tan(x_1-x_2)+2\tan(x_1+x_2),$}\\
{\scriptsize} & {\scriptsize$-\cot2 x_2-2\tan(x_1-x_2)$}\\
{\scriptsize} & {\scriptsize$+2\tan(x_1+x_2)+2\tan x_2)$}\\
{\scriptsize} & 
{\scriptsize$(\widetilde C\,:\,x_1>0,x_2>0,x_1+x_2<\frac{\pi}{2})$}\\
\hline
\end{tabular}$$

\vspace{0.3truecm}

\centerline{{\bf Table 3(continued).}}

\newpage

$$\begin{tabular}{|c|c|}
\hline
{\scriptsize$H\curvearrowright G/K$} & {\scriptsize$X\qquad(\widetilde C)$}\\
\hline
{\scriptsize$SO(2)^2\times SO(3)^2\curvearrowright$} & 
{\scriptsize$X_{(x_1,x_2)}=(-\cot x_1-\cot(x_1-x_2)-\cot(x_1+x_2)$}\\
{\scriptsize$(SO(5)\times SO(5))/SO(5)$} & 
{\scriptsize$+\tan x_1+\tan(x_1-x_2)+\tan(x_1+x_2),$}\\
{\scriptsize} & {\scriptsize$\cot(x_1-x_2)-\cot x_2-\cot(x_1+x_2)$}\\
{\scriptsize} & {\scriptsize$-\tan(x_1-x_2)+\tan x_2+\tan(x_1+x_2))$}\\
{\scriptsize} & 
{\scriptsize$(\widetilde C\,:\,x_1>0,x_2>0,x_2>x_1,x_1+x_2<\frac{\pi}{2})$}\\
\hline
{\scriptsize$\rho_6(SO(5))\curvearrowright$} & 
{\scriptsize$X_{(x_1,x_2)}=2(-\cot x_1+\tan(x_1-x_2)+\tan(x_1+x_2),$}\\
{\scriptsize$(SO(5)\times SO(5))/SO(5)$} & 
{\scriptsize$-\tan(x_1-x_2)+\tan x_2+\tan(x_1+x_2))$}\\
{\scriptsize} & 
{\scriptsize$(\widetilde C\,:\,x_1>0,x_2>x_1-\frac{\pi}{2},x_1+x_2<
\frac{\pi}{2})$}\\
\hline
{\scriptsize$\rho_7(U(2))\curvearrowright Sp(2)/U(2)$} & 
{\scriptsize$X_{(x_1,x_2)}=(-\cot x_1+\tan(x_1-x_2)+\tan(x_1+x_2),$}\\
{\scriptsize} & {\scriptsize$-\cot x_2-\tan(x_1-x_2)+\tan(x_1+x_2))$}\\
{\scriptsize} & 
{\scriptsize$(\widetilde C\,:\,x_1>0,x_2>0,x_1+x_2<\frac{\pi}{2})$}\\
\hline
{\scriptsize$U(q+2)\curvearrowright$} & 
{\scriptsize$X_{(x_1,x_2)}=(-(2q-4)\cot x_1-2\cot(x_1-x_2)
-2\cot(x_1+x_2)$}\\
{\scriptsize$Sp(q+2)/Sp(2)\times Sp(q)$} & 
{\scriptsize$-2\cot2x_1+(2q-4)\tan x_1+2\tan(x_1-x_2)$}\\
{\scriptsize$(q>2)$} & 
{\scriptsize$+2\tan(x_1+x_2)+4\tan 2x_1,$}\\
{\scriptsize} & 
{\scriptsize$2\cot(x_1-x_2)-(2q-4)\cot x_2-2\cot(x_1+x_2)$}\\
{\scriptsize} & 
{\scriptsize$-2\cot2x_2-2\tan(x_1-x_2)+(2q-4)\tan x_2$}\\
{\scriptsize} & 
{\scriptsize$+2\tan(x_1+x_2)+4\tan2x_2)$}\\
{\scriptsize} & 
{\scriptsize$(\widetilde C\,:\,x_1>0,x_2>x_1,x_2<\frac{\pi}{4})$}\\
\hline
{\scriptsize$SU(4)\curvearrowright$} & 
{\scriptsize$X_{(x_1,x_2)}=(-2\cot x_1-\cot(x_1-x_2)-\cot(x_1+x_2)$}\\
{\scriptsize$Sp(4)/Sp(2)\times Sp(2)$} & 
{\scriptsize$2\tan x_1+2\tan(x_1-x_2)+3\tan(x_1+x_2),$}\\
{\scriptsize} & 
{\scriptsize$\cot(x_1-x_2)-2\cot x_2-\cot(x_1+x_2)$}\\
{\scriptsize} & 
{\scriptsize$-2\tan(x_1-x_2)+\tan x_2+3\tan(x_1+x_2))$}\\
{\scriptsize} & 
{\scriptsize$(\widetilde C\,:\,x_1>0,x_2>x_1,x_1+x_2<\frac{\pi}{2})$}\\
\hline
{\scriptsize$U(4)\curvearrowright$} & 
{\scriptsize$X_{(x_1,x_2)}=(-2\cot x_1-2\cot(x_1-x_2)-2\cot(x_1+x_2)$}\\
{\scriptsize$Sp(4)/Sp(2)\times Sp(2)$} & 
{\scriptsize$2\tan x_1+\tan(x_1-x_2)+2\tan(x_1+x_2),$}\\
{\scriptsize} & 
{\scriptsize$2\cot(x_1-x_2)-2\cot x_2-2\cot(x_1+x_2)$}\\
{\scriptsize} & 
{\scriptsize$-\tan(x_1-x_2)+\tan x_2+2\tan(x_1+x_2))$}\\
{\scriptsize} & 
{\scriptsize$(\widetilde C\,:\,x_1>0,x_2>x_1,x_1+x_2<\frac{\pi}{2})$}\\
\hline
{\scriptsize$Sp(j+1)\times Sp(q-j+1)\curvearrowright$} & 
{\scriptsize$X_{(x_1,x_2)}=(-4(j-1)\cot x_1-6\cot2x_1$}\\
{\scriptsize$Sp(q+2)/Sp(2)\times Sp(q)$} & 
{\scriptsize$+4(q-j-1)\tan x_1+4\tan(x_1-x_2)+4\tan(x_1+x_2),$}\\
{\scriptsize$(q>2)$} & 
{\scriptsize$-4(q-j-1)\cot x_2-6\cot2x_2-4\tan(x_1-x_2)$}\\
{\scriptsize} & 
{\scriptsize$+4(j-1)\tan x_2+4\tan(x_1+x_2))$}\\
{\scriptsize} & 
{\scriptsize$(\widetilde C\,:\,x_1>0,x_2>0,x_1+x_2<\frac{\pi}{2})$}\\
\hline
\end{tabular}$$

\vspace{0.3truecm}

\centerline{{\bf Table 3(continued${}^2$).}}

\newpage

$$\begin{tabular}{|c|c|}
\hline
{\scriptsize$H\curvearrowright G/K$} & {\scriptsize$X\qquad(\widetilde C)$}\\
\hline
{\scriptsize$Sp(2)\times Sp(2)\curvearrowright$} & 
{\scriptsize$X_{(x_1,x_2)}=(-3\cot x_1+\tan x_1$}\\
{\scriptsize$Sp(4)/Sp(2)\times Sp(2)$} & 
{\scriptsize$+3\tan(x_1-x_2)+4\tan(x_1+x_2),$}\\
{\scriptsize} & 
{\scriptsize$-3\cot x_2-3\tan(x_1-x_2)+4\tan(x_1+x_2))$}\\
{\scriptsize} & 
{\scriptsize$(\widetilde C\,:\,x_1>0,x_2>0,x_1+x_2<\frac{\pi}{2})$}\\
\hline
{\scriptsize$SU(2)^2\cdot SO(2)^2\curvearrowright$} & 
{\scriptsize$X_{(x_1,x_2)}=(-\cot x_1-\cot(x_1-x_2)-\cot(x_1+x_2)$}\\
{\scriptsize$(Sp(2)\times Sp(2))/Sp(2)$} & 
{\scriptsize$+\tan x_1+\tan(x_1-x_2)+\tan(x_1+x_2),$}\\
{\scriptsize} & 
{\scriptsize$\cot(x_1-x_2)-\cot x_2-\cot(x_1+x_2)$}\\
{\scriptsize} & 
{\scriptsize$-\tan(x_1-x_2)+\tan x_2+\tan(x_1+x_2))$}\\
{\scriptsize} & 
{\scriptsize$(\widetilde C\,:\,x_1>0,x_2>x_1,x_1+x_2<\frac{\pi}{2})$}\\
\hline
{\scriptsize$\rho_8(Sp(2))\curvearrowright$} & 
{\scriptsize$X_{(x_1,x_2)}=2(-\cot x_1+\tan(x_1-x_2)+\tan(x_1+x_2),$}\\
{\scriptsize$(Sp(2)\times Sp(2))/Sp(2)$} & 
{\scriptsize$-\cot x_2-\tan(x_1-x_2)+\tan(x_1+x_2))$}\\
{\scriptsize} & 
{\scriptsize$(\widetilde C\,:\,x_1>0,x_2>0,x_1+x_2<\frac{\pi}{2})$}\\
\hline
{\scriptsize$\rho_9(Sp(2))\curvearrowright$} & 
{\scriptsize$X_{(x_1,x_2)}=2(-\cot x_1+\tan(x_1-x_2)+\tan(x_1+x_2),$}\\
{\scriptsize$(Sp(2)\times Sp(2))/Sp(2)$} & 
{\scriptsize$-\cot x_2-\tan(x_1-x_2)+\tan(x_1+x_2))$}\\
{\scriptsize} & 
{\scriptsize$(\widetilde C\,:\,x_1>0,x_2>0,x_1+x_2<\frac{\pi}{2})$}\\
\hline
{\scriptsize$Sp(4)\curvearrowright$} & 
{\scriptsize$X_{(x_1,x_2)}=(-4\cot x_1-3\cot(x_1-x_2)-4\cot(x_1+x_2)$}\\
{\scriptsize$E_6/Spin(10)\cdot U(1)$} & 
{\scriptsize$+4\tan x_1+3\tan(x_1-x_2)+\tan(x_1+x_2),$}\\
{\scriptsize} & 
{\scriptsize$3\cot(x_1-x_2)-3\cot x_2-4\cot(x_1+x_2)$}\\
{\scriptsize} & 
{\scriptsize$-3\tan(x_1-x_2)+6\tan x_2+\tan(x_1+x_2))$}\\
{\scriptsize} & 
{\scriptsize$(\widetilde C\,:\,x_1>0,x_2>x_1,x_1+x_2<\frac{\pi}{2})$}\\
\hline
{\scriptsize$SU(6)\cdot SU(2)\curvearrowright$} & 
{\scriptsize$X_{(x_1,x_2)}=(-4\cot x_1-2\cot(x_1-x_2)-2\cot(x_1+x_2)$}\\
{\scriptsize$E_6/Spin(10)\cdot U(1)$} & 
{\scriptsize$-2\cot 2x_1+4\tan x_1$}\\
{\scriptsize} & 
{\scriptsize$+4\tan(x_1-x_2)+3\tan(x_1+x_2),$}\\
{\scriptsize} & 
{\scriptsize$2\cot(x_1-x_2)-4\cot x_2-2\cot(x_1+x_2)$}\\
{\scriptsize} & 
{\scriptsize$-2\cot 2x_2-4\tan(x_1-x_2)$}\\
{\scriptsize} & 
{\scriptsize$+5\tan x_2+3\tan(x_1+x_2))$}\\
{\scriptsize} & 
{\scriptsize$(\widetilde C\,:\,x_1>0,x_2>x_1,x_1+x_2<\frac{\pi}{2})$}\\
\hline
{\scriptsize$\rho_{10}(SU(6)\cdot SU(2))\curvearrowright$} & 
{\scriptsize$X_{(x_1,x_2)}=(-4\cot x_1-4\cot(x_1-x_2)-4\cot(x_1+x_2)$}\\
{\scriptsize$E_6/Spin(10)\cdot U(1)$} & 
{\scriptsize$-2\cot 2x_1+4\tan x_1$}\\
{\scriptsize} & 
{\scriptsize$+2\tan(x_1-x_2)+\tan(x_1+x_2),$}\\
{\scriptsize} & 
{\scriptsize$4\cot(x_1-x_2)-4\cot x_2-4\cot(x_1+x_2)$}\\
{\scriptsize} & 
{\scriptsize$-2\cot 2x_2-2\tan(x_1-x_2)+5\tan x_2$}\\
{\scriptsize} & 
{\scriptsize$+\tan(x_1+x_2))$}\\
{\scriptsize} & 
{\scriptsize$(\widetilde C\,:\,x_1>0,x_2>x_1,x_1+x_2<\frac{\pi}{2})$}\\
\hline
\end{tabular}$$

\vspace{0.3truecm}

\centerline{{\bf Table 3(continued${}^3$).}}

\newpage

$$\begin{tabular}{|c|c|}
\hline
{\scriptsize$H\curvearrowright G/K$} & {\scriptsize$X\qquad(\widetilde C)$}\\
\hline
{\scriptsize$\rho_{11}(Spin(10)\cdot U(1))\curvearrowright$} & 
{\scriptsize$X_{(x_1,x_2)}=(-8\cot x_1-2\cot2x_1$}\\
{\scriptsize$E_6/Spin(10)\cdot U(1)$} & 
{\scriptsize$+6\tan(x_1-x_2)+5\tan(x_1+x_2),$}\\
{\scriptsize} & 
{\scriptsize$-2\cot2x_2-6\tan(x_1-x_2)$}\\
{\scriptsize} & 
{\scriptsize$+9\tan x_2+5\tan(x_1+x_2))$}\\
{\scriptsize} & 
{\scriptsize$(\widetilde C\,:\,x_1>0,x_2>0,x_1+x_2<\frac{\pi}{2})$}\\
\hline
{\scriptsize$\rho_{12}(Spin(10)\cdot U(1))\curvearrowright$} & 
{\scriptsize$X_{(x_1,x_2)}=(-6\cot x_1-\cot(x_1-x_2)-\cot(x_1+x_2)$}\\
{\scriptsize$E_6/Spin(10)\cdot U(1)$} & 
{\scriptsize$2\tan x_1+5\tan(x_1-x_2)$}\\
{\scriptsize} & 
{\scriptsize$+4\tan(x_1+x_2)+2\tan2x_1,$}\\
{\scriptsize} & 
{\scriptsize$\cot(x_1-x_2)-6\cot x_2-\cot(x_1+x_2)$}\\
{\scriptsize} & 
{\scriptsize$-5\tan(x_1-x_2)+3\tan x_2$}\\
{\scriptsize} & 
{\scriptsize$+4\tan(x_1+x_2)+2\tan2x_2)$}\\
{\scriptsize} & 
{\scriptsize$(\widetilde C\,:\,x_1>0,x_2<\frac{\pi}{4},x_2>x_1)$}\\
\hline
{\scriptsize$Sp(4)\curvearrowright E_6/F_4$} & 
{\scriptsize$X_{(x_1,x_2)}=(-8\cot2x_1-4\cot(x_1-\sqrt3x_2)
-4\cot(x_1+\sqrt3x_2)$}\\
{\scriptsize} & 
{\scriptsize$+8\tan2x_1+4\tan(x_1-\sqrt3x_2)+4\tan(x_1+\sqrt3x_2),$}\\
{\scriptsize} & 
{\scriptsize$4\sqrt3\cot(x_1-\sqrt3x_2)-4\sqrt3\cot(x_1+\sqrt3x_2)$}\\
{\scriptsize} & 
{\scriptsize$-4\sqrt3\tan(x_1-\sqrt3x_2)+4\sqrt3\tan(x_1+\sqrt3x_2))$}\\
{\scriptsize} & 
{\scriptsize$(\widetilde C\,:\,0<x_1<\frac{\pi}{4},\frac{x_1}{\sqrt 3}<x_2<
\frac{x_1}{\sqrt 3}+\frac{\pi}{2\sqrt3},-\frac{x_1}{\sqrt 3}<x_2<
-\frac{x_1}{\sqrt 3}+\frac{\pi}{2\sqrt3})$}\\
\hline
{\scriptsize$\rho_{13}(F_4)\curvearrowright E_6/F_4$} & 
{\scriptsize$X_{(x_1,x_2)}=(-16\cot2x_1+8\tan(x_1-\sqrt3x_2)
+8\tan(x_1+\sqrt3x_2),$}\\
{\scriptsize} & 
{\scriptsize$-8\sqrt3\tan(x_1-\sqrt3x_2)
+8\sqrt3\tan(x_1+\sqrt3x_2))$}\\
{\scriptsize} & 
{\scriptsize$(\widetilde C\,:\,x_1>0,x_1-\sqrt3x_2<\frac{\pi}{2},
x_1+\sqrt3x_2<\frac{\pi}{2})$}\\
\hline
{\scriptsize$\rho_{14}(SO(4))\curvearrowright G_2/SO(4)$} & 
{\scriptsize$X_{(x_1,x_2)}=(-2\cot2x_1+3\tan(3x_1-\sqrt3x_2)
+\tan(x_1-\sqrt3x_2)$}\\
{\scriptsize} & 
{\scriptsize$+\tan(x_1+\sqrt3x_2)+3\tan(3x_1+\sqrt3x_2),$}\\
{\scriptsize} & 
{\scriptsize$-2\sqrt3\cot2\sqrt3x_2-\sqrt3\tan(3x_1-\sqrt3x_2)
-\sqrt3\tan(x_1-\sqrt3x_2)$}\\
{\scriptsize} & 
{\scriptsize$+\sqrt3\tan(x_1+\sqrt3x_2)+\sqrt3\tan(3x_1+\sqrt3x_2))$}\\
{\scriptsize} & 
{\scriptsize$(\widetilde C\,:\,x_1>0,x_2>0,\sqrt3x_1+x_2<
\frac{\pi}{2\sqrt3})$}\\
\hline
{\scriptsize$\rho_{15}(SO(4))\curvearrowright G_2/SO(4)$} & 
{\scriptsize$X_{(x_1,x_2)}=(-2\cot2x_1+3\tan(3x_1-\sqrt3x_2)
+\tan(x_1-\sqrt3x_2)$}\\
{\scriptsize} & 
{\scriptsize$+\tan(x_1+\sqrt3x_2)+3\tan(3x_1+\sqrt3x_2),$}\\
{\scriptsize} & 
{\scriptsize$-2\sqrt3\cot2\sqrt3x_2-\sqrt3\tan(3x_1-\sqrt3x_2)
-\sqrt3\tan(x_1-\sqrt3x_2)$}\\
{\scriptsize} & 
{\scriptsize$+\sqrt3\tan(x_1+\sqrt3x_2)+\sqrt3\tan(3x_1+\sqrt3x_2))$}\\
{\scriptsize} & 
{\scriptsize$(\widetilde C\,:\,x_1>0,x_2>0,\sqrt3x_1+x_2<
\frac{\pi}{2\sqrt3})$}\\
\hline
\end{tabular}$$

\vspace{0.3truecm}

\centerline{{\bf Table 3(continued${}^4$).}}

\newpage

$$\begin{tabular}{|c|c|}
\hline
{\scriptsize$H\curvearrowright G/K$} & {\scriptsize$X\qquad(\widetilde C)$}\\
\hline
{\scriptsize$\rho_{16}(G_2)\curvearrowright(G_2\times G_2)/G_2$} & 
{\scriptsize$X_{(x_1,x_2)}=2(-2\cot2x_1+3\tan(3x_1-\sqrt3x_2)
+\tan(x_1-\sqrt3x_2)$}\\
{\scriptsize} & 
{\scriptsize$+\tan(x_1+\sqrt3x_2)+3\tan(3x_1+\sqrt3x_2),$}\\
{\scriptsize} & 
{\scriptsize$-2\sqrt3\cot2\sqrt3x_2-\sqrt3\tan(3x_1-\sqrt3x_2)
-\sqrt3\tan(x_1-\sqrt3x_2)$}\\
{\scriptsize} & 
{\scriptsize$+\sqrt3\tan(x_1+\sqrt3x_2)+\sqrt3\tan(3x_1+\sqrt3x_2))$}\\
{\scriptsize} & 
{\scriptsize$(\widetilde C\,:\,x_1>0,x_2>0,\sqrt3x_1+x_2<
\frac{\pi}{2\sqrt3})$}\\
\hline
{\scriptsize$SU(2)^4\curvearrowright(G_2\times G_2)/G_2$} & 
{\scriptsize$X_{(x_1,x_2)}=(-2\cot2x_1-3\cot(3x_1-\sqrt3x_2)
-\cot(x_1-\sqrt3x_2)$}\\
{\scriptsize} & 
{\scriptsize$-\cot(x_1+\sqrt3x_2)-3\cot(3x_1+\sqrt3x_2)+2\tan2x_1$}\\
{\scriptsize} & 
{\scriptsize$+3\tan(3x_1-\sqrt3x_2)+\tan(x_1-\sqrt3x_2)$}\\
{\scriptsize} & 
{\scriptsize$+\tan(x_1+\sqrt3x_2)+3\tan(3x_1+\sqrt3x_2),$}\\
{\scriptsize} & 
{\scriptsize$\sqrt3\cot(3x_1-\sqrt3x_2)+\sqrt3\cot(x_1-\sqrt3x_2)
-\sqrt3\cot(x_1+\sqrt3x_2)$}\\
{\scriptsize} & 
{\scriptsize$-\sqrt3\cot(3x_1+\sqrt3x_2)-2\sqrt3\cot2\sqrt3 x_2
-\sqrt3\tan(3x_1-\sqrt3x_2)$}\\
{\scriptsize} & 
{\scriptsize$-\sqrt3\tan(x_1-\sqrt3x_2)+\sqrt3\tan(x_1+\sqrt3x_2)$}\\
{\scriptsize} & 
{\scriptsize$+\sqrt3\tan(3x_1+\sqrt3x_2)+2\sqrt3\tan2\sqrt3x_2)$}\\
{\scriptsize} & 
{\scriptsize$(\widetilde C\,:\,x_1>0,x_2<\frac{\pi}{4\sqrt3},x_2>\sqrt3x_1)$}\\
\hline
\end{tabular}$$

\vspace{0.3truecm}

\centerline{{\bf Table 3(continued${}^5$).}}

\vspace{1truecm}

\centerline{{\bf References}}

\vspace{0.5truecm}

{\small 
\noindent
[BV] J. Berndt and L. Vanhecke, 
Curvature-adapted submanifolds, 
Nihonkai Math. J. 

{\bf 3} (1992) 177-185.

\noindent
[BCO] J. Berndt, S. Console and C. Olmos, Submanifolds and holonomy, 
Research Notes 

in Mathematics 434, CHAPMAN $\&$ HALL/CRC Press, Boca Raton, London, New York 

Washington, 2003.

\noindent
[B] N. Bourbaki, 
Groupes et alg$\acute e$bres de Lie, Chapitre 4,5 et 6, Hermann, 
1968.

\noindent
[Ch] U. Christ, 
Homogeneity of equifocal submanifolds, J. Differential Geometry {\bf 62} 

(2002) 1-15.

\noindent
[Co] L. Conlon, 
Remarks on commuting involutions, 
Proc. Amer. Math. Soc. {\bf 22} (1969) 

255-257.

\noindent
[GT] O. Goertsches and G. Thorbergsson, 
On the Geometry of the orbits of Hermann actions, 

Geom. Dedicata {\bf 129} (2007) 101-118.

\noindent
[He] S. Helgason, 
Differential geometry, Lie groups and symmetric spaces, Academic Press, 

New York, 1978.

\noindent
[Ha] R. S. Hamilton, Three-manifolds with 
positive Ricci curvature, J. Differential Geom. 

{\bf 17} (1982) 255-306.

\noindent
[Hu1] G. Huisken, Flow by mean curvature of 
convex surfaces into spheres, J. Differential 

Geom. {\bf 20} (1984) 237-266.

\noindent
[Hu2] G. Huisken, Contracting convex hypersurfaces 
in Riemannian manifolds by their mean 

curvature, Invent. math. {\bf 84} (1986) 463-480.

\noindent
[HLO] E. Heintze, X. Liu and C. Olmos, Isoparametric submanifolds and a 
Chevalley-

type restriction theorem, Integrable systems, geometry, and topology, 151-190, 

AMS/IP Stud. Adv. Math. 36, Amer. Math. Soc., Providence, RI, 2006.
%%arXiv:math.DG/0004028.

\noindent
[HOT] E. Heintze, C. Olmos and G. Thorbergsson, 
Submanifolds with constant principal 

curvatures and normal holonomy groups, 
Intern. J. Math. {\bf 2} (1991) 167-175.

\noindent
[HPTT] E. Heintze, R.S. Palais, C.L. Terng and G. Thorbergsson, 
Hyperpolar actions 

on symmetric spaces, Geometry, topology and physics for Raoul Bott 
(ed. S. T. Yau), 

Conf. Proc. Lecture Notes Geom. Topology {\bf 4}, Internat. Press, Cambridge, 
MA, 1995 

pp214-245.

\noindent
[HTST] D. Hirohashi, H. Tasaki, H. Song and R. Takagi, Minimal orbits 
of the isotropy 

groups of symmetric spaces of compact type, Differential Goem. and 
its Appl. {\bf 13} 

(2000) 167-177.

\noindent
[Koi1] N. Koike, On proper Fredholm submanifolds in a Hilbert space arising 
from subma-

nifolds in a symmetric space, Japan. J. Math. {\bf 28} (2002) 61-80.

\noindent
[Koi2] N. Koike, Actions of Hermann type and proper complex equifocal 
submanifolds, 

Osaka J. Math. {\bf 42} (2005) 599-611.

\noindent
[Kol] A. Kollross, A Classification of hyperpolar and cohomogeneity one 
actions, Trans. 

Amer. Math. Soc. {\bf 354} (2001) 571-612.

\noindent
[LT] X. Liu and C. L. Terng, The mean curvature flow for isoparametric 
submanifolds, 

arXiv:math.DG/0706.3550v1.

\noindent
[P] R.S. Palais, 
Morse theory on Hilbert manifolds, Topology {\bf 2} (1963) 299-340.

\noindent
[PT] R.S. Palais and C.L. Terng, Critical point theory and submanifold 
geometry, Lecture 

Notes in Math. {\bf 1353}, Springer, Berlin, 1988.

\noindent
[S] G. Schwarz, 
Smooth functions invariant under the action of a compact Lie group, 

Topology {\bf 14} (1975) 63-68.

\noindent
[T1] C.L. Terng, 
Isoparametric submanifolds and their Coxeter groups, 
J. Differential Ge-

ometry {\bf 21} (1985) 79-107.

\noindent
[T2] C.L. Terng, 
Proper Fredholm submanifolds of Hilbert space, 
J. Differential Geometry 

{\bf 29} (1989) 9-47.

\noindent
[TT] C.L. Terng and G. Thorbergsson, 
Submanifold geometry in symmetric spaces, J. Diff-

erential Geometry {\bf 42} (1995) 665-718.

\noindent
[Z] X. P. Zhu, Lectures on mean curvature flows, Studies in Advanced Math., 
AMS/IP, 2002.

\vspace{0.5truecm}

{\small 
\rightline{Department of Mathematics, Faculty of Science}
\rightline{Tokyo University of Science, 1-3 Kagurazaka}
\rightline{Shinjuku-ku, Tokyo 162-8601 Japan}
\rightline{(koike@ma.kagu.tus.ac.jp)}
}

\end{document}